\documentclass[12pt]{amsart}
\usepackage{epsf,overpic,amssymb,amscd,times,stmaryrd,euscript}
\usepackage{tikz-cd} 
\usepackage[utf8]{inputenc}
\usepackage[english]{babel}
\usepackage{amsmath}
\usepackage{amsthm}
\usepackage{lipsum}
\usepackage{amsfonts}
\usepackage{amssymb}
\usepackage{amscd}
\usepackage{bezier}
\usepackage{enumerate}
\usepackage{fancyhdr}
\usepackage{graphicx}
\usepackage{epstopdf}
\usepackage{booktabs}
\usepackage{rotating}
\usepackage{listings}
\usepackage[square, numbers, comma, sort&compress]{natbib} % Use the natbib reference package - read up on this to edit the reference style; if you want text (e.g. Smith et a
\usepackage{fancyhdr}
\usepackage{amsmath,amsfonts,amssymb,amscd,amsthm,xspace}
\usepackage[centerlast,small,sc]{caption}
\usepackage[pdfpagemode={UseOutlines},bookmarks=true,bookmarksopen=true,
   bookmarksopenlevel=0,bookmarksnumbered=true,hypertexnames=false,
   colorlinks,linkcolor={blue},citecolor={blue},urlcolor={red},
   pdfstartview={FitV},unicode,breaklinks=true]{hyperref}
\usepackage{upgreek}
\usepackage[margin=3cm]{geometry}
\usepackage{lipsum}
\usepackage[page]{appendix}

\newcommand{\Vol}{\operatorname{Vol}} %Volume

\newtheorem{lemma}{Lemma}
\newtheorem{prop}{Proposition}
\newtheorem{thrm}{Theorem}
\newtheorem{coro}{Corollary}
\newtheorem*{theorem*}{Theorem}
\newtheorem*{corollary*}{Corollary}
\newtheorem{defi}{Definition}
\theoremstyle{definition}

\newtheorem{remark}{Remark}[section]

\begin{document}

\title{ Legendrian Contact Homology and Topological Entropy }

\author{Marcelo R.R. Alves}
\address{Institut de Math\'ematiques,
Universit\'e de Neuch\^atel,
Rue \'Emile Argand 11,
CP~158,
2000 Neuch\^atel,
Switzerland}
\email{marcelo.ribeiro@unine.ch} \urladdr{}

\date{}



\begin{abstract}
In this paper we study the growth rate of a version of Legendrian contact homology, which we call strip Legendrian contact homology, in 3-dimensional contact manifolds and its relation to the topological entropy of Reeb flows. We show that: if for a pair of Legendrian knots in a contact 3-manifold $(M,\xi)$ the strip Legendrian contact homology is defined and has exponential homotopical growth with respect to the action, then every Reeb flow on $(M,\xi)$ has positive topological entropy. This has the following dynamical consequence: for all Reeb flows (even degenerate ones) on $(M,\xi)$ the number of hyperbolic periodic orbits grows exponentially with respect to the period. We show that for an infinite family of 3-manifolds, infinitely many different contact structures exist that possess a pair of Legendrian knots for which the strip Legendrian contact homology has exponential growth rate.
\end{abstract}

\maketitle

\section{Introduction}

The objective of this article is to study the growth rate of a version of Legendrian contact homology and its implications to the dynamics of Reeb flows on contact 3-manifolds. It fits in a larger program which aims at understanding the relationship between SFT-invariants of a contact structure and global dynamical invariants of Reeb flows.
In this work we focus our attention on one dynamical invariant, the topological entropy. The topological entropy $h_{top}$ is a non-negative number that one associates to a dynamical system and which measures the complexity of the system.  Positivity of the topological entropy for a dynamical system implies that the system has some type of exponential instability. It is a deep result of Yomdin (see \cite{PAT}) that $h_{top}(\phi)$ for a $C^{\infty}$-flow $\phi = (\phi^t)_{t\in\mathbb{R}}$ on a compact manifold satisfies 
\begin{equation}
h_{top}(\phi) \geq v(\phi)
\end{equation}
where
\[
v(\phi) = \sup_{N\subset M} v(\phi,N), \text{ and }
\]
\[
v(\phi, N) = \limsup_{t\to \infty} \frac{\log \Vol_g^n(\phi^t(N))}{t}.
\]
Here, $n = \dim N$, the  supremum is taken over all submanifolds $N \subset M$, and $\Vol^n_g$ is the $n$-dimensional volume with respect to some Riemannian metric $g$ on $M$.

\subsection{Basic definitions and main results}

We first recall some basic definitions from contact geometry. A 1-form $\lambda$  on a $(2n+1)$-dimensional manifold $Y$ called a \textit{contact form} if $ \lambda \wedge (d\lambda)^n $ is a volume form on $Y$. The hyperplane $\xi= \ker \lambda$ is called the \textit{contact structure}. For us a \textit{contact manifold} will be a pair $(Y,\xi)$ such that $\xi$ is the kernel of some contact form $\lambda$ on $Y$ (these are usually called co-oriented contact manifolds in the literature). When $\lambda$ satisfies $\xi = \ker\lambda$, we will say that $\lambda$ is a contact form on $(Y,\xi)$. On any contact manifold there always exist infinitely many different contact forms. Given a contact form $\lambda$, its
\textit{Reeb vector field} is the unique vector field $X_\lambda$ satisfying $\lambda(X_\lambda)=1$ and $i_{X_\lambda}d\lambda=0$. The Reeb flow $\phi_{X_{\lambda}}$ of $\lambda$ is the flow generated by the vector field $X_\lambda$. We will refer to the periodic orbits of $\phi_{X_{\lambda}}$ as \textit{Reeb orbits} of $\lambda$. The action $A(\gamma)$ of a Reeb orbit is defined by $A(\gamma):=\int_{\gamma} \lambda$. A contact form is called \textit{hypertight} if its Reeb flow has no contractible periodic orbits.

An isotropic submanifold of $(Y,\xi)$ is a submanifold of $Y$ whose tangent bundle is contained in $\xi$. When this submanifold is of maximal possible dimension, it is called a \textit{Legendrian submanifold} of $(Y,\xi)$. It turns out that the maximal possible dimension of an isotropic submanifold is $n$. Given a contact form $\alpha$ and a pair of Legendrian submanifolds $(\Lambda,\widehat{\Lambda})$, a \textit{Reeb chord} of $\alpha$ from $\Lambda$ to $\widehat{\Lambda}$ is a trajectory $c$ of the Reeb flow of $\alpha$ that starts in $\Lambda$ and ends in $\widehat{\Lambda}$. We define the action $A(c)$ of a Reeb chord $c$ as $A(c)=\int_c \alpha$. A Reeb chord $c$ is said to be \textit{transverse} if the intersection $\phi^{A(c)}_{X_\alpha}(\Lambda) \cap \widehat{\Lambda}$ is transverse at the endpoint of $c$. If all Reeb chords of $\alpha$ from $\Lambda$ to $\widehat{\Lambda}$ are transverse and do not intersect Reeb orbits of $\alpha$ we say that $(\alpha,\Lambda \to \widehat{\Lambda})$ is \textit{regular}.

We study the topological entropy of Reeb flows from the point of view of contact topology. More precisely, we search for conditions on the topology of a contact manifold $(M,\xi)$ that force \textbf{all} Reeb flows on $(M,\xi)$ to have positive topological entropy. The condition we impose is on the behaviour of a contact topological invariant called strip Legendrian contact homology. We show that if a contact manifold $(M,\xi)$ admits a contact form $\lambda_0$ and a pair of Legendrian knots $(\Lambda,\widehat{\Lambda})$ for which the strip Legendrian contact homology is well defined and has \textit{exponential homotopical growth}, then all Reeb flows on  $(M,\xi)$ have positive topological entropy.

The inspiration for studying the topological entropy of Reeb flows from this perspective comes from the study of geodesic flows.
It is an important result in the theory of geodesic flows that if the based loop space of a manifold $Q$ has a rich homology, then for every Riemannian metric $g$ on $Q$ the geodesic flow of $g$ on the unit tangent bundle $T_1 Q $ of $Q$ has positive topological entropy; see Paternain's book \cite{PAT}. Geodesic flows on the unit tangent bundle are examples of Reeb flows. The reason for this is that there exists a contact structure $\xi_{geo}$ on $T_1 Q$ such that for every Riemannian metric $g$ on $Q$, the geodesic flow $\phi_g$ of $g$ coincides with the Reeb flow of a contact form $\lambda_g$ on $(T_1 Q, \xi_{geo})$. In \cite{MS} Macarini and Schlenk showed that if the based loop space of $Q$ has rich homology then every Reeb flow on $(T_1 Q, \xi_{geo})$ has positive topological entropy.
 Their result can be seen as an extension of the result previously known for geodesic flows in $T_1 Q$ to the larger class of Reeb flows on $ (T_1 X, \xi_{geo} )$. The strategy to estimate the topological entropy used in \cite{MS} can be sketched as follows:
 \begin{center}
Exponential growth of Lagrangian Floer homology of the cotangent fiber $(T^*Q)_{|_p}$ \\
$\Rightarrow$ \\
Exponential volume growth of the unit tangent fiber $(T_1 Q)_{|_p}$ for all Reeb flows in $(T_1 Q,\xi_{geo})$ \\
$\Rightarrow$ \\
Positivity of the topological entropy for all Reeb flows in $(T_1 Q,\xi_{geo})$.
\end{center}
To obtain the first implication Macarini and Schlenk use the fact that the contact manifold $(T_1 Q,\xi_{geo})$ has the structure of a Legendrian fibration, and apply the geometric idea of \cite{FS1,FS2} to show that the number of trajectories connecting a  Legendrian fiber to  another Legendrian fiber can be used to obtain a volume growth estimate.

Inspired by \cite{MS} we present in this article techniques that can be used to estimate the topological entropy for Reeb flows on contact 3-manifolds which are not unit cotangent bundles. Because most contact $3$-manifolds are not Legendrian fibrations, it is usually not possible to apply the scheme presented above. However a sufficiently small neighbourhood of a given Legendrian knot in a contact $3$-manifold can always be given the structure of a Legendrian fibration. This local structure of a Legendrian fibration allows us to use a localised version of the geometric idea of \cite{FS1,FS2}, and to prove that the exponential growth of the strip Legendrian contact homology of a pair of Legendrian knots on a contact 3-manifold $(M,\xi)$, implies exponential growth of the length of these Legendrian knots for any Reeb flow on $(M^3,\xi)$. Applying Yomdin's Theorem (see \cite{PAT}) we obtain the following result, which is the main structural result of this paper.
\begin{thrm} \label{theorem1}
Let $(Y,\xi )$ be a contact 3-manifold, and $\lambda_0$ be a hypertight contact form on $(Y,\xi)$. Assume that $\lambda_0$ is adapted to the pair of disjoint Legendrian knots $(\Lambda,\widehat{\Lambda})$.  Then, if the strip Legendrian contact homology $LC\mathbb{H}_{st}(\lambda_0,\Lambda \to\widehat{\Lambda})$ has exponential homotopical growth rate with exponential weight $a>0$, it follows that the Reeb flow of any $C^\infty$ contact form $\lambda$ on $(Y,\xi)$ has positive topological entropy.
Moreover, if we denote by $f_\lambda$ the positive function such that $\lambda=f_\lambda \lambda_0$, we have
\begin{equation}
h_{top}(\phi_{X_{\lambda}}) \geq \frac{a}{\max(f_\lambda)}.
\end{equation}
\end{thrm}
We refer the reader to Definition \ref{defi:adapted} for the definition of contact forms \textit{adapted} to a pair of disjoint Legendrian knots, and to Definition \ref{defigrowth}  to the definition of the exponential homotopical growth of the strip Legendrian contact homology. The fact that $\lambda_0$ is adapted to the pair of disjoint Legendrian knots $(\Lambda,\widehat{\Lambda})$ implies that $(\lambda_0, \Lambda \to \widehat{\Lambda})$ is regular. Notice however that we do not require  $\lambda_0$ to be non-degenerate. The reason for this is that the definition of $LC\mathbb{H}_{st}$ only uses pseudoholomorphic strips which detect $\lambda_0$-Reeb chords from $\Lambda$ to $\widehat{\Lambda}$. Since these strips do not ``see'' Reeb orbits of $\lambda_0$ their non-degeneracy is irrelevant for the construction of $LC\mathbb{H}_{st}$; see a more detailed discussion in Remark \ref{rmk:nondegeneracy}.

Theorem \ref{theorem1} gives a criterion implying that all Reeb flows on a contact 3-manifold have positive topological entropy. To explain the significance of this result we quote a result due to Katok \cite{K,K2} and Lima and Sarig \cite{LS,S} which shows that positivity of topological entropy of the 3-dimensional flow of a vector field without singularities has deep implications for the dynamics of such flows.
\begin{theorem*}
If $\phi$ is a smooth flow on a closed oriented 3-manifold generated by a non-vanishing vector field, then $\phi$ has positive topological entropy if, and only if, there exists a Smale ``horseshoe'' as a subsystem of the flow. As a consequence, the number of hyperbolic periodic orbits of $\phi$ grows exponentially with respect to the period.
\end{theorem*}
Here a ``horseshoe'' is a compact invariant set where the dynamics is semi-conjugate to that of the suspension of a finite shift by a finite-to-one map. In particular, the number of hyperbolic periodic orbits on a ``horseshoe'' of a 3-dimensional flow $\phi$ grows exponentially with respect to the period. The recent work of Lima and Sarig \cite{LS} proves a stronger result: there exists an invariant set $\mathcal{K}$ of $\phi$ of full topological entropy where the dynamics of the restriction $\phi |_{\mathcal{K}}$ is semi-conjugate to the suspension of a countable shift by a finite-to-one map. It is worthy noting that their set $K$ is not necessarily a horseshoe, since it may not be compact. Nevertheless, since the topological entropy on countable shifts is the supremum of the topological entropy of its finite subshifts, \cite{LS} implies that $\phi$ has horseshoes with topological entropy as close to $h_{top}(\phi)$ as we wish.
\footnote{We remark that in \cite{K} Katok proves similar results for diffeomorphisms on surfaces, and in \cite{K2} he states the result for flows on 3-manifolds mentioned above. The complete proofs for 3-dimensional flows with positive topological entropy only appeared in \cite{LS}, which builds on the ideas of \cite{K,K2,S}.}

It follows from the result of Katok, and Lima and Sarig, that if a contact 3-manifold $(M,\xi)$ satisfies the hypothesis of Theorem \ref{theorem1}, then for every Reeb flow on $(M,\xi)$ the number of hyperbolic periodic orbits grows exponentially with the period. An interesting feature of this result (and also of a result in \cite{MS}) is that it gives a method of proving the existence of many Reeb orbits of a Reeb flow which is not based on the use of a contact topological invariant defined using Reeb orbits. The appearance of the Reeb orbits is obtained as a consequence of a multiplicity result for the existence of \textbf{Reeb chords}.

The other main result of the paper is the existence of a large family of contact 3-manifolds which satisfy the hypothesis of Theorem \ref{theorem1}.
\begin{thrm} \label{theorem2}
Let $M$ be a closed oriented connected 3-manifold which can be cut along a nonempty family of incompressible tori into a family $\{M_i, 0 \leq i \leq q\}$ of irreducible manifolds with boundary such that the component $M_0$ satisfies:
\begin{itemize}
\item{$M_0$ is the mapping torus of a diffeomorphism $h: S \to S$ with pseudo-Anosov monodromy on a surface $S$ with non-empty boundary.}
\end{itemize}
Then $M$ can be given infinitely many different tight contact structures $\xi_k$, such that:
\begin{itemize}
\item There exist disjoint Legendrian knots $\Lambda$, $\widehat{\Lambda}$ on $(M,\xi_k)$, and a contact form $\lambda_k$ on $(M,\xi_k)$ adapted to the pair $\Lambda$ and $\widehat{\Lambda}$, such that $LC\mathbb{H}_{st}(\lambda_k,\Lambda \to \widehat{\Lambda})$ has exponential homotopical growth rate. It follows that every Reeb flow on $(M,\xi_k)$ has positive topological entropy.
\end{itemize}
\end{thrm}
The contact manifolds covered by this theorem are among the examples of tight contact manifolds constructed in \cite{Colin} and coincide with the ones studied in \cite{A,Vaugon}. In particular, Theorem \ref{theorem2} implies that on these contact 3-manifolds every Reeb flow has positive topological entropy, a result which was proved by the author using different methods in \cite{A}.

It is important to mention that the methods developed in this paper and in \cite{A} are independent. We believe that they are complementary. One advantage of the methods of \cite{A} is that they give estimates for the topological entropy for $C^2$ Reeb flows while the method we use here is only valid for $C^\infty$ Reeb flows. On the other hand, there are examples of high-dimensional contact manifolds which are simply connected and on which every Reeb flow has positive topological entropy: the unit tangent bundles of simply connected manifolds presented in \cite{MS}. On simply connected manifolds there is no hope of applying the methods of \cite{A} but the volume growth technique for estimating topological entropy developed here can still be used. In general, we expect that detecting exponential volume growth of Legendrian submanifolds  should be the most robust method to estimate topological entropy of Reeb flows in high dimensional contact manifolds.
\subsection{Related developments}
The techniques developed in this article are fundamental for several projects being developed by the author with collaborators.
We list some of these projects.
\begin{itemize}
\item In a joint project with Pedro A.S. Salom\~ao \cite{A2} we combine ideas of this article with those of \cite{HMS,Mo} to study the forcing of topological entropy in contact 3-manifolds. More precisely, we establish a criterion to determine if a transverse knot $\mathcal{K}$ on a contact 3-manifold $(M,\xi)$ forces topological entropy, in the sense that every Reeb flow on $(M,\xi)$ which has $\mathcal{K}$ as a Reeb orbit has positive topological entropy. This criterion is then used to obtain a contact topological generalisation of a theorem of Denvir and Mackay \cite{DM} that says that if a Riemannian metric of the two-dimensional torus that has a simple closed contractible geodesic, then its geodesic flow has positive topological entropy.
\item The attentive reader will notice that the argument used for Theorem \ref{theorem1} remains true in high-dimensions provided that the Legendrian submanifolds being considered are spheres. The reason for this is that a Legendrian sphere $\Lambda$ always admits a neighbourhood which is foliated by other Legendrian spheres in the same isotopic class as $\Lambda$. In joint projects with Leonardo Macarini \cite{AM}, Felix Schlenk and Matthias Meiwes \cite{AlvesMeiwes} we will use this result to give examples of high-dimensional contact manifolds on which every Reeb flow has positive topological entropy, and which are not unit tangent bundles or quotients of unit tangent bundles. 
\item In a joint work with Vincent Colin and Ko Honda we combine the techniques developed in this paper with those of \cite{CH2} to use openbook decompositions to prove positivity of topological entropy for all Reeb flows on many contact 3-manifolds. This work aims at establishing that on ``most'' tight contact 3-manifolds every Reeb flow has positive topological entropy.
\item In our proof of Theorem \ref{theorem1} we have showed that if, on a contact 3-manifold $(Y,\xi)$, there exists a contact form $\lambda_0$ and a pair of disjoint Legendrian knots $\Lambda$ and $\widehat{\Lambda}$ satisfying the hypothesis of Theorem \ref{theorem1} then for any contact form $\lambda$ on $(Y,\xi)$ the length of $\phi^t_{X_\lambda}(\Lambda)$ grows exponentially as $t\to +\infty$. It is natural to ask if, under the same hypothesis, this growth property is shared by all Legendrian knots $\Lambda'$ belonging to the same Legendrian isotopy class of $\Lambda$. In a forthcoming paper \cite{A5} we will show that this is in fact the case. In this same paper we will also show that, under the hypothesis of Theorem \ref{theorem1}, we have that for any contact form $\lambda$ on $(M,\xi)$ the topological entropy of $\phi_{X_\lambda}$ is positive when restricted to $\upomega$-limit $\upomega_\lambda(\Lambda')$ of a generic Legendrian knot $\Lambda'$ in the Legendrian isotopy class of $\Lambda$; a proof of this result in the case where the contact 3-manifold $(M,\xi)$ is the unit tangent bundle of a surface of genus $\geq 2$ is given in \cite{Orsay}.
 The $\upomega$-limit\footnote{This is a natural generalisation of the notion of $\upomega$-limit set of a point, which is central in the theory of dynamical systems.} $\upomega_\lambda(\Lambda')$  of $\Lambda'$ is a compact set invariant by the the flow $\phi^t_{X_\Lambda}$ defined as:
     \begin{center}
       $\upomega_\lambda(\Lambda'):= \{p \in Y \ | \ \exists \mbox{ sequences } x_n \in \Lambda' \mbox{ and } t_n \to +\infty \mbox{ such that } \phi^{t_n}_{X_\lambda}(x_n) \to p \}.$
     \end{center}
     The positivity of the topological entropy of $\phi^t_{X_\Lambda}$ restricted to $\upomega_\lambda(\Lambda')$ can be interpreted in the following way: asymptotically, the path of curves $\phi^t_{X_\lambda}(\Lambda')$ detects the positive topological entropy of $\phi^t_{X_\lambda}$. We believe that this result helps to clarify the role that Legendrian curves can have in the study of the global dynamics of Reeb flows.
\end{itemize}

\subsection{Geometric idea of the proof of Theorem \ref{theorem1}}

We will give an intuitive idea of the proof of Theorem \ref{theorem1}. The basic idea, which is also used in \cite{MS}, is to use the number of Reeb chords of $\lambda$ from a Legendrian submanifold $\Lambda$ to other Legendrian submanifolds to estimate the growth rate of the volume of $\phi_{X_\lambda}^t(\Lambda)$.

In \cite{MS} the authors study Reeb chords from one fixed unit tangent fiber $\Lambda_0$ in the unit tangent bundle $T_1Q$ to the unit tangent fibers $\Lambda_q$, where $q$ belongs to a set of full measure in $Q$. Using Lagrangian Floer homology they show that if the homology of the based loop space of $Q$ is rich, then there exist numbers $C_0 \geq 0$, $\overline{a}(\lambda)>0$, and $\overline{d}(\lambda)$, and a set $\mathcal{U}$ of full measure in $Q$, such that for $q\in \mathcal{U}$ we have:
\begin{equation}
 N_C(\lambda,\Lambda_0,{\Lambda}_q) \geq e^{C\overline{a} + \overline{d}} \mbox{ for all } C\geq C_0,
\end{equation}
where $N_C(\lambda,\Lambda_0,{\Lambda}_q)$ is the number of Reeb chords of $\lambda$ from $\Lambda_0$ to $\Lambda_q$. Using the canonical projection from $T_1Q$ to $Q$ they estimate the area of the cylinder $Cyl^C_{X_\lambda}(\Lambda_q):=\{\phi_{X_\lambda}^t(\Lambda_0); t\in[0,C]\}$ via the counting functions $N_C(\lambda,\Lambda_0,{\Lambda}_q)$. The idea of using counting functions for such area estimates is due to Gabriel Paternain; see \cite{PAT}. The result is an inequality of the form:
\begin{equation}
Area(Cyl^C_{X_\lambda}(\Lambda_q)) \geq \int_{Q} N_C(\lambda,\Lambda_0,{\Lambda}_q) d\mu_g \geq \mu_g(Q)e^{C\overline{a} + \overline{d}}
\end{equation}
where $\mu_g$ is the measure induced by a Riemannian metric $g$ on $Q$, and $\mu_g(Q)$ is volume of $Q$ for the measure $\mu_g$. We point out, that in this last inequality, the fact that $T_1 Q$ is a Legendrian fibration is used in a crucial way.

\

Most contact 3-manifolds do not have the structure of a Legendrian fibration. In fact, it was proved by Giroux (see \cite[Proposition 1.1.7]{confoliations}) that the only contact 3-manifolds with such structures are unit tangent bundles of surfaces and their coverings. However, a sufficiently small neighbourhood of a Legendrian knot $\widehat{\Lambda}$ on a contact 3-manifold $(Y,\xi)$, is always a Legendrian fibration. This is a consequence of the Weinstein Legendrian neighbourhood Theorem, whose 3-dimensional version asserts that sufficiently small neighbourhoods of Legendrian knots are always contactomorphic; i.e there exists a normal form for small neighbourhoods of Legendrian knots.

In the hypotheses of Theorem \ref{theorem1} we have a pair of disjoint Legendrian knots $(\Lambda, \widehat{\Lambda})$ in $(Y,\xi)$ and a contact form $\lambda_0$, for which the strip Legendrian contact homology $LC\mathbb{H}_{st}(\lambda_0,\Lambda \to \widehat{\Lambda})$ has exponential homotopical growth rate. We begin by choosing a small neighbourhood $\mathcal{V}_\epsilon(\widehat{\Lambda})$ of  $\widehat{\Lambda}$ which is contactomorphic to $(S^1 \times \mathbb{D},\ker(\cos(\theta)dx +\sin(\theta) dy))$; where $(\theta,x,y) \in S^1\times \mathbb{D}$ and $\widehat{\Lambda}$ is identified with $S^1 \times \{0\}$. It is clear that for all $z:=(x,y) \in \mathbb{D}$ the curve $\widehat{\Lambda}^z:=S^1 \times \{z\}$ is Legendrian in $\mathcal{V}_\epsilon(\widehat{\Lambda})$.

Let $a>0$ be the exponential weight of the growth of $LC\mathbb{H}_{st}(\lambda_0,\Lambda \to \widehat{\Lambda})$. Using the invariance properties of the Legendrian contact homology we show that given $\delta >0$, if we choose $\mathcal{V}_\epsilon(\widehat{\Lambda})$ to be sufficiently small, then there exist numbers $C_0 \geq 0$ and $d$, and a subset $\mathcal{U}$ of $\mathbb{D}$ of full measure, such that for all $C \geq C_0$ and $z  \in \mathcal{U}$ the number $N_C(\lambda,\Lambda,\widehat{\Lambda}^z)$  of Reeb chords of $\lambda$ from $\Lambda$ to $\widehat{\Lambda}^z$ satisfies:
\begin{equation}
N_C(\lambda,\Lambda,\widehat{\Lambda}^z) \geq e^{\frac{aC}{(1+4\delta)\max{f_\lambda}}+d},
\end{equation}
where $f_\lambda$ is the function such that $\lambda=f_\lambda \lambda_0$.

\begin{figure}
\centering
\includegraphics[width=9cm]{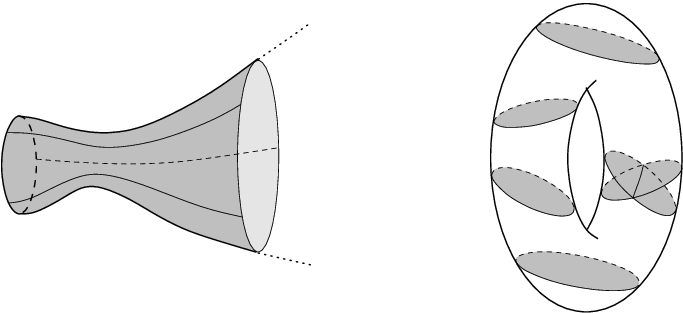}
\caption{\footnotesize{On the left side of the picture we have a piece of the immersed cylinder $Cyl^C_{X_\lambda}(\Lambda)$.
The shaded discs on the right side of the picture are pieces of the intersection $Cyl^C_{X_\lambda}(\Lambda) \cap \mathcal{V}_\epsilon(\widehat{\Lambda})$. Under the hypotheses of Theorem \ref{theorem1}, the area of this intersection grows exponentially with respect $C$.}}\label{figure}
\end{figure}

We can use the counting functions $N_C(\lambda,\Lambda,\widehat{\Lambda}^z) $ to estimate the area of the intersection $Cyl^C_{X_\lambda}(\Lambda) \cap \mathcal{V}_\epsilon(\widehat{\Lambda})$ between the cylinder $Cyl^C_{X_\lambda}(\Lambda)$ and the neighbourhood $\mathcal{V}_\epsilon(\widehat{\Lambda})$ (see figure \ref{figure}).  We first choose a metric $g_0$ on $Y$ which restricts to the ``flat'' metric $d\theta \otimes d\theta + dx \otimes dx + dy \otimes dy$ in the coodinates $(\theta,x,y)$ on $\mathcal{V}_\epsilon(\widehat{\Lambda})$.
Because $\mathcal{V}_\epsilon(\widehat{\Lambda})$ has the structure of a Legendrian fibration, we can apply a local version of Paternain's idea to obtain that for $C\geq C_0$
\begin{equation}
Area(Cyl^C_{X_\lambda}(\Lambda) \cap \mathcal{V}_\epsilon(\widehat{\Lambda})) \geq \int_{\mathbb{D}} N_C(\lambda,\Lambda,\widehat{\Lambda}^z) dx dy \geq e^{\frac{aC}{\max(f_\lambda)(1+4\delta)} + {d'}},
\end{equation}
where the constant $d'=2\uppi d$.

This local picture shows that the area $Area(Cyl^C_{X_\lambda}(\Lambda) \cap \mathcal{V}_\epsilon(\widehat{\Lambda}))$ (of the intersection between $Cyl^C_{X_\lambda}(\Lambda)$ and $\mathcal{V}_\epsilon(\widehat{\Lambda})$) already grows exponentially fast and allows us to estimate the topological entropy of the Reeb flow of $X_\lambda$ by using Yomdin's theorem.

\begin{remark}
 It is likely that the local picture we used here can also be applied to estimate the intermediate and slow entropies which were studied in \cite{FLS,FS3} in contact manifolds which are not unit tangent bundles.
\end{remark}

\subsection{Structure of the paper}
In Section \ref{section2} we define the strip Legendrian contact homology and show under which conditions it can be defined. In Section \ref{sectiongrowth} we define the exponential homotopical growth of the strip Legendrian contact homology and prove Theorem \ref{theorem1}. In Section \ref{section4} we prove Theorem \ref{theorem2}. In Section \ref{section5} we propose some open questions which we believe are interesting for future developments.

\

\textbf{Acknowledgements:} special thanks to my professors Fr\'ed\'eric Bourgeois and Chris Wendl for their guidance and support during the development of this work, which is a part of the author's Ph.D. thesis that was developed under their supervision. My special thanks to Fr\'ed\'eric for his guidance in the world of Legendrian contact homology, without which this work would not have been possible. I would also like to thank Pedro Salom\~ao for helpful discussions and for explaining me the techniques used in \cite{HMS} which made it possible to avoid dealing with multiply covered pseudoholomorphic curves. I also thank the anonymous referee for many helpful suggestions for improvement. My gratitude to Felix Schlenk for his interest in this work and for suggesting many future directions of research. My personal thanks to Ana Nechita, Andr\'e Alves, Hilda Ribeiro and Lucio Alves for their unconditional personal support. I warmly thank FNRS-Belgium for the financial support which allowed me to focus on mathematical research.

\section{The strip Legendrian contact homology} \label{section2}

\setcounter{thrm}{0}

The strip Legendrian contact homology is defined in the spirit of the SFT invariants introduced in \cite{SFT}. To define it we use pseudoholomorphic curves in symplectizations of contact manifolds. Pseudo-holomorphic curves were introduced in symplectic manifolds by Gromov in \cite{Gr} and adapted to symplectizations and symplectic cobordisms by Hofer \cite{H}; see also \cite{CPT} as a general reference for pseudoholomorphic curves in symplectic cobordisms.

\subsection{Pseudoholomorphic curves in symplectizations and symplectic cobordisms}

\subsubsection{Cylindrical almost complex structures and Lagrangian cylinders} \label{sec:symp}

Let $(Y,\xi)$ be a contact manifold and $\lambda$ a contact form on $(Y,\xi)$. The symplectization of $(Y,\xi)$ is the product $\mathbb{R} \times Y$ with the symplectic form $d(e^s \lambda)$ (where $s$ denotes the $\mathbb{R}$ coordinate in $\mathbb{R} \times Y$). The 2-form $d\lambda$ restricts to a symplectic form on the vector bundle $\xi$ and it is well known that the set  $\mathfrak{j}(\lambda)$ of $d\lambda$-compatible almost complex structures on the symplectic vector bundle $\xi$ is non-empty and contractible. Notice that if $Y$ is 3-dimensional, the set $\mathfrak{j}(\lambda)$ does not depend on the contact form $\lambda$ on $(Y,\xi)$.

For $j \in \mathfrak{j}(\lambda)$ we can define an $\mathbb{R}$-invariant almost complex structure $J$ on $\mathbb{R} \times Y$ by demanding that
\begin{equation} \label{eq16}
J \partial_s = X_\lambda, \ \ J\mid_\xi = j.
\end{equation}
We will denote by $\mathcal{J}(\lambda)$ the set of almost complex structures in $\mathbb{R} \times Y$ that are $\mathbb{R}$-invariant, $d(e^s\lambda)$-compatible and satisfy \eqref{eq16} for some $j \in \mathfrak{j}(\lambda)$.

If $\overline{\Lambda}$ is a Legendrian submanifold of $(Y,\xi)$ we denote by $Z_{\overline{\Lambda}}$ the Lagrangian cylinder $\mathbb{R}\times {\overline{\Lambda}}$ over ${\overline{\Lambda}}$ in the symplectization of $\lambda$. It is immediate to check that $Z_{\overline{\Lambda}}$ is an exact Lagrangian submanifold in the symplectization $(\mathbb{R} \times Y,d(e^s\lambda))$.

\subsubsection{Straight exact symplectic cobordisms}

An exact symplectic cobordism is, roughly, an exact symplectic manifold $(W,\varpi)$ that outside a compact subset is like the union of  cylindrical ends of symplectizations. We restrict our attention to a special case of exact symplectic cobordisms which we call \textit{straight}.

We first fix a contact manifold $(Y,\xi)$ and let $\lambda^+$ and $\lambda^-$ be contact forms on $(Y,\xi)$. It follows that there exists a function $f^{\lambda^+}_{\lambda^-}:Y \to (0,+\infty)$ such that $  f^{\lambda^+}_{\lambda^-}\lambda^- = \lambda^+$. Assume that  $f^{\lambda^+}_{\lambda^-}>1$. Then one can chose a function $h: \mathbb{R}\times Y \to (0,+\infty)$ such that:
\begin{eqnarray}
\mbox{ there exists } R^+ \in \mathbb{R} \mbox{ such that } h(s,p) = e^{s-R^+}f^{\lambda^+}_{\lambda^-}(p) \mbox{ for } s\geq R^+, \\
\mbox{ there exists } R^-< R^+ \mbox{ such that } h(s,p) = e^{s-R^-} \mbox{ for } s\leq R^-, \\
\partial_s h > 0 \mbox{ at all points.} \ \ \ \ \ \ \ \ \ \ \ \ \ \ \ \ \ \ \ \ \ \ \ \ \ \ \ \ \ \ \ \ \ \ \ 
\end{eqnarray}
It then follows that $d (h\lambda^-)$ is an exact symplectic form on $\mathbb{R}\times Y$. We call the pair $ ( \mathbb{R}\times Y, h\lambda^-)$ an \textit{straight exact symplectic cobordism} from $\lambda^+$ to $\lambda^-$.
To simplify notation we write $W= (\mathbb{R}\times Y, h\lambda^-)$ and divide $W$ in three regions:
\begin{eqnarray*}
W^+ = [R^+,+\infty) \times Y, \\
W^- = (-\infty, R^-] \times Y, \\
W(\lambda^+,\lambda^-) = [R^-,R^+] \times Y
\end{eqnarray*}

In such $W$, we say that an almost complex structure $\overline{J}$ is asymptotically cylindrical if
\begin{eqnarray}
\overline{J} \ \mbox{coincides with } \ J^+ \in \mathcal{J}(K^+ \lambda^+) \ \mbox{in the region} \ W^+, \label{eq18} \\
\overline{J} \ \mbox{coincides with } \ J^- \in \mathcal{J}(K^- \lambda^-) \ \mbox{in the region} \ W^-, \label{eq19} \\
\overline{J} \ \mbox{is compatible with} \ \varpi \ \mbox{in} \ W(\lambda^+,\lambda^-),
\end{eqnarray}
where $K^+>0$ and $K^->0$ are constants.

For fixed $J^+ \in \mathcal{J}(K^+ \lambda^+)$ and $J^- \in \mathcal{J}(K^- \lambda^-)$, we denote by $\mathcal{J}(J^-,J^+)$ the
 set of asymptotically cylindrical almost complex structures in $(\mathbb{R} \times Y,\varpi)$ coinciding with $J^+$ on $W^+$ and $J^-$ on $W^-$.  It is well known that $\mathcal{J}(J^-,J^+)$ is non-empty and contractible.

\

\textbf{Remark:} In many references in the literature, a slightly different definition of asymptotically cylindrical almost complex structures is used: instead of demanding that $\overline{J}$ satisfies \eqref{eq18} and \eqref{eq19}, the stronger condition that $\overline{J}$ coincides with $J^{\pm} \in \mathcal{J}( \lambda^{\pm}) \ \mbox{in the region} \ W^{\pm}$ is required. We need to consider this more relaxed definition of cylindrical almost complex structures when we study the cobordism maps for strip Legendrian contact homologies.

\subsubsection{Splitting symplectic cobordisms} \label{splitting}

Let $\lambda^+$, $\lambda$ and $\lambda^-$ be contact forms on a contact manifold $(Y,\xi)$, and $f^{\lambda^+}_{\lambda^-}:Y \to (0,+\infty)$ and $f^{\lambda}_{\lambda^-}:Y \to (0,+\infty)$ be the functions such that $  f^{\lambda^+}_{\lambda^-}\lambda^- = \lambda^+$ and $  f^{\lambda}_{\lambda^-}\lambda^- = \lambda$. Assume that  $f^{\lambda^+}_{\lambda^-}>1$ and $f^{\lambda^+}_{\lambda^-}>  f^{\lambda}_{\lambda^-}\ >1$. 

For each $R>1$, let $\chi_R :\mathbb{R} \times Y \to \mathbb{R} $ such that:
\begin{eqnarray}
\partial_s \chi_R > 0 \ \mbox{in} \ [-R, R] \times Y,  \\
\chi_R(s) = e^{s-R-2} f^{\lambda^+}_{\lambda^-} \ \mbox{in} \ [R+2,+ \infty) \times Y, \\
\chi_R(s) = (1+\frac{s}{R}\epsilon)f^{\lambda}_{\lambda^-} \ \mbox{in} \ [-R, R] \times Y, \\
\chi_R(s) = e^{s+R+2} \ \mbox{in} \ (- \infty, -R-2] \times Y.
\end{eqnarray}
Then $(\mathbb{R} \times Y, d(\chi_R \lambda^-))$ is a straight exact symplectic cobordism from $\lambda^+$ to $\lambda^-$. As $R \to +\infty$ the region where the symplectic form $d(\chi_R \lambda^-)$ becomes similar to the symplectization of $(Y,\lambda)$ becomes arbitrarily large. To simplify notation we denote $(\mathbb{R} \times Y, d(\chi_R \lambda^-))$ by $W_R$.

% Then, for each $R>0$, we can construct an exact symplectic form $\varpi_R = d\kappa_R$ on $W = \mathbb{R} \times Y$ where:
%\begin{eqnarray}
%\kappa_R = e^{s-R-2}\lambda^+ \ \mbox{in} \ [R+2, + \infty) \times Y, \\
%\kappa_R = f(s)\lambda \ \mbox{in} \ [-R,R] \times Y, \\
%\kappa_R = e^{s+R+2}\lambda^- \ \mbox{in} \ (-\infty, -R-2] \times Y,
%\end{eqnarray}
%and $f: [-R,R]  \to [1-\epsilon,1+\epsilon]$ satisfies $f(-R) = 1-\epsilon$, $f(R) = 1+\epsilon$ and $f'>0$.
%In $ (\mathbb{R} \times Y, \varpi_R)$ we consider a compatible cylindrical almost complex structure $\widetilde{J}_R$; but we require an extra condition on $\widetilde{J}_R$:
%\begin{eqnarray}
%\widetilde{J}_R \ \mbox{coincides with} \ J \in \mathcal{J}(\lambda) \ \mbox{in} \ [-R,R] \times Y.
%\end{eqnarray}

Again we divide $W$ into regions:
$W_R^+ = [R+2, + \infty) \times Y$, $W_R(\lambda^+,\lambda)= [ R, R+2] \times Y$, $W_R(\lambda)=[-R, R] \times Y$, $W_R(\lambda,\lambda^-)=[-R-2,-R] \times Y$ and $W_R^-=(-\infty, -R-2] \times Y$. In $W_R$ we consider a compatible cylindrical almost complex structure $\widetilde{J}_R$; but we require an extra condition on $\widetilde{J}_R$:
\begin{eqnarray}
\widetilde{J}_R \ \mbox{coincides with} \ J \in \mathcal{J}(\lambda) \ \mbox{in} \ [-R,R] \times Y.
\end{eqnarray}
The family $(W_R, \widetilde{J}_R)_{R\in (0,+\infty)}$ of exact symplectic cobordisms with cylindrical almost complex structures is called a splitting family from $\lambda^+$ to $\lambda^-$ along $\lambda$.

%\begin{exmp}
%Let $\lambda^+$, $\lambda$ and $\lambda^-$ be contact forms in $Y$. Suppose now that $\lambda^+ = \zeta^+ \lambda^-$ and $\lambda = \zeta \lambda^-$ for functions $\zeta : Y \to \mathbb{R}$ and $\zeta^+ : Y \to \mathbb{R}$ satisfying $ \zeta^+ > \zeta + \epsilon > \zeta - \epsilon > 1$ for some $\epsilon>0$. 
%To gain an intuition about this construction, one can initially think that in the limit as $R \to +\infty$ the sequence $ (\mathbb{R} \times Y, d(\chi_R \lambda^-))$ splits into two exact symplectic cobordisms, $V(\lambda^+,\lambda)$ from $\lambda^+$ to $\lambda$, followed by $V(\lambda,\lambda^-)$ from $\lambda$ to $\lambda^-$. Actually, when one studies pseudoholomorphic curves in such cobordisms, the limiting object is more complicated then just the pair of two cobordims we mentioned; levels of symplectizations have to be inserted above $V(\lambda^+,\lambda)$, between $V(\lambda^+,\lambda)$ and $V(\lambda,\lambda^-)$, and below $V(\lambda,\lambda^-)$ to complete the picture. We refer again to the paper \cite{CPT} for a complete discussion of this topic.
%\end{exmp}

\subsubsection{Conical exact Lagrangian cobordisms}

Exact Lagrangian cobordisms play an important role in Symplectic Field Theory because they induce maps of topological invariants of this theory; see \cite{SFT}. In this paper we only use the simplest type of exact Lagrangian cobordisms which we call \textit{conical}.

Letting $\lambda^+$ and $\lambda^-$ be contact forms on a contact manifold $(Y,\xi)$, we assume that the function 
$f^{\lambda^+}_{\lambda^-}:Y \to (0,+\infty)$ such that $  f^{\lambda^+}_{\lambda^-}\lambda^- = \lambda^+$ satisfies $f^{\lambda^+}_{\lambda^-}>1$. Let then $W$ be a straight exact symplectic cobordism from $\lambda^+$ to $\lambda^-$. To each Legendrian submanifold ${\overline{\Lambda}}$ in $(Y,\xi)$ we can associate the cylinder $L_{\overline{\Lambda}} :=\mathbb{R}\times {\overline{\Lambda}} \subset W$. It is immediate to check that $L_{\overline{\Lambda}}$ is an exact Lagrangian submanifold of $W$ and we call it a \textit{conical exact Lagrangian cobordism} in $W$ from ${\overline{\Lambda}}$ to itself.

%Let $\Lambda^+$ be Legendrian submanifolds in, respectively, $(Y,\ker(\lambda^+))$ and $(Y,\ker(\lambda^-))$ and let $(\mathbb{R} \times Y,\varpi=d\kappa)$ be an exact symplectic cobordism from $\lambda^+$ to $\lambda^-$. We call a Lagrangian submanifold in $(\mathbb{R} \times Y,\varpi)$ a Lagrangian cobordism if there exists Legendrian submanifolds $\overline{\Lambda}^+$ in $(Y,\ker(\lambda^+))$ and $\overline{\Lambda}^-$ in $(Y,\ker(\lambda^-))$, and $N>0$ such that:

%\begin{eqnarray}
%L \cap ([N,+\infty) \times Y) = ([N,+\infty) \times \overline{\Lambda}^+), \\
%L \cap ((-\infty,-N] \times Y)= ((-\infty,-N] \times \overline{\Lambda}^-).
%\end{eqnarray}

%In this case we say that $L$ is a Lagrangian cobordism from $\overline{\Lambda}^+$ to $\overline{\Lambda}^-$. If such an $L$ is an exact Lagrangian submanifold of $(\mathbb{R} \times Y,d\kappa)$, we call it an exact Lagrangian cobordism from $\overline{\Lambda}^+$ to $\overline{\Lambda}^-$.

%\begin{exmp}
%If we take a Legendrian submanifold $\Lambda$ in $(Y,\ker(\lambda^-))$ then $\mathbb{R} \times \Lambda$ is an  exact Lagrangian submanifold in the symplectization of $(Y,\lambda^-)$. It is also an exact Lagrangian cobordism in $(\mathbb{R} \times Y, d(\chi_R \lambda^-))$ from $\Lambda $ to itself.
%\end{exmp}

\subsubsection{Pseudoholomorphic curves} \label{pseudohol}

Let $(S,i)$ be a closed Riemann surface, possibly with boundary, and $\Gamma \subset S$ be a finite set. We define $\Gamma_\partial:= \partial S  \cap \Gamma$.

Let $\lambda$ be a contact form in $(Y,\xi)$, $J \in \mathcal{J}(\lambda)$ and ${\overline{\Lambda}}$ be a Legendrian submanifold of $(Y,\xi)$. A finite energy pseudoholomorphic curve in the symplectization $(\mathbb{R} \times Y,J)$ with boundary in the Lagrangian cylinder $Z_{\overline{\Lambda}}:=\mathbb{R}\times {\overline{\Lambda}}$ over ${\overline{\Lambda}}$ is a  map $\widetilde{w}=(s,w): (S \setminus \Gamma; \partial(S) \setminus \Gamma_\partial) \to (\mathbb{R} \times Y;L)$ that satisfies
\begin{eqnarray}
\overline{\partial}_J(\widetilde{w}) := d\widetilde{w} \circ i - J \circ d\widetilde{w}=0, \\
\widetilde{w}(\partial S \setminus \Gamma_\partial) \subset Z_{\overline{\Lambda}},
\end{eqnarray}
and
\begin{equation}
0\leq E(\widetilde{w})<+ \infty,
\end{equation}
where $E(\widetilde{w}) := \sup_{q \in \mathcal{E}} \int_{S \setminus \Gamma} \widetilde{w}^*d(q\lambda)$ with $\mathcal{E}= \{ q: \mathbb{R} \to [0,1]; q' \geq 0\}$. The quantity $E(\widetilde{w})$ is called the Hofer energy and was introduced in \cite{H}. The operator $\overline{\partial}_J$ above is called the Cauchy-Riemann operator for the almost complex structure $J$.

For us, particularly important will be the case where $(S \setminus \Gamma,i)$ is biholomorphic to $(\mathbb{R} \times [0,1],i_0)$ (here $i_0$ denotes the complex structure in $\mathbb{C}$), $L=Z_{\overline{\Lambda}}$ for a Legendrian link $\overline{\Lambda}$ which is the union of two disjoint Legendrian connected Legendrian submanifolds $\Lambda$ and $\widehat{\Lambda}$, and $\widetilde{w}$ satisfies $\widetilde{w}(\{0\}\times \mathbb{R}) \subset Z_\Lambda$ and $\widetilde{w}(\{1\}\times \mathbb{R}) \subset Z_{\widehat{\Lambda}}$.  In this case $\widetilde{w}$ is called a pseudoholomorphic strip. By using a biholomorphism $\varphi: (\overline{\mathbb{D}^2} \setminus \{-1,1\},i_0) \to (\mathbb{R} \times [0,1],i_0)$ satisfying $ \varphi(H^+)= \{1\}\times \mathbb{R}$ (where $H_+ \subset (S^1 \setminus \{-1,1\})$ is the northern hemisphere of $S^1$) and $ \varphi(H^-)= \{0\}\times \mathbb{R}$ (where $H_- \subset (S^1 \setminus \{-1,1\})$ is the southern hemisphere of $S^1$) we can also view pseudoholomorphic strips as maps having as domain the closed disc with two punctures on the boundary. 

For a straight exact symplectic cobordism $W$ from a contact form $\lambda^+$ on $(Y,\xi)$ to a contact form $\lambda^-$ on $(Y,\xi)$, and $\overline{J} \in \mathcal{J}(J^-,J^+)$ where $J^+ \in \mathcal{J}(\lambda^+)$ and $J^- \in \mathcal{J}(\lambda^-)$, a finite energy pseudoholomorphic curve with boundary in a conical exact Lagrangian cobordism $L_{\overline{\Lambda}}$ over a Legendrian submanifold ${\overline{\Lambda}}$ of $(Y,\xi)$ is again a map $\widetilde{w}: (S \setminus \Gamma, \partial(S) \setminus \Gamma_\partial) \to (\mathbb{R} \times Y,L_{\overline{\Lambda}})$ that satisfies
\begin{eqnarray}
\overline{\partial}_J(\widetilde{w}) := d\widetilde{w} \circ i - J \circ d\widetilde{w}=0, \\
\widetilde{w}(\partial S \setminus \Gamma_\partial) \subset L_{\overline{\Lambda}},
\end{eqnarray}
and
\begin{equation}
0<E_{\lambda^-}(\widetilde{w}) + E_c(\widetilde{w}) + E_{\lambda^+} (\widetilde{w}) < +\infty,
\end{equation}
where
\\
$E_{\lambda^-}(\widetilde{w}) = \sup_{q \in \mathcal{E}} \int_{\widetilde{w}^{-1}(W^-))} \widetilde{w}^*d(q\lambda^-)$,
\\
$E_{\lambda^+}(\widetilde{w}) = \sup_{q \in \mathcal{E}} \int_{\widetilde{w}^{-1}(W^+)} \widetilde{w}^*d(q\lambda^+)$,
\\
$E_c(\widetilde{w}) = \int_{\widetilde{w}^{-1}(W(\lambda^+,\lambda^-))} \widetilde{w}^*d(h\lambda^-)$.
\\
These energies were also introduced in \cite{H}.

In splitting symplectic cobordisms we use a slightly modified version of energy. To simplify notation we denote by $\varpi_R$ the symplectic form $d(\chi\lambda^-)$ in $W_R$. Instead of demanding $0<E_-(\widetilde{w}) + E_c(\widetilde{w}) + E_+ (\widetilde{w}) < +\infty$ we demand
\begin{equation}
0<E_{\lambda^-}(\widetilde{w}) + E_{\lambda^-,\lambda}(\widetilde{w}) + E_{\lambda}(\widetilde{w})+ E_{\lambda,\lambda^+}(\widetilde{w})+ E_{\lambda^+} (\widetilde{w}) < +\infty
\end{equation}
where
\\
$ E_{\lambda}(\widetilde{w})= \sup_{q \in \mathcal{E}} \int_{\widetilde{w}^{-1}W(\lambda)} \widetilde{w}^*d(q\lambda)$,
\\
$E_{\lambda^-,\lambda}(\widetilde{w})= \int_{\widetilde{w}^{-1}(W(\lambda,\lambda^-))} \widetilde{w}^*\varpi_R$, \
\\
$E_{\lambda,\lambda^+}(\widetilde{w}=\int_{\widetilde{w}^{-1}(W(\lambda^+,\lambda))} \widetilde{w}^*\varpi_R$,
\\
and $E_{\lambda^-}(\widetilde{w})$ and $E_{\lambda^+}(\widetilde{w})$ are as above.

The elements of the set $\Gamma \subset S$ are called punctures of the pseudoholomorphic $\widetilde{w}$. We first divide $\Gamma$ in two classes: we call the elements of $\Gamma_\partial$ boundary punctures and the elements in $\Gamma \setminus \Gamma_\partial$ interior punctures. The work of Hofer \cite{H}, Hofer et al. \cite{HWZ} and Abbas \cite{Ab} allows us do classify the punctures in four different types\footnote{We assume here that our punctures are not removable singularities because pseudoholomorphic curves can be extended over punctures of this kind, as shown in \cite{Ab,H}.}. We will describe the behaviour of punctures for pseudoholomorphic curves in exact symplectic cobordisms, since the case of symplectizations is contained in this one as a particular case.

Before presenting this classification we introduce some notation: we let $B_\delta (z)$ be the ball of radius $\delta$ centered at the puncture $z$, and denote by $\mathfrak{b}_\delta (z)$ the set defined as the closure $\overline{\partial(B_\delta (z)) \cap int(S)}$ of the intersection  of the boundary of $B_\delta (z)$ with the interior of $S$. Notice that $\mathfrak{b}_\delta (z)$ is a circle or an interval, depending on whether $z$ is an interior or a boundary puncture. We can describe the types of punctures as follows:
\begin{itemize}
\item{$z \in \Gamma$ is a positive boundary puncture if $z \in \Gamma_\partial$ and $\lim_{z' \to z} s(z') = +\infty$, and there exists a sequence $\delta_n \to 0$ and a Reeb chord $c^+$ of  $X_{\lambda^+}$ from $\overline{\Lambda}$ to itself, such that $w(\mathfrak{b}_{\delta_n} (z))$ converges in $C^\infty $ to $c^+$ as $n\to +\infty$;}
\item{$z \in \Gamma$ is a negative boundary puncture if $z \in \Gamma_\partial$ and $\lim_{z' \to z} s(z') = -\infty$, and there exists a sequence $\delta_n \to 0$ and Reeb chord $c^-$ of  $X_{\lambda^-}$ from $\overline{\Lambda}$ to itself, such that $w(\mathfrak{b}_{\delta_n} (z))$ converges in $C^\infty $ to $c^-$ as $n\to +\infty$;}
\item{$z \in \Gamma$ is a positive interior puncture if $z \in \Gamma \setminus \Gamma_\partial$ and $\lim_{z' \to z} s(z') = +\infty$, and there exists a sequence $\delta_n \to 0$ and  Reeb orbit $\gamma^+$ of  $X_{\lambda^+}$, such that $w(\mathfrak{b}_{\delta_n} (z))$ converges in $C^\infty $ to $\gamma^+$ as $n\to +\infty$;}
\item{$z \in \Gamma$ is a negative interior puncture if $z \in \Gamma \setminus \Gamma_\partial$ and $\lim_{z' \to z} s(z') = -\infty$, and there exists a sequence $\delta_n \to 0$ and  Reeb orbit $\gamma^-$ of  $X_{\lambda^-}$, such that $w(\mathfrak{b}_{\delta_n} (z))$ converges in $C^\infty $ to $\gamma^-$ as $n\to +\infty$.}
\end{itemize}
The results in \cite{Ab,H,HWZ} imply that these are indeed the only possibilities we need to consider for the behaviour of the $\widetilde{w}$ near punctures. Intuitively, we have that at the punctures the pseudoholomorphic curve $\widetilde{w}$ detects Reeb chords and Reeb orbits. For a boundary (interior) puncture $z$, if there is a subsequence $\delta_n$ such that  $w(\mathfrak{b}_{\delta_n} (z))$ converges to a given Reeb chord $c$ (orbit $\gamma$), we will say that $\widetilde{w}$ is asymptotic to this Reeb chord $c$ (orbit $\gamma$) at the puncture $z$.

If a pseudoholomorphic curve is asymptotic to a transverse Reeb chord or a non-degenerate Reeb orbit at a puncture, more can be said about its asymptotic behaviour in neighbourhoods of this puncture.
In order to describe this behaviour for a boundary puncture $z$, we take a closed neighbourhood $U$ of $z$ that admits a holomorphic chart $\psi_U : (U \setminus \{z\}) \to \overline{\mathbb{R^+}} \times [0,1] \ \subset \mathbb{C}$, such that $\psi_U ((U \cap \partial(S)) \setminus \{z\}) = \mathbb{R^+} \times \{0\} \cup \mathbb{R^+} \times \{1\}$ and $\psi_U(\partial U \setminus \partial S) = \{0\} \times [0,1]$. In coordinates $(r,t) \in \mathbb{R} \times [0,1]$ we have $r(x) \to +\infty$ when $x$ tends to the puncture $z$. With this notation, it is shown in \cite{Ab}, that if $z$ is a positive boundary puncture on which $\widetilde{w}$ is asymptotic to a transverse Reeb chord $c^+$ of $X_{\lambda^+}$ from $\Lambda$ to\footnote{Here $\Lambda^+$ and  $\widehat{\Lambda}^+$ denote connected components of $\overline{\Lambda}^+$.} $\widehat{\Lambda}$, then $\widetilde{w} \circ \psi_u^{-1} (r,t) = (s(r,t),w(r,t))$ satisfies:

\begin{itemize}
\item{$w^r(t)=w(r,t)$ converges in $C^{\infty}$ to the Reeb chord $c^+$, uniformly in $t$ and exponentially in $r$.}
\end{itemize}
Similarly, if $z$ is a negative boundary puncture on which $\widetilde{w}$ is asymptotic to a transverse Reeb chord $c^-$ of $X_{\lambda^-}$ from $\Lambda$ to $\widehat{\Lambda}$, then $\widetilde{w} \circ \psi_u^{-1} (r,t) = (s(r,t),w(r,t))$ satisfies:
\begin{itemize}
\item{$w^r(t)=w(r,t)$ converges in $C^{\infty}$ to the Reeb chord $c^-$, uniformly in $t$ and exponentially in $r$.}
\end{itemize}

We discuss now the case where $z$ is an interior puncture. Take a neighbourhood $U \subset S$ of $z$ that admits a holomorphic chart $\psi_U : (U,z) \to (\mathbb{D},0)$. Using polar coordinates $(r,t) \in (0,+\infty) \times S^1$ we can write $x \in (\mathbb{D} \setminus 0)$ as $x = e^{-r}t $.
With this notation, it is shown in \cite{H,HWZ} that if $z$ is a positive interior puncture at which $\widetilde{w}$ is asymptotic to a non-degenerate Reeb orbit $\gamma^+$ of $X_{\lambda^+}$, then $\widetilde{w} \circ \psi_U^{-1} (r,t) = (s(r,t),w(r,t))$ satisfies
\begin{itemize}
\item{$w^r(t)=w(r,t)$ converges in $C^{\infty}$ to a Reeb orbit $\gamma^+$ of $X_{\lambda^+}$, exponentially in $r$ and uniformly in $t$.}
\end{itemize}
Similarly, if $z$ is a negative interior puncture at which $\widetilde{w}$ is asymptotic to a non-degenerate Reeb orbit $\gamma^-$ of $X_{\lambda^-}$, then $\widetilde{w} \circ \psi_u^{-1} (r,t) = (s(r,t),w(r,t))$ satisfies
\begin{itemize}
\item{$w^r(t)=w(r,t)$ converges in $C^{\infty}$ to a Reeb orbit $\gamma^-$ of $-X_{\lambda^-}$ as  $r \to +\infty$, exponentially in $r$ and uniformly in $t$.}
\end{itemize}

\textit{Remark: The exponential rate of convergence of pseudoholomorphic curves near punctures to Reeb orbits and Reeb chords is of crucial importance for the Fredholm theory (see \cite{Ab1,HWZ1}) that gives the dimension of the space of pseudoholomorphic curves with fixed asymptotic data, and is a consequence of asymptotic formulas obtained in \cite{Ab,HWZ}.}

The discussion above can be summarised by saying that near punctures the finite energy pseudoholomorphic curves detect Reeb orbits and Reeb chords. It is exactly this behavior that makes these objects useful for the study of dynamics of Reeb vector fields.

\textbf{Fact: As a consequence of the exactness of the symplectic cobordisms and the Lagrangian submanifolds considered above we obtain that the energy $E(\widetilde{w})$ of $\widetilde{w}$ satisfies $E(\widetilde{w}) \leq 5A(\widetilde{w})$ where $ A(\widetilde{w})$ is the sum of the action of the Reeb orbits and Reeb chords detected by the punctures of $\widetilde{w}$ counted with multiplicity; see for example \cite{CPT}.}

\subsection{Strip Legendrian contact homology} \label{sectionstrip}

We are now ready to define the strip Legendrian contact homology. First, we introduce the following notation: for a given contact form $\lambda$ we denote by $\mathcal{T}_{\Lambda \to \widehat{\Lambda}}(\lambda)$ the set of Reeb chords of $X_{\lambda}$ starting at $\Lambda$ and ending at $\widehat{\Lambda}$. To simplify notation we let $\mathcal{T}_{\Lambda}(\lambda)= \mathcal{T}_{\Lambda \to \Lambda}(\lambda)$.

\begin{defi} \label{defi:adapted}
Let $(Y,\xi)$ be a contact manifold, $\lambda_0$ be a contact form on $(Y,\xi)$ and $\Lambda$  and $\widehat{\Lambda}$ be two disjoint connected Legendrian submanifolds of $(Y\xi)$.
The contact form $\lambda_0$ is said to be \textit{adapted}   to the pair $(\Lambda,\widehat{\Lambda})$ if 
\begin{itemize}
\item (a) the Reeb flow of $\lambda_0$ has no contractible Reeb orbits,
\item (b) no Reeb chord in $\mathcal{T}_{\Lambda}(\lambda_0)$ vanishes in $\pi_1(Y,\Lambda)$,
\item (c)  no Reeb chord in $\mathcal{T}_{\widehat{\Lambda}}(\lambda_0)$ vanishes in $\pi_1(Y,\widehat{\Lambda})$,
\item (d) $(\lambda_0, \Lambda \to \widehat{\Lambda})$ is regular, i.e all Reeb chords in $\mathcal{T}_{\Lambda \to \widehat{\Lambda}}(\lambda_0)$  are transverse and do not intersect Reeb orbits of $\lambda_0$.
\end{itemize}
\end{defi}
\begin{remark} \label{rmk:nondegeneracy}
The attentive reader will notice that we do not demand that the Reeb orbits of $\lambda_0$ are non-degenerate or that the Reeb chords in  $\mathcal{T}_{\Lambda}(\lambda_0)$ and $\mathcal{T}_{\widehat{\Lambda}}(\lambda_0)$ are transverse. The reason for this is that because of conditions (a), (b) and (c) the pseudoholomorphic curves which we will use to define the strip Legendrian contact homology will not have any punctures asymptotic to Reeb orbits of $\lambda_0$ or Reeb chords in $\mathcal{T}_{\Lambda}(\lambda_0)$ and $\mathcal{T}_{\widehat{\Lambda}}(\lambda_0)$. Note however that, as we explain below, it is crucial for the definition of the the strip Legendrian contact homology that the elements in $\mathcal{T}_{\Lambda \to \widehat{\Lambda}}(\lambda_0)$ are transverse and do not intersect Reeb orbits of $\lambda_0$.
\end{remark}

For the remainder of section \ref{sectionstrip} we fix a contact manifold $(Y,\xi)$ and a contact form $\lambda_0$ on $(Y,\xi)$ adapted to a pair $(\Lambda,\widehat{\Lambda})$ of disjoint connected Legendrian submanifolds of $(Y,\xi)$.
Recall that $Z_\Lambda$ and $Z_{\widehat{\Lambda}}$ denote the exact Lagrangian cylinders over the Legendrian submanifolds $\Lambda$ and $\widehat{\Lambda}$ in the symplectization $(\mathbb{R}\times Y, d(e^s\lambda_0)$ of $\lambda_0$. 

In order to assign a $\mathbb{Z}_2$-grading to the Reeb chords in $\mathcal{T}_{\Lambda \to \widehat{\Lambda}}(\lambda_0)$, we use the Conley-Zehnder index. 
 We recall the definition of Conley-Zehnder index for Reeb chords as presented in \cite{E,EES}: we will use the fact that all Reeb chords in $\mathcal{T}_{\Lambda \to \widehat{\Lambda}}(\lambda_0)$ are transverse as stated in condition (d). For the definition, we first fix, once and for all, orientations for $\Lambda$ and $\widehat{\Lambda}$. Then, for each Reeb chord  $c \in \mathcal{T}_{\Lambda \to \widehat{\Lambda}}(\lambda_0)$, let $\Psi_c$ be a nowhere vanishing section of the vector bundle $\xi \mid_c$ that:
\begin{itemize}
\item{is tangent to $\Lambda$ on the initial point of $c$ and coincides with the orientation we fixed for $\Lambda$ at this initial point,}
\item{is tangent to $\widehat{\Lambda}$ on the final point of $c$ and coincides with the orientation we fixed for $\widehat{\Lambda}$ at this final point.}
 \end{itemize}
The section $\Psi_c$ induces a (unique up to homotopy) symplectic trivialisation of $(\xi \mid_c, d\lambda_0)$, which we also denote by $\Psi_c$.

Using the Reeb flow $\phi_{X_{\lambda_0}}$ we define a path of Lagrangian subspaces $\mathcal{Z} $ of $(\xi \mid_c,d\lambda_0)$. We consider the parametrisation $c:[0,T_c] \to Y$ of the Reeb chord $c$ given by the Reeb flow. Letting $D\phi_{X_{\lambda_0}}$ denote the linearisation of the Reeb flow, we define $\mathcal{Z}(t)$ to be the unique Lagragian subspace of $(\xi \mid_{c(t)},d\lambda_0)$  that contains $D\phi^t_{X_{\lambda_0}}(c(0))(q) $, where $q\in \xi \mid_{c(0)}$ is a vector tangent to $\Lambda$ and giving the orientation we chose for $\Lambda$.  After arriving at the endpoint $c(T_c)$, we complete $\mathcal{Z} $ to obtain a Lagrangian loop by making a continuous left-rotation of $\mathcal{Z}(T_c)$ (among Legendrian subspaces of $\xi\mid_{c(T_c)}$) till it meets the tangent space to $\widehat{\Lambda}$. With this completion and using our trivialisation $\Psi_c$, we associate to $\mathcal{Z} $ a path of Lagragian subspaces of the standard symplectic plane. The Conley-Zehnder index $\mu^{\Psi_c}_{CZ}(c)$ is defined to be the Maslov index of this path.

It is clear from the constructions above that, because we fixed the orientations of $\Lambda$ and $\widehat{\Lambda}$, the parity of $\mu^{\Psi_c}_{CZ}(c)$ is independent of the trivialisation ${\Psi_c}$. This allows us to define, for each $c \in \mathcal{T}_{\Lambda \to \widehat{\Lambda}}(\lambda_0)$, its $\mathbb{Z}_2$-grading by $\mid c \mid = (\mu^{\Psi_c}_{CZ}(c) + 1)\mod2$. We call chords with grading $0$ even chords, and chords with grading $1$ odd chords.

Let $LCH_{st}(\lambda_0,\Lambda \to \widehat{\Lambda})$ be the $\mathbb{Z}_2$ vector-space generated by $\mathcal{T}_{\Lambda \to \widehat{\Lambda}}(\lambda_0)$. We denote by $LCH_{st,odd}(\lambda_0,\Lambda \to \widehat{\Lambda})$ the subspace of $LCH_{st}(\lambda_0,\Lambda \to \widehat{\Lambda})$ generated by odd chords, and $LCH_{st,even}(\lambda_0,\Lambda \to \widehat{\Lambda})$ the subspace generated by even chords.
Given two Reeb chords $c_1$ and $c_2$, and an almost complex structure $J \in \mathcal{J}(\lambda_0)$, we know from the previous section that it makes sense to consider the moduli space $\mathcal{M}(c_1,c_2;J;\Lambda, \widehat{\Lambda})$ whose elements are equivalence classes of finite energy pseudoholomorphic strips $\widetilde{w}: (\overline{D} \setminus \{-1,1\},i_0) \to (\mathbb{R} \times Y,J)$ satisfying:
\begin{itemize}
\item{$1$ is a positive boundary puncture, and $\widetilde{w}$ is asymptotic to $c_1$ at $1$,}
\item{$-1$ is a negative boundary puncture, and $\widetilde{w}$ is asymptotic to $c_2$ at $-1$,}
\item{ $\widetilde{w}(H_-) \subset Z_\Lambda$,}
\item{ $\widetilde{w}(H_+) \subset Z_{\widehat{\Lambda}}$.}
\end{itemize}
Two pseudoholomorphic strips $\widetilde{w}$ and $\widetilde{w}'$ satisfying these four conditions represent the same element in $\mathcal{M}(c_1,c_2;J; \Lambda, \widehat{\Lambda})$ if there exists a biholomorphism $\varphi$ of the disc that fixes $\{-1,1\}$ such that $\widetilde{w} \circ \varphi =\widetilde{w}'$.

It follows from Abbas' asymptotic analysis \cite{Ab} explained in Section \ref{pseudohol} that all the elements of the moduli space $\mathcal{M}(c_1,c_2;J;\Lambda, \widehat{\Lambda})$ are somewhere injective pseudoholomorphic curves. It is well known that the linearization $D\overline{\partial}_J$ at any element $\mathcal{M}(c_1,c_2;J;\Lambda, \widehat{\Lambda})$ is a Fredholm map (we remark that this property is valid for more general moduli spaces of curves with prescribed asymptotic behaviour). 
In \cite[Proposition 3.15]{DR} Dimitroglou Rizell  combines the techniques of \cite{Ab1}, Dragnev \cite{Dr} and Lazzarini \cite{Laz} and showed that for a generic set $\mathcal{J}_{reg}(\lambda_0,\Lambda,\widehat{\Lambda}) \subset \mathcal{J}(\lambda_0)$ all the elements in $\mathcal{M}(c_1,c_2;J;\Lambda, \widehat{\Lambda})$ are transverse in the sense that the linearization $D\overline{\partial}_J$ of the Cauchy-Riemann operator $\overline{\partial}_J$ at the elements of $\mathcal{M}(c_1,c_2;J;\Lambda,\widehat{\Lambda})$ is surjective, this being valid for all Reeb chords $c_1$ and $c_2$ from $\Lambda$ to $\widehat{\Lambda}$. The key reason for this to hold is that the Reeb chords $c_1$ and $c_2$ are embedded, which implies that the elements  $\mathcal{M}(c_1,c_2;J;\Lambda,\widehat{\Lambda})$ are somewhere injective pseudoholomorphic curves.

Thus, in the case where $J \in \mathcal{J}_{reg}(\lambda_0,\Lambda,\widehat{\Lambda})$ one can use the implicit function theorem, and obtain that any connected component of the moduli space $\mathcal{M}(c_1,c_2;J;\Lambda, \widehat{\Lambda})$ is a finite dimensional manifold, and its dimension is given by the Fredholm index $I_F$ of $D\overline{\partial}_J$ computed at any element of this connected component of $\mathcal{M}(c_1,c_2;J;\Lambda, \widehat{\Lambda})$. We let $\mathcal{M}^k (c_1,c_2;J;\Lambda, \widehat{\Lambda})\subset \mathcal{M} (c_1,c_2;J;\Lambda, \widehat{\Lambda})$ be the moduli space of pseudoholomorphic strips with Fredholm index $k$.

It follows from the formula in \cite{Ab1} for the Fredholm index $I_F$ of the linearised $D\overline{\partial}_J$ operator over a strip in $\mathcal{M}(c_1,c_2;J;\Lambda, \widehat{\Lambda})$, that $I_F$ has the same parity as $\mid c_1 \mid + \mid c_2 \mid $.
It follows from Stokes' Theorem that $0 \leq \int_{\overline{D} \setminus \{-1,1\}} \widetilde{w}^{*} (d\lambda_0) = A(c_1) - A(c_2) $, and therefore $\mathcal{M}(c_1,c_2;J;\Lambda, \widehat{\Lambda})$ can only be non-empty if $A(c_1) \geq A(c_2)$. Because of the $\mathbb{R}$-invariance of the almost complex structure $J$ there is an $\mathbb{R}$-action on the spaces $\mathcal{M}(c_1,c_2;J;\Lambda, \widehat{\Lambda})$, and we let $\widetilde{\mathcal{M}}(c_1,c_2;J;\Lambda, \widehat{\Lambda})= \mathcal{M}(c_1,c_2;J;\Lambda, \widehat{\Lambda})/\mathbb{R}$.

We are now ready to define a differential $d_J$ in $LCH_{st}(\lambda_0,\Lambda \to \widehat{\Lambda})$.
\begin{defi}
Let $c \in \mathcal{T}_{\Lambda \to \widehat{\Lambda}}(\lambda_0)$ and $J \in \mathcal{J}_{reg}(\lambda_0) \subset \mathcal{J}(\lambda_0)$. We define:
\begin{equation}
d_J (c):= \sum_{c' \in \mathcal{T}_{\Lambda \to \widehat{\Lambda}}(\lambda_0)}  [n_{c,c'}\mod2]  c'
\end{equation}
where $n_{c,c'}$ is the cardinality of the moduli space $\mathcal{M}^1 (c,c';J;\Lambda, \widehat{\Lambda})$ of pseudoholomorphic strips of Fredholm index $1$ modulo the $\mathbb{R}$-action.

The differential is extended to $LCH_{st}(\lambda_0,\Lambda \to \widehat{\Lambda})$ by linearity.
\end{defi}

To complete the construction of the strip Legendrian contact homology, we  must prove that $d_J$ is well-defined and that $d_J \circ d_J=0$. Before proceeding with the proofs of these results we discuss the intuition behind the definition of this homology. The strip Legendrian contact homology can be seen as a relative version of the cylindrical contact homology (see \cite{B,SFT}). For cylindrical contact homology to be well-defined for a contact form, this contact form has to have some special properties; for example, for a hypertight contact form (i.e. one that does not have contractible periodic orbits) cylindrical contact homology is well-defined. As we will see later, the non-existence of contractible Reeb orbits precludes the ``bubbling'' of pseudoholomorphic planes. This, together with SFT-compactness, implies that if its asymptotic orbits are in a primitive homotopy class, a sequence of pseudoholomorphic cylinders of Fredholm index 2 can only break in a pseudoholomorphic building of 2 levels, each containing a cylinder of Fredholm index 1. This implies that only such buildings can appear in the boundary of the compactified moduli space of pseudoholomorphic cylinders of Fredholm index $2$. This description of the compactified moduli spaces of pseudoholomorphic cylinders of index $2$, is the crucial step that allows us to define cylindrical contact homology with coefficients in $\mathbb{Z}_2$.

The strip Legendrian contact homology is the relative analogue of the cylindrical contact homology. This time the differential involves pseudoholomorphic strips with boundary conditions on Lagrangian cylinders over Legendrian submanifolds. For such a theory to be well-defined we have to preclude not only ``bubbling'' of planes but also of pseudoholomorphic half-planes . The conditions (b) and (c) above serve exactly to make impossible such ``bubbling'' phenomena, and the condition (d) is a non-degeneracy condition. Under these hypotheses it is possible to define the strip Legendrian contact homology, and to carry out this construction one uses results on the analytical properties of pseudoholomorphic strips and discs. For these results we refer to: \cite{Ab} for the necessary results on the asymptotic behaviour of punctures, \cite{Ab1} for the necessary results on Fredholm theory,  \cite{Ab1,Dr} for the necessary transversality results, \cite{EES} for (essentially) the necessary techniques to perform gluing. We now proceed to prove:
\begin{lemma} \label{lemma1}
For $J \in \mathcal{J}_{reg}(\lambda_0,\Lambda, \widehat{\Lambda}) \subset \mathcal{J}(\lambda_0)$, and $d_J$ defined above we have:
\begin{itemize}
\item{(1) $d_J$ is well defined,}
\item{(2) $d_J$ decreases the action of Reeb chords,}
\item{(3) for each $c \in \mathcal{T}_{\Lambda \to \widehat{\Lambda}}(\lambda_0)$, $d_J(c)$ is a finite sum,}
\item{(4) $d_J:LCH_{st,odd}(\lambda_0,\Lambda \to \widehat{\Lambda}) \to LCH_{st,even}(\lambda_0,\Lambda \to \widehat{\Lambda})$ and $d_J:LCH_{st,even}(\lambda_0,\Lambda \to \widehat{\Lambda}) \to LCH_{st,odd}(\lambda_0,\Lambda \to \widehat{\Lambda}).$}
\end{itemize}
\end{lemma}

\textit{Proof:} In order for $d_J$ to be well-defined we have to prove that the moduli space $\widetilde{\mathcal{M}}^1 (c,c';J;\Lambda, \widehat{\Lambda})$ is finite for every $c$ and $c'$. Because $J \in \mathcal{J}_{reg}(\lambda_0,\Lambda, \widehat{\Lambda})$, it follows that $\widetilde{\mathcal{M}}^1 (c,c';J;\Lambda, \widehat{\Lambda})$ is a $0$-dimensional manifold. If we show that it is compact then it has to be a finite set.
To obtain the compactness we will apply the ``bubbling of'' analysis for pseudoholomorphic curves in symplectizations of \cite{Ab,H} and the SFT-compactness results of \cite{CPT}.

Let $\widetilde{w}_n$ be a sequence of elements of $\widetilde{\mathcal{M}}^1 (c,c';J;\Lambda, \widehat{\Lambda})$. Because of the assumptions we made on the contact form $\lambda_0$, the sequence $\widetilde{w}_n$ cannot have interior bubbling points: by \cite{H} we know that an interior bubbling point would imply the existence of a finite energy plane and thus of a contractible periodic orbit of $X_{\lambda_0}$, something that contradicts (a). Boundary bubbling points are also forbidden: it follows from \cite{Ab,CPT} that they would give rise to either a pseudoholomorphic disc with boundary in $Z_\Lambda$, a pseudoholomorphic disc with boundary on $Z_{\widehat{\Lambda}}$, a pseudoholomorphic disc with only one puncture asymptotic to a Reeb chord from $\Lambda$ to itself, or a a pseudoholomorphic disc with only one puncture asymptotic to a Reeb chord from $\widehat{\Lambda}$ to itself. The first two possibilities are impossible because $Z_{{\Lambda}}$ and $Z_{\widehat{\Lambda}}$ are exact Lagrangian submanifolds of the symplectization $(\mathbb{R}\times Y, d(e^s\lambda_0))$ and the later two because they would contradict conditions (b) and (c) satisfied by $\lambda_0$. Combining this information with the SFT-compactness results of \cite{CPT} we have that $\widetilde{w}_n$ converges in the SFT sense to a pseudoholomorphic building $\widetilde{w}$ with k-levels $\widetilde{w}^l$, with all levels $\widetilde{w}^l$ being pseudoholomorphic strips with non-vanishing Hofer energy that satisfy:
\begin{itemize}
\item{$1$ is a positive boundary puncture, and $\widetilde{w}^l$ is asymptotic to $c_l \in \mathcal{T}_{\Lambda \to \widehat{\Lambda}}(\lambda_0) $ at $1$,}
\item{$-1$ is a negative boundary puncture, and $\widetilde{w}$ is asymptotic to $c_{l+1} \in \mathcal{T}_{\Lambda \to \widehat{\Lambda}}(\lambda_0)$ at $-1$,}
\item{$\widetilde{w}(H_+) \subset Z_{\widehat{\Lambda} }$, where $H_+ \subset (S^1 \setminus \{-1,1\})$ is the northern hemisphere,}
\item{$\widetilde{w}(H_-) \subset Z_{\Lambda}$, where $H_- \subset (S^1 \setminus \{-1,1\})$ is the southern hemisphere,}
\end{itemize}
where $c_1 = c$ and $c_{k+1} = c'$. Because every $\widetilde{w}^l$ is somewhere injective and has non-zero Hofer energy, we have that the indices $I_F(\widetilde{w}^l)$ satisfy $I_F(\widetilde{w}^l) \geq 1$. This implies that $I_F (\widetilde{w})= \sum (I_F(\widetilde{w}^l) \geq l$. On the other hand as $\widetilde{w}$ is the limit of a sequence of pseudoholomorphic strips of Fredholm index 1, it has to satisfy $I_F(\widetilde{w})=1$. We thus conclude that $l=1$, and $\widetilde{w}\in \widetilde{\mathcal{M}}^1 (c,c';J;\Lambda, \widehat{\Lambda})$, which implies the desired compactness. This proves that $n_{c,c'}$ is finite for every $c,c' \in \mathcal{T}_{\Lambda \to \widehat{\Lambda}}(\lambda_0)$, and thus that $d_J$ is well defined.

To verify item (2), we remark that given $c \in \mathcal{T}_{\Lambda \to \widehat{\Lambda}}(\lambda_0)$, the number $n_{c,c'}$ can only be non-zero for Reeb chords $c'$ that satisfy $A(c') < A(c)$. This implies that $d_J$ decreases the action of Reeb chords.

By the non-degeneracy condition (d) one obtains that the set of Reeb chords with action smaller then $A(c)$ is finite, and so $n_{c,c'}$ is non-zero only for a finite number of $c'$. This implies that $d_J(c)$ is a finite sum, proving item (3).

 Item (4) follows easily from the fact that the Fredholm index of a strip connecting two chords $c$ and $c'$ has the same parity as $\mid c \mid + \mid c' \mid $, as this implies  $\widetilde{\mathcal{M}}^1 (c,c';J;\Lambda, \widehat{\Lambda})$  can only be non-empty if $c$ and $c'$ have different parity.
\qed

In \cite{CPT} it is proved that the moduli spaces $\widetilde{\mathcal{M}}^k (c,c';J;\Lambda, \widehat{\Lambda})$ admit a compactification. The compactfied moduli space is composed not only of pseudoholomorphic curves, but also of pseudoholomorphic buildings. We will denote this compactification of $\widetilde{\mathcal{M}}^k (c,c';J;\Lambda, \widehat{\Lambda})$ by $\overline{\mathcal{M}}^k (c,c';J;\Lambda, \widehat{\Lambda})$.

\begin{lemma} \label{lemma2}
For $J \in \mathcal{J}_{reg}(\lambda_0) \subset \mathcal{J}(\lambda_0)$, and $d_J$ as defined before we have
$d_J \circ d_J = 0$.
\end{lemma}
\textit{Proof:} The lemma will be a consequence of the description we will give of the compactified moduli space $\overline{\mathcal{M}}^2 (c,c';J;\Lambda, \widehat{\Lambda})$ of pseudoholomorphic strips with Fredholm index $2$. Because of regularity of $J$, it will follow that for all $c,c' \in  \mathcal{T}_{\Lambda \to \widehat{\Lambda}}(\lambda_0)$, $\overline{\mathcal{M}}^2 (c,c';J;\Lambda, \widehat{\Lambda})$ is either empty, or the finite union of disjoint circles and closed intervals. We summarise this in the following claim.

\textit{Compactness Claim: Suppose $\overline{\mathcal{M}}^2 (c,c';J;\Lambda, \widehat{\Lambda})$ is non-empty. Then, each connected component $I$ of $\overline{\mathcal{M}}^2 (c,c';J;\Lambda, \widehat{\Lambda})$ is either a circle or a closed interval. Moreover, when $I$ is homeomorphic to a closed interval, its boundary is composed by pseudoholomorphic buildings $\widetilde{w}$ with 2 levels $\widetilde{w}^1$ and $\widetilde{w}^2$ satisfying:}
\begin{center}
{\textit{$\widetilde{w}^1 \in \overline{\mathcal{M}}^1 (c,\check{c};J;\Lambda, \widehat{\Lambda})$ and $\widetilde{w^2} \in \overline{\mathcal{M}}^1 (\check{c},c';J;\Lambda, \widehat{\Lambda})$ for some $\check{c} \in \mathcal{T}_{\Lambda \to \widehat{\Lambda}}(\lambda_0)$}.}
\end{center}

{\flushleft Before proving this claim we use it to prove Lemma~\ref{lemma2}. For this, we write}
\begin{equation}
d_J \circ d_J(c) = \sum_{c' \in \mathcal{T}_{\Lambda \to \widehat{\Lambda}}(\lambda_0)} (m_{c,c'} \mod 2)c'.
\end{equation}
It is clear that the lemma will follow if we can prove that $m_{c,c'}$ is even. On one hand, it follows from our definition of $d_J$ that $m_{c,c'}$ counts the number of 2-level pseudoholomorphic buildings  whose levels $\widetilde{w}^1$ and $\widetilde{w}^2$ satisfy
$\widetilde{w}^1 \in \overline{\mathcal{M}}^1 (c,\check{c};J;\Lambda, \widehat{\Lambda})$ and $\widetilde{w^2} \in \overline{\mathcal{M}}^1 (\check{c},c';J;\Lambda, \widehat{\Lambda})$ for some $\check{c} \in \mathcal{T}_{\Lambda \to \widehat{\Lambda}}(\lambda_0)$. This together with the compactness claim implies that the number of boundary components of $\overline{\mathcal{M}}^2 (c,c';J)$ is smaller or equal to  $m_{c,c'}$.

On the other hand, because of the regularity of $J \in \mathcal{J}_{reg}(\lambda_0,\Lambda, \widehat{\Lambda})$, we can apply the gluing theorem for Fredholm regular pseudoholomorphic strips: it implies that if $\widetilde{w}$ is 2-level pseudoholomorphic building whose levels $\widetilde{w}^1$ and $\widetilde{w}^2$ satisfy
$\widetilde{w}^1 \in \overline{\mathcal{M}}^1 (c,\check{c};J;\Lambda, \widehat{\Lambda})$ and $\widetilde{w}^2 \in \overline{\mathcal{M}}^1 (\check{c},c';J;\Lambda, \widehat{\Lambda})$ for some $\check{c} \in \mathcal{T}_{\Lambda \to \widehat{\Lambda}}(\lambda_0)$, then $\widetilde{w}$ represents exactly one point in the boundary of $\overline{\mathcal{M}}^2 (c,c';J;\Lambda, \widehat{\Lambda})$.
We thus have that $m_{c,c'}$  is bigger or equal to the number of boundary points of $\overline{\mathcal{M}}^2 (c,c';J)$.

Summarising, the combination of the Compactness Claim and the gluing theorem allows us to conclude that the number $m_{c,c'}$ is exactly the number of boundary components of the moduli space $\overline{\mathcal{M}}^2 (c,c';J;\Lambda, \widehat{\Lambda})$. Because $\overline{\mathcal{M}}^2 (c,c';J;\Lambda, \widehat{\Lambda})$ is a finite union of  disjoint intervals and circles, this number is even. This finishes the proof of the lemma modulo the Compactness claim.

\

\textit{Proof of Compactness Claim:}  Suppose $\overline{\mathcal{M}}^2 (c,c';J;\Lambda, \widehat{\Lambda})$ is non-empty and let $I$ be one of its connected components. It follows from the regularity of $J$ that the interior $\dot{I}$ of $I$ is a 1-dimensional manifold. We have now two possibilities, either $I$ is compact or not. If $\dot{I}$ is compact then it is a circle.

If that is not the case, let $\widetilde{w}_n$ be a sequence of elements of $\dot{I}$ converging to the boundary of $I$. Reasoning as we did in the proof of the Lemma \ref{lemma1}, we obtain that no ``bubbling'' can occur. Thus the SFT-compactness theorem of \cite{CPT} implies that $\widetilde{w}_n$ converges to a pseudoholomorphic building $\widetilde{w}$ with k-levels $\widetilde{w}^l$, such that all levels $\widetilde{w}^l$ are pseudoholomorphic strips satisfying:
\begin{itemize}
\item{$1$ is a positive boundary puncture, and $\widetilde{w}^l$ is asymptotic to $c_l \in \mathcal{T}_{\Lambda \to \widehat{\Lambda}}(\lambda_0) $ at $1$,}
\item{$-1$ is a negative boundary puncture, and $\widetilde{w}$ is asymptotic to $c_{l+1} \in \mathcal{T}_{\Lambda \to \widehat{\Lambda}}(\lambda_0)$ at $-1$,}
\item{$\widetilde{w}(H^+) \subset Z_{\widehat{\Lambda}}$ where $H_- \subset (S^1 \setminus \{-1,1\})$ is the northern hemisphere,}
\item{$\widetilde{w}(H_-) \subset Z_\Lambda$ where $H_- \subset (S^1 \setminus \{-1,1\})$ is the southern hemisphere,}
\end{itemize}
where $c_1 = c$ and $c_{k+1} = c'$. Again, because every $\widetilde{w}^l$ is somewhere injective, we have that the Fredholm indices of the levels satisfy $F(\widetilde{w}^l) \geq 1$, and it follows that $I_F (\widetilde{w})= \sum (I_F(\widetilde{w}^l) \geq l$. On the other hand as $\widetilde{w}$ is the limit of a sequence of pseudoholomorphic strips of Fredholm index 2, it has to satisfy $I_F(\widetilde{w})=2$.

We have then two possibilities: either $l=1$ and $\widetilde{w} \in \overline{\mathcal{M}}^2 (c,c';J;\Lambda, \widehat{\Lambda})$; or $l=2$ which forces $I_F(\widetilde{w}^1)=I_F(\widetilde{w}^2)=1$, for $\widetilde{w}^1 \in \overline{\mathcal{M}}^1 (c,c_2;J;\Lambda, \widehat{\Lambda})$ and $\widetilde{w}^2 \in \overline{\mathcal{M}}^1 (c_2,c';J;\Lambda, \widehat{\Lambda})$ for some $c_2 \in \mathcal{T}_{\Lambda \to \widehat{\Lambda}}$.  The first case is ruled out because we assumed that $\widetilde{w}_n$ converges to the boundary of $I$.
We have thus obtained that all the elements on the boundary of  $I$ are 2-level pseudoholomorphic buildings with the properties claimed, and an analysis identical to the one done in Lemma \ref{lemma1} shows that the moduli space of such buildings is a compact finite set.  We conclude that the boundary of $I$ is a compact 0-dimensional manifold.

The gluing theorem gives a description of a small neighbourhood $U$ of a 2-level pseudoholomorphic building $\widetilde{w}$  appearing in the boundary of a connected component $I$ of the moduli space $\overline{\mathcal{M}}^2 (c,c';J;\Lambda, \widehat{\Lambda})$. The neighbourhood $U$ admits a homeomorphism to the interval $[0, +\infty)$, that takes $0$ to the 2-level building $\widetilde{w}$  and all other values to pseudoholomorphic strips in the interior of $I$.

Summing up, we have that the connected component $I$ of $\overline{\mathcal{M}}^2 (c,c';J;\Lambda, \widehat{\Lambda})$ has the structure of a one-dimensional manifold with boundary, i.e a closed interval. This finishes the proof of the compactness claim.
\qed

 We  denote  by  $LC\mathbb{H}_{st}(\lambda_0,\Lambda \to \widehat{\Lambda})$  the  homology  associated to  $(LCH_{st}(\lambda_0,\Lambda \to \widehat{\Lambda}), d_J)$.

\subsubsection{Strip Legendrian contact homology in special homotopy classes}

The free homotopy classes of paths starting at $\Lambda$ and ending at $\widehat{\Lambda}$ generate subcomplexes of $LC\mathbb{H}_{st}(\lambda_0,\Lambda \to \widehat{\Lambda})$. To formalize this we denote by $\Sigma_{\Lambda \to \widehat{\Lambda}}$ the set of homotopy classes of paths starting at $\Lambda$ and ending at $\widehat{\Lambda}$. For our contact form $\lambda_0$ and an element $\rho \in \Sigma_{\Lambda \to \widehat{\Lambda}}$ we denote by $\mathcal{T}^{\rho}_{\Lambda \to \widehat{\Lambda}}(\lambda_0)$ the set of Reeb chords from $\Lambda$ to $\widehat{\Lambda}$ that belong to $\rho$.

It is clear that for all $c \in \mathcal{T}^{\rho}_{\Lambda \to \widehat{\Lambda}}(\lambda_0)$, the terms $[n_{c,c'}\mod2]$ appearing in the differential $d_J(c)  = \sum_{c' \in \mathcal{T}_{\Lambda \to \widehat{\Lambda}}(\lambda_0)} [n_{c,c'}\mod2] c'$, can only be non-zero if $c' \in \mathcal{T}^{\rho}_{\Lambda \to \widehat{\Lambda}}(\lambda_0)$. This implies that the vector spaces $LCH^{\rho}_{st}(\lambda_0,\Lambda \to \widehat{\Lambda})$ are subcomplexes of the chain complex $(LCH_{st}(\lambda_0,\Lambda \to \widehat{\Lambda}),d_J)$ and that:
\begin{equation}
LC\mathbb{H}_{st}(\lambda_0,\Lambda \to \widehat{\Lambda}) = \bigoplus_{\rho \in \Sigma_{\Lambda \to \widehat{\Lambda}}}LC\mathbb{H}^{\rho}_{st}(\lambda_0,\Lambda \to \widehat{\Lambda}),
\end{equation}
where $LC\mathbb{H}^{\rho}_{st}(\lambda_0,\Lambda \to \widehat{\Lambda})$ denotes the homology of the chain-complex $(LCH^{\rho}_{st}(\lambda_0,\Lambda \to \widehat{\Lambda}),d_J)$.

\subsubsection{Cobordism maps}

Symplectic cobordisms play a crucial role in SFT because they induce maps between the SFT invariants. We now explain such maps in our setting.

\begin{prop} \label{propcobordism}
Let $\lambda^+$ and $\lambda^-$ be contact forms on  $(Y,\xi)$ adapted to a pair of disjoint connected Legendrian submanifolds $\Lambda$  and $\widehat{\Lambda}$. Let $W$ be a straight exact symplectic cobordism from $\lambda^+$ to $\mathrm{k}\lambda^-$ for some constant $\mathrm{k}>0$, $L_\Lambda \subset W$ be the conical exact Lagrangian cobordism from $\Lambda$ to itself, and ${L}_{\widehat{\Lambda}} \subset W$ be the conical exact Lagrangian cobordism from $\widehat{\Lambda}$ to itself. Then, these cobordisms induce a map $\Phi_{W,L_\Lambda,{L}_{\widehat{\Lambda}}}$ from $LC\mathbb{H}_{st}(\lambda^+,\Lambda \to \widehat{\Lambda})$ to $LC\mathbb{H}_{st}(\mathrm{k}\lambda^-,\Lambda \to \widehat{\Lambda})$.
Moreover the map $\Phi_{W,L_\Lambda,{L}_{\widehat{\Lambda}}}$ respects the filtration of $LC\mathbb{H}_{st}(\lambda^+,\Lambda \to \widehat{\Lambda})$ by  homotopy classes, i.e for each homotopy class $\rho\in \Sigma_{\Lambda \to \widehat{\Lambda}}$ we have that  $\Phi_{W,L_\Lambda,{L}_{\widehat{\Lambda}}}$ restricts to a map from $LC\mathbb{H}^\rho_{st}(\lambda^+,\Lambda \to \widehat{\Lambda})$ to $LC\mathbb{H}^\rho_{st}(\mathrm{k}\lambda^-,\Lambda \to \widehat{\Lambda})$.
\end{prop}

\textit{Proof:} Taking almost complex structures $J^+\in \mathcal{J}_{reg}(\lambda^+)$ and $J^- \in \mathcal{J}_{reg}(\lambda^-)$, we can define the homologies $LC\mathbb{H}_{st}(\lambda^+,\Lambda \to \widehat{\Lambda})$ and $LC\mathbb{H}_{st}(\lambda^-,\Lambda \to \widehat{\Lambda})$. The idea to define the cobordism map is to choose $J_W \in \mathcal{J}(J^+,J^-)$ and define $\Phi_{W,L_\Lambda,{L}_{\widehat{\Lambda}}}$ by counting pseudoholomorphic strips $\widetilde{w}:(\overline{D} \setminus \{-1,1\},i_0) \to (W,J_W)$ with Fredholm index $0$, having $1$ as a positive boundary puncture asymptotic to a Reeb chord $c^+ \in \mathcal{T}_{\Lambda \to \widehat{\Lambda}}(\lambda^+)$, $-1$ as a negative boundary puncture asymptotic to a Reeb chord $c^- \in \mathcal{T}_{\Lambda \to \widehat{\Lambda}}(\lambda^-)$, and with boundary conditions $\widetilde{w}(H^-) \subset L_\Lambda$ and $\widetilde{w}(H^+) \subset L_{\widehat{\Lambda}}$.
In order for this approach to work we must first choose $J_W$ that satisfy certain regularity assumptions.

For a pair of Reeb chords $c^+ \in \mathcal{T}_{\Lambda \to \widehat{\Lambda}}(\lambda^+)$ and $c^- \in \mathcal{T}_{\Lambda \to \widehat{\Lambda}}(\lambda^-)$, we consider the moduli space $\mathcal{M}(c^+,c^-;J_W;L_\Lambda,L_{\widehat{\Lambda}})$ whose elements are equivalence classes of pseudoholomorphic strips $\widetilde{w}:(\overline{D} \setminus \{-1,1\},i_0) \to (W,J_W)$ having $1$ as a positive boundary puncture asymptotic to a Reeb chord $c^+ \in \mathcal{T}_{\Lambda \to \widehat{\Lambda}}(\lambda^+)$, $-1$ as a negative boundary puncture asymptotic to a Reeb chord $c^- \in \mathcal{T}_{\Lambda \to \widehat{\Lambda}}(\lambda^-)$, and with boundary conditions $\widetilde{w}(H^-) \subset L_\Lambda$ and $\widetilde{w}(H^+) \subset L_{\widehat{\Lambda}}$. Two such pseudoholomorphic strips $\widetilde{w}$ and $\widetilde{w}'$ represent the same element in $\mathcal{M}(c^+,c^-;J_W;L_\Lambda,L_{\widehat{\Lambda}})$ if one can be obtained from the other via a biholomorphic reparametrization of the disk that fixes $\{-1,1\}$.
 Again, we let $\mathcal{M}^k(c^+,c^-;J_W;L_\Lambda,L_{\widehat{\Lambda}}) \subset \mathcal{M}(c^+,c^-;J_W;L_\Lambda,L_{\widehat{\Lambda}})$ be the set of pseudoholomorphic strips of  Fredholm index $k$. All the strips in these moduli spaces are somewhere injective curves. More precisely it shown in \cite[Theorem 2.14]{CRGG} that the set of injective points in each pseudoholomorphic curve in $ \mathcal{M}(c^+,c^-;J_W;L_\Lambda,L_{\widehat{\Lambda}})$ is open and dense in the domain. It the follows from the techniques of \cite{Ab1,Dr} that for a generic set $ \mathcal{J}_{reg}(J^+,J^-;L_\Lambda,L_{\widehat{\Lambda}}) \subset \mathcal{J}(J^+,J^-)$, all the pseudoholomorphic strips in $\mathcal{M}^k(c^+,c^-;J_W;L_\Lambda,L_{\widehat{\Lambda}})$, for all choices of  $c^+ \in \mathcal{T}_{\Lambda \to \widehat{\Lambda}}(\lambda^+)$ and $c^- \in \mathcal{T}_{\Lambda \to \widehat{\Lambda}}(\lambda^-)$, are Fredholm regular. We assume from now on that $J_W \in \mathcal{J}_{reg}(J^+,J^-;L_\Lambda,L_{\widehat{\Lambda}})$.

Initially, the map  $\Phi_{W,L_\Lambda,{L}_{\widehat{\Lambda}}}$ is obtained by counting elements of $\mathcal{M}^0(c^+,c^-;J_W;L_\Lambda,L_{\widehat{\Lambda}})$, and is defined from $LCH_{st}(\lambda^+,\Lambda \to \widehat{\Lambda})$ to $LCH_{st}(\lambda^-,\Lambda \to \widehat{\Lambda})$. More precisely
\begin{equation}
\Phi_{W,L_\Lambda,{L}_{\widehat{\Lambda}}}(c^+) = \sum_{c^- \in \mathcal{T}_{\Lambda^- \to \widehat{\Lambda}^-}(\lambda^-)} \bigg( \#(\mathcal{M}^0(c^+,c^-;J_W;L_\Lambda,L_{\widehat{\Lambda}})\mod 2 \bigg)c^-.
\end{equation}
To see that it actually induces a map on the homology level one has to check that $d_{J^-} \circ \Phi_{W,L_\Lambda,{L}_{\widehat{\Lambda}}}=\Phi_{W,L_\Lambda,{L}_{\widehat{\Lambda}}} \circ d_{J^+}$. The proof of this fact is a combination of compactness and gluing, and is identical to similar statements for Lagrangian Floer homology or cylindrical contact homology; see for example \cite{B,MS}. A compactness argument identical to the one in Lemma \ref{lemma1} shows that $\Phi_{W,L_\Lambda,{L}_{\widehat{\Lambda}}}$ is well defined and is a finite sum.

Because of the regularity for all pseudoholomorphic curves involved in the definitions of the maps $d_{J^+}$, $d_{J^-}$ and $\Phi_{W,L_\Lambda,{L}_{\widehat{\Lambda}}}$, it is possible to perform gluing for the pseudoholomorphic curves involved in these definitions.
More precisely, the map $d_{J^-} \circ\Phi_{W,L_\Lambda,{L}_{\widehat{\Lambda}}}(c^+)$ counts the number of 2-level pseudoholomorphic buildings $\widetilde{w}$ whose levels $(\widetilde{w}^1_-, \widetilde{w}^2_-)$ satisfy $\widetilde{w}^1_- \in \mathcal{M}^0(c^+,\check{c}^-;J_W;L_\Lambda,L_{\widehat{\Lambda}})$ and $\widetilde{w}^2_- \in \mathcal{M}^1(\check{c}^-,c^-;J^-,\Lambda,\widehat{\Lambda})$ for some $\check{c}^-,c^- \in \mathcal{T}_{\Lambda \to \widehat{\Lambda}}(\lambda^-)$. Analogously, $\Phi_{W,L_\Lambda,{L}_{\widehat{\Lambda}}} \circ d_{J^+}(c^+)$ counts the number of 2-level pseudoholomorphic buildings $\widetilde{w}$  whose levels $(\widetilde{w}^1_+, \widetilde{w}^2_+)$ satisfy $\widetilde{w}^1_+ \in \mathcal{M}^0(c^+,\check{c}^+;J^+;\Lambda,\widehat{\Lambda})$ and $\widetilde{w}^2_+ \in \mathcal{M}^1(\check{c}^+,c^-;J_W;L_\Lambda,L_{\widehat{\Lambda}})$ for some $\check{c}^+ \in \mathcal{T}_{\Lambda \to \widehat{\Lambda}}(\lambda^+)$ and $c^- \in \mathcal{T}_{\Lambda \to \widehat{\Lambda}}(\lambda^-)$. The gluing theorem implies that each 2-level building that appears in the definition of the map $d_{J^-} \circ \Phi_{W,L_\Lambda,{L}_{\widehat{\Lambda}}}= \Phi_{W,L_\Lambda,{L}_{\widehat{\Lambda}}} \circ d_{J^+}$ is one point in the boundary of the moduli space $\overline{\mathcal{M}}^1(c^+,c^-;J_W;L_\Lambda,L_{\widehat{\Lambda}})$.

Because of the exactness of $W$, of the Lagrangians $L_\Lambda$ and $L_{\widehat{\Lambda}}$, and using that $\lambda^+$ and $\lambda^-$ are adapted to the pair of Legendrians $(\Lambda,\widehat{\Lambda})$ we have that a sequence of elements in  $\mathcal{M}^1(c^+,c^-;J_W;L_\Lambda,L_{\widehat{\Lambda}})$ can only break in 2-level buildings that appear in the definition of the maps $d_{J^-} \circ  \Phi_{W,L_\Lambda,{L}_{\widehat{\Lambda}}}$ and $ \Phi_{W,L_\Lambda,{L}_{\widehat{\Lambda}}} \circ d_{J^+}$. The complete proof of this compactness fact uses just the additivity of the Fredholm index for pseudoholomorphic buildings and the Fredholm regularity of the pseudoholomorphic buildings involved, and is identical to a similar argument we used in Lemma \ref{lemma2}.

This combination of compactness and gluing implies that $\overline{\mathcal{M}}^1(c^+,c^-;J_W;L_\Lambda,L_{\widehat{\Lambda}})$ is a one-dimensional manifold whose boundary is composed exactly of the 2-level buildings that appear in the definition of the map $d_{J^-} \circ  \Phi_{W,L_\Lambda,{L}_{\widehat{\Lambda}}} -  \Phi_{W,L_\Lambda,{L}_{\widehat{\Lambda}}} \circ d_{J^+}$.

Summarising, we obtain for each $c^+\in \mathcal{T}_{\Lambda^+ \to \widehat{\Lambda}^+}(\lambda^+)$:
\begin{equation}
(d_{J^+} \circ \Phi_{W,L_\Lambda,{L}_{\widehat{\Lambda}}} -  \Phi_{W,L_\Lambda,{L}_{\widehat{\Lambda}}} \circ d_{J^-})(c^+) = \sum_{c^- \in \mathcal{T}_{\Lambda^- \to \widehat{\Lambda}^-}(\lambda^-) } [a(c^+,c^-)\mod2] c^-
\end{equation}
where $a(c^+,c^-)$ is the number of pseudoholomorphic buildings appearing on the boundary of $\overline{\mathcal{M}}^1(c^+,c^-;J_W;L_\Lambda,L_{\widehat{\Lambda}})$. As $\overline{\mathcal{M}}^1(c^+,c^-;J_W;L_\Lambda,L_{\widehat{\Lambda}})$ is a 1-dimensional manifold, the number of its boundary components is even which implies that $[a(c^-)\mod2]=0$, and finishes the proof of the proposition.

The last assertion of the proposition is a direct consequence of  the fact that if the moduli space $\mathcal{M}^0(c^+,c^-;J_W;L_\Lambda,L_{\widehat{\Lambda}})$ is non-empty then $c^+$ and $c^-$ must belong to the same homotopy class in $\Sigma_{\Lambda \to \widehat{\Lambda}}$.
\qed

In the case where $W$ is the symplectization of a contact manifold with an $\mathbb{R}$-invariant regular almost complex structure the induced cobordism map is the identity. The reason for this is that in this situation the only pseudoholomorphic strips with Fredholm index $0$ are the trivial strips over Reeb chords.

\subsubsection{A special type of cobordisms}

We now restrict our attention to a special situation. We assume let will study cobordisms maps induced by straight exact symplectic cobordism from $\lambda_0$ to $\mathrm{k}\lambda_0$. We now list the properties we impose on our setup.
\begin{itemize}
\item  we let $\lambda_0$ be a contact form on $(Y,\xi)$ adapted to a pair of disjoint connected Legendrian submanifolds $\Lambda$ and $\widehat{\Lambda}$,
\item we let $W$ be a straight exact symplectic cobordism from $\lambda_0$ to  $\mathrm{k}\lambda_0$ where  $0<\mathrm{k}<1$,
\item  we consider the conical exact Lagrangian cobordisms $L_\Lambda$ and $L_{\widehat{\Lambda}}$ in $W$ from $\Lambda$ to itself and from $\widehat{\Lambda}$ to itself, respectively.
\end{itemize}
Notice that it follows that $L_\Lambda$ and $L_{\widehat{\Lambda}}$ are disjoint.
It is well known that $W$ can deformed continuously to the symplectization of $\lambda_0$ though straight exact symplectic cobordisms; see \cite{HMS} for a proof. More precisely we have
\begin{itemize}
\item  there exists a homotopy $(W_{t})_{t\in [0,1]}$ of straight exact symplectic cobordisms from $\lambda_0$ to $\mathrm{k}\lambda_0$ such that $W_{0}=W$ and $W_{1}$ is the symplectization of $\lambda_0$.
\end{itemize}
It is clear that
\begin{itemize}
\item  for every $t\in [0,1]$ we have that $L_\Lambda$ and $L_{\widehat{\Lambda}}$ are conical exact Lagrangian cobordisms in $W$ from $\Lambda$ to itself and from $\widehat{\Lambda}$ to itself, respectively.
\end{itemize}

Like in Proposition \ref{propcobordism} our cobordism will actually induce for each $\rho \in \Sigma_{\Lambda \to \widehat{\Lambda}}$, a map $\Phi_{W,L_\Lambda,{L}_{\widehat{\Lambda}}} $ from $LC\mathbb{H}^{\rho}_{st}(\lambda,\Lambda \to \widehat{\Lambda})$ to $LC\mathbb{H}^{\rho}_{st}(\lambda,\Lambda \to \widehat{\Lambda})$. In order to define the cobordism maps for the subcomplex $LC\mathbb{H}^{\rho}_{st}(\lambda^+,\Lambda \to \widehat{\Lambda})$ the assumptions on the regularity of the almost complex structure $J_W$ are slightly weaker. We will now explain this.

Fix $\rho \in \Sigma_{\Lambda \to \widehat{\Lambda}}$ and $J\in \mathcal{J}(\lambda_0)$. For each pair $c$ and $c'$ of Reeb chords in $\mathcal{T}^{\rho}_{\Lambda \to \widehat{\Lambda}}(\lambda_0)$ and almost complex structure $J_W \in \mathcal{J}(J,J)$, we consider the moduli space $\mathcal{M}^k(c,c';J_W; L_\Lambda,L_{\widehat{\Lambda}})$ of pseudoholomorphic strips with Fredholm index $k$. As all the strips in these moduli spaces are somewhere injective curves, one combines \cite[Theorem 2.14]{CRGG} and the techniques of \cite{Ab1,Dr} to conclude that there exists a dense set $ \mathcal{J}_{reg,\rho}(J,J) \subset \mathcal{J}(J,J)$, such that if $J_W \in \mathcal{J}_{reg,\rho}(J,J)$ all the pseudoholomorphic strips in all the moduli spaces $\mathcal{M}^k(c,c';J_W; L_\Lambda,L_{\widehat{\Lambda}})$, for every $c$ and $c'$, are Fredholm regular. Applying the same reasoning as the one in the proof of Proposition \ref{propcobordism} one shows that for $J_W \in \mathcal{J}_{reg,\rho}(J,J)$ there exists a map $\Phi_{V,\varpi,L,\widehat{L}}$ from $LC\mathbb{H}^{\rho}_{st}(\lambda,\Lambda \to \widehat{\Lambda})$ to $LC\mathbb{H}^{\rho}_{st}(\lambda,\Lambda \to \widehat{\Lambda})$. Notice that $ \mathcal{J}_{reg,\rho}(J,J)$ might contain elements that are not in $ \mathcal{J}_{reg}(J,J)$.

We define for a fixed choice of $J\in \mathcal{J}_{reg}(\lambda_0)$ and $J_W \in \mathcal{J}_{reg,\rho}(J,J)$ the space $\widetilde{\mathcal{J}}(J,J)$ of smooth homotopies
\begin{equation}
J_t \in \mathcal{J}(J,J),  \ t \in [0,1], 
\end{equation}
of almost complex structures such that $J_0=J_W$, $J_1 = J$, and $J_t$ is compatible with the symplectic form on $W_t$ for every $t\in [0,1]$. For Reeb chords $c,c' \in \mathcal{T}_{\Lambda \to \widehat{\Lambda}}(\lambda_0)$ we define the moduli space
\begin{equation}
\widehat{\mathcal{M}}^k(c,c';J_t;L_\Lambda,L_{\widehat{\Lambda}}):= \{(t,\widetilde{w}) | t \in [0,1] \mbox{ and } \widetilde{w}\in  \mathcal{M}^k(c,c';J_t;L_\Lambda,L_{\widehat{\Lambda}})\}.
\end{equation}

By the results in \cite{Ab1,CRGG,Dr}, we know that there is a generic subset $\widetilde{\mathcal{J}}_{reg}(J,J) \subset \widetilde{\mathcal{J}}(J,J)$  such that all elements of $\widehat{\mathcal{M}}^k(c,c';J_t;L_\Lambda,L_{\widehat{\Lambda}})$ are Fredholm regular for every choice of $c,c' \in \mathcal{T}^{\rho}_{\Lambda \to \widehat{\Lambda}}(\lambda_0)$ and every number $k$. This implies that $\widehat{\mathcal{M}}^k(c,c';J_t;L_\Lambda,L_{\widehat{\Lambda}})$ is a $k+1$-dimensional manifold with boundary. The crucial condition that makes this valid is again the fact that all the pseudoholomorphic curves that make part of this moduli space are somewhere injective. Notice that $\widehat{\mathcal{M}}^k(c,c';J_t,L_\Lambda,L_{\widehat{\Lambda}})$ is not necessarily compact since sequences of elements in $\widehat{\mathcal{M}}^k(c,c';J_t;L_\Lambda,L_{\widehat{\Lambda}})$ might break. However, as shown in \cite{CPT}, it can be compactfied with the addition of pseudoholomorphic buildings.

\begin{prop} \label{propchainhom}
Let $(Y,\xi)$ be a contact manifold and $\lambda_0$ be a contact form on $(Y,\xi)$ adapted to a pair of disjoint connected Legendrian submanifolds $\Lambda$ and $\widehat{\Lambda}$ on $(Y,\xi)$.
Let $W$ be a straight exact symplectic cobordism from $\lambda_0$ to $\mathrm{k}\lambda_0$ where $0<\mathrm{k}<1$ is a constant, and let $L_\Lambda$ and $L_{\widehat{\Lambda}}$ be the conical exact Lagrangian cobordisms in $W$ from $\Lambda$ to itself and from $\widehat{\Lambda}$ to itself, respectively. 
Then, if $J \in \mathcal{J}_{reg}(\lambda_0)$, we have that for all $J_W \in \mathcal{J}_{reg,\rho}(J,J)$ the map $\Phi_{W,L_\Lambda,{L}_{\widehat{\Lambda}}} $ from $LC\mathbb{H}^{\rho}_{st}(\lambda,\Lambda \to \widehat{\Lambda})$ to itself is chain homotopic to the identity.
\end{prop}

\textit{ Proof:} The proof is a standard argument in SFT, and we direct the reader to the original source \cite{SFT} for an exposition of this argument for general SFT invariants and \cite{B,Bt} where the very similar case of cylindrical contact homology is treated. We first take an almost complex structure $J \in \mathcal{J}_{reg}(\lambda_0) $ and choose an almost complex structure $J_W \in \mathcal{J}_{reg}(J,J)$ compatible with $\varpi$. The map $\Phi_{W,L_\Lambda,{L}_{\widehat{\Lambda}}} $ will count pseudoholomorphic strips in $(W,J_W)$ satisfying boundary conditions as the ones used in Proposition \ref{propcobordism}.
From our discussion above we know that there exists an homotopy  $(W_{t})_{t\in [0,1]}$ straight exact symplectic cobordisms given from $W$ to the symplectization of $\lambda_0$, and that $L_\Lambda$ and $L_{\widehat{\Lambda}}$ can be considered as conical exact symplectic cobordims on $W_t$ for each $t\in [0,1]$.
 For the homotopy $W_t$ we can take a homotopy $J_t$ of almost complex structures in $\widetilde{\mathcal{J}}_{reg}(J,J)$. 

The most important ingredient of the proof is the description of the compactification  $\overline{\widehat{\mathcal{M}}}^0(c,c';J_t,L_\Lambda,L_{\widehat{\Lambda}}))$ of the moduli space $\widehat{\mathcal{M}}^0(c,c';J_t,L_\Lambda,L_{\widehat{\Lambda}}))$. We want to understand the boundary of $\overline{\mathcal{M}}^0(c,c';J_t,L_\Lambda,L_{\widehat{\Lambda}}))$. For this we take a sequence $(t_n,\widetilde{w}_n)$ of elements in $\widehat{\mathcal{M}}^0(c,c';J_t,L_\Lambda,L_{\widehat{\Lambda}}))$ converging to the boundary of $\overline{\widehat{\mathcal{M}}}^0(c,c';J_t,L_\Lambda,L_{\widehat{\Lambda}}))$. There are three possibilities: $t_n$ goes to $0$, $t_n$ goes to $1$ or the limit $t_{\infty}$ of $t_n$ belongs to $(0,1)$. In all three possibilities we know that the conditions (a), (b), (c) and (d) satisfied by $\lambda_0$, $\Lambda$ and $\widehat{\Lambda}$ prevent any ``bubbling'' in the sequence. If $t_n \to 0$ regularity of $J_W$ and $J$ imply that the $\widetilde{w}_n$ must converge to an element $\widetilde{w}$ of $\mathcal{M}^0(c,c';J_W;L_\Lambda,L_{\widehat{\Lambda}}))$. If $t_n \to 1$ it follows from the regularity of $J$ that $\widetilde{w}_n$ must  converge to $\widetilde{w}$ of $\mathcal{M}^0(c,c';J,\Lambda,\widehat{\Lambda})$. If $t_n \to t_{\infty} \in (0,1)$ the sequence $\widetilde{w}_n$ converges to a pseudoholomorphic building $\widetilde{w}$ whose levels $\widetilde{w}^l$ are all pseudoholomorphic strips. The fact that for $J_{t_{\infty}}$ all the strips appearing have Fredholm index $\geq -1$ gives us more information about this pseudoholomorphic building. Because $(t_n,\widetilde{w})$ is in the boundary of $\overline{\widehat{\mathcal{M}}}^0(c,c';J_t;L_\Lambda,L_{\widehat{\Lambda}}))$ and $t_n \in (0,1)$, we know that $\widetilde{w}$ has at least 2-levels. It also has Fredholm index $0$ since all $\widetilde{w}_n$ have Fredholm index $0$. Of all levels of $\widetilde{w}$ there is precisely one in the cobordism $W_{t_\infty}$, and all the others are in symplectizations of $\lambda_0$ above and below the cobordism. The levels in the symplectizations have Fredholm index $\geq 1$ and the one on $
(W_{t_\infty},J_{t_{\infty}})$ has Fredholm index $\geq -1$. Combining this with the fact that $\widetilde{w}$ has Fredholm index $0$ we conclude that $\widetilde{w}$ has exactly 2 levels: one with Fredholm index $-1$ in $W_{t_\infty}$ and one with Fredholm index $1$ in the symplectization of $\lambda_0$. There are then two possible cases: either the building $\widetilde{w}$ has the top level $\widetilde{w}^1 \in \widetilde{\mathcal{M}}^1(c,\check{c};J;\Lambda,\widehat{\Lambda})$ and the lower level  $\widetilde{w}^2 \in \mathcal{M}^{-1}(\check{c},c';J_{t_\infty},L_\Lambda,L_{\widehat{\Lambda}}))$ for some $\check{c} \in \mathcal{T}^\rho_{\Lambda \to \widehat{\Lambda}}(\lambda_0)$; or the building $\widetilde{w}$ has the top level $\widetilde{w}^1 \in \mathcal{M}^{-1}(c,\check{c};J_{t_\infty};L_\Lambda,L_{\widehat{\Lambda}}))$ and the lower level  $\widetilde{w}^2 \in \widetilde{\mathcal{M}}^1(\check{c},c';J;\Lambda,\widehat{\Lambda})$ for some $\check{c} \in \mathcal{T}^\rho_{\Lambda \to \widehat{\Lambda}}(\lambda_0)$.

On the other hand, the gluing theorem implies that: if $\widetilde{w}$ is a 2-level building whose top level $\widetilde{w}^1 \in \widetilde{\mathcal{M}}^1(c,\check{c};J;\Lambda,\widehat{\Lambda})$ and whose lower level  $\widetilde{w}^2 \in \mathcal{M}^{-1}(\check{c},c';J_{t_\infty};L_\Lambda,L_{\widehat{\Lambda}}))$, then $(t_\infty,\widetilde{w})$ is an element of the boundary of $\overline{\widehat{\mathcal{M}}}^0(c,c';J_t;L_\Lambda,L_{\widehat{\Lambda}}))$. The same is valid for $(t_\infty,\widetilde{w})$ if the building $\widetilde{w}$ has the top level $\widetilde{w}^1 \in \mathcal{M}^{-1}(c,\check{c};J_{t_\infty},L_\Lambda,L_{\widehat{\Lambda}}))$ and the lower level  $\widetilde{w}^2 \in \widetilde{\mathcal{M}}^1(\check{c},c';J;L_\Lambda,L_{\widehat{\Lambda}}))$. Finally, if $\widetilde{w} \in \mathcal{M}^0(c,c';J_V;L_\Lambda,L_{\widehat{\Lambda}}))$ then the regularity of the homotopy $J_t$ implies that $(0,\widetilde{w})$ is in the boundary of $\overline{\widehat{\mathcal{M}}}^0(c,c';J_t;L_\Lambda,L_{\widehat{\Lambda}}))$. The same reasoning implies that  $(1,\widetilde{w})$ is in the boundary of $\overline{\widehat{\mathcal{M}}}^0(c,c';J_t;L_\Lambda,L_{\widehat{\Lambda}}))$ if $\widetilde{w} \in \widetilde{\mathcal{M}}^0(c,c';J;\Lambda,\widehat{\Lambda})$.
We have thus obtained a complete description of the boundary of $\overline{\widehat{\mathcal{M}}}^0(c,c';J_t;L_\Lambda,L_{\widehat{\Lambda}}))$.

We now define a map $\mathcal{K}: LCH_{st}^{\rho}(\lambda_0,\Lambda \to \widehat{\Lambda}) \to LCH_{st}^{\rho}(\lambda_0,\Lambda \to \widehat{\Lambda})$ that counts finite energy Fredholm index $-1$ pseudoholomorphic strips in the homotopy $(W_t,J_t)$, with one boundary component in $L_\Lambda$ and one in $L_{\widehat{\Lambda}}$. In order to define $\mathcal{K}$ precisely we first define for $c,c' \in \mathcal{T}^\rho_{\Lambda \to \widehat{\Lambda}}(\lambda_0)$ the set $Q_{c,c'} := \{t\in [0,1] | \mathcal{M}^{-1}(c,c';J_{t_\infty};L_\Lambda,L_{\widehat{\Lambda}})) \neq \emptyset \}$. From the regularity of the homotopy $J_t$ we know that $Q_{c,c'}$ is always a finite set. We can now define for $c \in \mathcal{T}^\rho_{\Lambda \to \widehat{\Lambda}}(\lambda_0)$
\begin{equation}\label{mapK}
\mathcal{K}(c) := \sum_{c' \in \mathcal{T}^\rho_{\Lambda \to \widehat{\Lambda}}(\lambda_0)} ((\sum_{t\in Q_{c,c'}} \# \mathcal{M}^{-1}(c,c';J_{t_\infty},L_{t_\infty},\widehat{L}_{t_\infty}))\mod 2 ) \ c'.
\end{equation}

From the regularity of the homotopy $J_t$, we know that $\mathcal{M}^{-1}(c,c';J_{t_\infty};L_\Lambda,L_{\widehat{\Lambda}}))$ is finite for every $t \in [0,1]$. For action reasons, we also know that for a fixed $c \in \mathcal{T}^\rho_{\Lambda \to \widehat{\Lambda}}(\lambda_0)$ the term $a_{\mathcal{K}}(c,c'):= \sum_{t\in Q_{c,c'}} \# \mathcal{M}^{-1}(c,c';J_{t_\infty};L_\Lambda,L_{\widehat{\Lambda}}))$ can only be non-zero for a finite number of chords $c' \in \mathcal{T}^\rho_{\Lambda \to \widehat{\Lambda}}(\lambda_0)$. These considerations imply that $\mathcal{K}(c)$ as defined above is well-defined and is a finite sum. We then obtain $\mathcal{K}$ by extending it linearly to $LCH_{st}^{\rho}(\lambda_0,\Lambda \to \widehat{\Lambda})$.

Consider the map $Id + \Phi_{W,L_\Lambda,{L}_{\widehat{\Lambda}}}  + \mathcal{K}\circ d_J + d_J \circ \mathcal{K}$ from $LCH_{st}^{\rho}(\lambda_0,\Lambda \to \widehat{\Lambda})$ to itself. It follows from our discussion so far, that the pseudoholomorphic curves which are counted in the definition of $Id +\Phi_{W,L_\Lambda,{L}_{\widehat{\Lambda}}}  + \mathcal{K}\circ d_J + d_J \circ \mathcal{K}(c)$, are exactly the ones that appear in the boundary of the moduli spaces $\overline{\widehat{\mathcal{M}}}^0(c,c';J_t;L_\Lambda,L_{\widehat{\Lambda}}))$. As the compactified moduli space $\overline{\widehat{\mathcal{M}}}^0(c,c';J_t;L_\Lambda,L_{\widehat{\Lambda}}))$ is a finite union of compact intervals it has an even number of boundary components. It then follows that $Id + \Phi_{W,L_\Lambda,{L}_{\widehat{\Lambda}}}  + \mathcal{K}\circ d_J + d_J \circ \mathcal{K}(c)=0$.  As this is valid for all Reeb chords $c \in \mathcal{T}^\rho_{\Lambda \to \widehat{\Lambda}}(\lambda_0)$ we have proved that $\Phi_{W,L_\Lambda,{L}_{\widehat{\Lambda}}} $ is chain-homotopic to the identity, as claimed.
\qed

Combining Proposition \ref{propchainhom} with a gluing and compactness argument which is standard in SFT  \cite{SFT} one can prove that the strip Legendrian contact homology $LC\mathbb{H}^{\rho}_{st}(\lambda,\Lambda \to \widehat{\Lambda})$ does not depend on the regular almost complex structure $J\in  \mathcal{J}_{reg}(\lambda_0)$ used to define $d_J$, something which is not obvious from the definition of $LC\mathbb{H}^{\rho}_{st}(\lambda_0,\Lambda \to \widehat{\Lambda})$.

\section{Growth rate of $LC\mathbb{H}_{st}(\lambda_0,\Lambda \to \widehat{\Lambda})$ and lower bounds for $h_{top}$} \label{sectiongrowth}
Let $(Y,\xi)$ be a contact manifold. Given a contact form $\lambda_0$ on $(Y,\xi)$ adapted to a pair of disjoint connected Legendrian submanifolds $\Lambda$ and $\widehat{\Lambda}$ we define the exponential homotopical growth of the strip Legendrian contact homology $LC\mathbb{H}_{st}(\lambda_0,\Lambda \to \widehat{\Lambda})$ with respect to the action. We then use it to estimate the growth of the number of Reeb chords from $\Lambda$ to $\Lambda'$ of other contact forms on $(Y,\xi)$, for a Legendrian submanifold $\Lambda'$ sufficiently close to $\widehat{\Lambda}$.

We start defining for each number $C > 0$ the subset $\Sigma^{C}_{\Lambda \to \widehat{\Lambda}}(\lambda_0)$ of $\Sigma_{\Lambda \to \widehat{\Lambda}}$ of homotopy classes $\rho$ satisfying:
\begin{itemize}
\item{all the chords in $\mathcal{T}^{\rho}_{\Lambda \to \widehat{\Lambda}}(\lambda_0)$ have action $\leq C$,}
\item{$LC\mathbb{H}^{\rho}_{st}(\lambda_0,\Lambda \to \widehat{\Lambda})\neq 0$.}
\end{itemize}

\begin{defi} \label{defigrowth}
Let $\lambda_0$ be a contact form on $(Y,\xi)$ adapted to a pair of disjoint Legendrian knots $\Lambda$ and $\widehat{\Lambda}$. We say that $LC\mathbb{H}_{st}(\lambda_0,\Lambda \to \widehat{\Lambda})$ has exponential homotopical growth with exponential weight $a>0$ if there exist real numbers $d$ and $C_0\geq 0$ such that:
\begin{equation}
\#(\Sigma^{C}_{\Lambda \to \widehat{\Lambda}}(\lambda_0)) > e^{aC+d}
\end{equation}
for all $C\geq C_0$.
\end{defi}

If $\Lambda$ and $\Lambda'$ are Legendrian submanifolds of a contact manifold $(Y,\xi)$ and $\lambda$ is a contact form on $(Y,\xi)$, we will say that $\Lambda'$ is $(\lambda,\Lambda)$-transverse if all the Reeb chords in $\mathcal{T}_{\Lambda \to \Lambda'}(\lambda)$ are transverse. In this case we will denote by $N_C(\lambda,\Lambda,\Lambda')$ the number of Reeb chords in $\mathcal{T}_{\Lambda \to \Lambda'}(\lambda)$ with action $\leq C$. 

\subsection{The growth of the number of Reeb chords}

We start with some auxiliary results that we will need. We fix $(Y,\xi)$ be a contact manifold $\Lambda$ and $\widehat{\Lambda}$ be disjoint connected Legendrian submanifolds in $(Y,\xi)$, and $\lambda$ be a contact form on $(Y,\xi)$. 

Given $\delta>0$ we say that a Legendrian submanifold $\widehat{\Lambda}'$ in the same Legendrian isotopy class of $\widehat{\Lambda}$ is $\delta$-close to $\widehat{\Lambda}$ in the $C^3$-sense if there exists a Legendrian isotopy $\mathcal{L}: [0,1] \times \widehat{\Lambda} \to Y$ satisfying $\mathcal{L}( \{0\} \times \widehat{\Lambda}) = \widehat{\Lambda}$  and $\mathcal{L}( \{1\} \times \widehat{\Lambda}) = \widehat{\Lambda}'$, and which is $\delta$-close to the constant isotopy stationary at  $\widehat{\Lambda}$ in the $C^3$-topology. Notice that we demand the time derivatives of the isotopy to be $<\delta$.

We consider the space $Diff^1_{diff}(Y)$ of $C^1$ diffeomorphisms of $Y$ endowed with a metric that generates the canonical topology on $Diff^1_{diff}(Y)$.
Given $\epsilon>0$ we say that a contactomorphism $\Upsilon: (Y,\xi) \to (Y,\xi)$ is $\epsilon$-close to the identity in the $C^1$-sense if $\epsilon$-close to the identity in this metric. 
The following lemma is elementary. The proof is essentially contained in that of \cite[Theorem 2.6.2]{Ge} and for the convenience of the reader we present it in  Appendix B.
\begin{lemma} \label{lemma:changeofposition}
Let $\lambda$ be a contact form on $(Y,\xi)$ and $\Lambda$ and $\widehat{\Lambda}$ be disjoint connected Legendrian submanifolds in $(Y,\xi)$. Let $\mathcal{V}(\widehat{\Lambda})$ be a tubular neighbourhood of $\widehat{\Lambda}$ that does not intersect $\Lambda$. Then, given $\epsilon>0$ there exists $\delta>0$ such that for every $\widehat{\Lambda}'$ that is $\delta$-close to $\widehat{\Lambda}$ in the $C^3$-sense there exists a contactomorphism $\Upsilon_{\widehat{\Lambda}'}: (Y,\xi) \to (Y,\xi)$ which satisfies
\begin{itemize}
\item [(1)] $\Upsilon_{\widehat{\Lambda}'}(\widehat{\Lambda}) = \widehat{\Lambda}'$,
\item [(2)] $\Upsilon_{\widehat{\Lambda}'}$ is $\epsilon$-close to the identity in the $C^1$-sense,
\item [(3)] $\Upsilon_{\widehat{\Lambda}'}$ coincides with the identity in the complement of $\mathcal{V}(\widehat{\Lambda})$.
\end{itemize}
\end{lemma}

We are ready to prove
\begin{prop}\label{mainproposition}
Let $(Y,\xi)$ be a contact  manifold and $\Lambda$ and $\widehat{\Lambda}$ be two disjoint connected Legendrian submanifolds on $(Y,\xi)$. Suppose that $\lambda_0$ is a contact form on $(Y,\xi)$ adapted to the pair $(\Lambda,\widehat{\Lambda})$,
and that the strip contact homology $LC\mathbb{H}_{st}(\lambda_0,\Lambda \to \widehat{\Lambda})$ has exponential homotopical growth with exponential weight $a>0$.  Let $\lambda$ be another contact form associated to $(Y,\xi)$, and denote by $f_{\lambda}$ the function such that $\lambda=f_{\lambda}\lambda_0$. Then, given $\mu>0$ there exists $\delta>0$ such that for any every Legendrian submanifold $\widehat{\Lambda}'$ which is $(\lambda,\Lambda)$ transverse and $\delta$-close to $\widehat{\Lambda}$ in the $C^3$ sense we have that the numbers $N_C(\lambda,\Lambda,\widehat{\Lambda}')$ satisfy
\begin{equation}
e^{\frac{aC}{\max(f_\lambda) + \mu}+d}<N_C(\lambda,\Lambda,\widehat{\Lambda}')
\end{equation}
for all $C\geq C_0$, where $C_0$ is the constant in Definition \ref{defigrowth}.
\end{prop}
\textit{Proof:} We divide the proof in steps.

\textbf{Step 1:} We fix $\mu>0$. Then there exists a straight exact symplectic cobordism from $(\max(f_\lambda) + \mu) \lambda_0$ to $\lambda$. 
Therefore, there exists an $\epsilon>0$ such that if $\Upsilon: (Y,\xi) \to (Y,\xi)$ is a contactomorphism which is $\epsilon$-close to the identity in the $C^1$ sense then there exists a straight exact symplectic cobordism from $(\max(f_\lambda) + \mu) \lambda_0$ to the pull back contact form  $\Upsilon^*\lambda$. We fix $\epsilon>0$ with this property.

\textbf{Step 2:}
Given $\epsilon>0$ from step 1 we choose $\delta>0$ as in the statement of Lemma \ref{lemma:changeofposition}.
Let $\widehat{\Lambda}'$ be $\delta$-close to $\widehat{\Lambda}$ in the $C^3$-sense: we know that there exists a contactmorphism $\Upsilon_{\widehat{\Lambda}'}: (Y,\xi) \to (Y,\xi)$ satisfying (1), (2) and (3) in the statement of Lemma \ref{lemma:changeofposition}.
Let $\lambda_{\widehat{\Lambda}'}= \Upsilon_{\widehat{\Lambda}'}^* \lambda$. A direct computation shows that the Reeb flows of $\lambda_{\widehat{\Lambda}'}$ and $\lambda$ are conjugate by $\Upsilon_{\widehat{\Lambda}'}^{-1}$, i.e
\begin{equation} \label{eqconjug}
\Upsilon_{\widehat{\Lambda}'}^{-1} \circ \phi^t_{X_\lambda}= \phi^t_{X_{\lambda_{\widehat{\Lambda}'}}}\circ \Upsilon_{\widehat{\Lambda}'}^{-1}.
\end{equation}
By (3) we know that $\Upsilon_{\widehat{\Lambda}}(\Lambda)= \Lambda$. Combining this with condition (1) of Lemma \ref{lemma:changeofposition} and \eqref{eqconjug}  implies that $\Upsilon_{\widehat{\Lambda}'}^{-1}$ takes Reeb chords of $\lambda$ from $\Lambda$ to $\widehat{\Lambda}'$ to Reeb chords of $\lambda_{\widehat{\Lambda}'}$ from $\Lambda$ to $\widehat{\Lambda}$.
We thus have a bijection between the sets $\mathcal{T}_{\Lambda \to \widehat{\Lambda}'}(\lambda)$ and $\mathcal{T}_{\Lambda \to \widehat{\Lambda}}(\lambda_{\widehat{\Lambda}'})$. This bijection clearly preserves the action of the Reeb chords. We therefore obtain
\begin{equation} \label{eqbijec}
N_C(\lambda,\Lambda,\widehat{\Lambda}')= N_C(\lambda_{\widehat{\Lambda}'},\Lambda,\widehat{\Lambda}).
\end{equation}

\textbf{Step 3:}
Because of \eqref{eqbijec} we know that to prove the proposition it is sufficient to estimate $N_C(\lambda_{\widehat{\Lambda}'},\Lambda,\widehat{\Lambda})$. For this we consider the function $f_{\lambda_{\widehat{\Lambda}'}}$ such that $f_{\lambda_{\widehat{\Lambda}'}}\lambda_0=\lambda_{\widehat{\Lambda}'}$, and let $\mathrm{k}_{\lambda_{\widehat{\Lambda}'}}= \frac{\min f_{\lambda_{\widehat{\Lambda}'}}}{2}$. It is clear that there exists a straight exact symplectic cobordism $W(\lambda_{\widehat{\Lambda}'},\mathrm{k}_{\lambda_{\widehat{\Lambda}'}}\lambda_0)$ from $\lambda_{\widehat{\Lambda}'}$ to $\mathrm{k}_{\lambda_{\widehat{\Lambda}'}}\lambda_0$. By step 1 and our choice of $\delta$ in step 2 we know that there also exists a straight exact symplectic cobordism $W(\lambda_0,\lambda_{\widehat{\Lambda}'})$ from $(\max(f_\lambda) + \mu) \lambda_0$ to $\lambda_{\widehat{\Lambda}'}$. 

Using the existence of $W(\lambda_{\widehat{\Lambda}'},\mathrm{k}_{\lambda_{\widehat{\Lambda}'}}\lambda_0)$ and $W(\lambda_0,\lambda_{\widehat{\Lambda}'})$  we can apply the construction in section \ref{splitting} to obtain a splitting family $(W(\lambda_{\widehat{\Lambda}'})_R)_{R\in (0,+\infty)}$ of straight exact symplectic cobordisms from $\lambda_0$ to $\mathrm{k}_{\lambda_{\widehat{\Lambda}'}}\lambda_0$ along $\lambda_{\widehat{\Lambda}'}$.

\textbf{Step 4.}

On each $W(\lambda_{\widehat{\Lambda}'})_R$ we consider the conical exact Lagrangian cobordisms $L_\Lambda$ and $L_{\widehat{\Lambda}}$. The triple $(W(\lambda_{\widehat{\Lambda}'})_R, L_\Lambda,L_{\widehat{\Lambda}})$ satisfies the hypothesis of Proposition \ref{propchainhom}. 

Let $J\in \mathcal{J}_{reg}(\lambda)$ and $\rho\in \Sigma_{\Lambda \to \widehat{\Lambda}}$.  
 Applying Proposition \ref{propchainhom} we obtain that for each almost complex structure $\overline{J}_R \in \mathcal{J}_{reg,\rho}(J,J)$ in $W(\lambda_{\widehat{\Lambda}'})_R$ the map $\Phi_{W(\lambda_{\widehat{\Lambda}'})_R,L_\Lambda,L_{\widehat{\Lambda}}}$ from $LC\mathbb{H}^{\rho}_{st}(\lambda,\Lambda \to \widehat{\Lambda})$ to itself is the identity.
 
 \textbf{Step 5.} We prove the proposition in the case $\lambda$ is non-degenerate and that all Reeb chords in $\mathcal{T}_{\Lambda}(\lambda)$, $\mathcal{T}_{\widehat{\Lambda}'}(\lambda)$ and $\mathcal{T}_{\Lambda \to \widehat{\Lambda}'}(\lambda)$ are transverse. We remark that this is equivalent to demanding that $\lambda_{\widehat{\Lambda}'}$ is non-degenerate and that all Reeb chords in $\mathcal{T}_{\Lambda}(\lambda_{\widehat{\Lambda}'})$, $\mathcal{T}_{\widehat{\Lambda}}(\lambda_{\widehat{\Lambda}'})$ and $\mathcal{T}_{\Lambda \to \widehat{\Lambda}}(\lambda_{\widehat{\Lambda}'})$ are transverse.

We first pick for $W(\lambda_{\widehat{\Lambda}'})_R$ a splitting almost complex structure $J_R$ as in Section \ref{splitting} and take $\rho \in \Sigma^{\frac{C}{\max(f_\lambda) + \mu}}_{\Lambda \to \widehat{\Lambda}}(\lambda_0)$. We claim that for such an almost complex structure, there exist Reeb chords $c,c' \in \mathcal{T}^{\rho}_{\Lambda \to \widehat{\Lambda}}(\lambda_0)$ such that $\mathcal{M}(c,c';J_R; L_\Lambda,L_{\widehat{\Lambda}})$ is non-empty.

We argue by contradiction: if no such strip exists we obtain that $J_R \in \mathcal{J}_{reg,\rho}(J,J)$ and it follows from Step 4 that the map $\Phi_{W(\lambda_{\widehat{\Lambda}'})_R,L_\Lambda,L_{\widehat{\Lambda}}}$ from  $LC\mathbb{H}^{\rho}_{st}(\lambda_0,\Lambda \to \widehat{\Lambda})$ to itself is an isomorphism. Because $\mathcal{M}(c,c';J_R; L_\Lambda,L_{\widehat{\Lambda}})$ is empty for all Reeb chords $c,c' \in \mathcal{T}^{\rho}_{\Lambda \to \widehat{\Lambda}}(\lambda)$, we conclude that the cobordism map $\Phi_{W(\lambda_{\widehat{\Lambda}'})_R,L_\Lambda,L_{\widehat{\Lambda}}}$ is the zero map. On the other hand, as $LC\mathbb{H}^{\rho}_{st}(\lambda_0,\Lambda \to \widehat{\Lambda})\neq 0$ and $\Phi_{W(\lambda_{\widehat{\Lambda}'})_R,L_\Lambda,L_{\widehat{\Lambda}}}$ is an isomorphism it cannot be the zero map. We have thus reached a contradiction.  Therefore that there must exist Reeb chords $c,c' \in \mathcal{T}^{\rho}_{\Lambda \to \widehat{\Lambda}}(\lambda)$ such that  $\mathcal{M}(c,c';J_R; L_\Lambda,L_{\widehat{\Lambda}})$ is non-empty.

We send $R\to +\infty$ and take a sequence $\widetilde{u}_R \in  \mathcal{M}(c,c';J_R; L_\Lambda,L_{\widehat{\Lambda}})$. Because there is a bound on the energy of all elements of  $\mathcal{M}(c,c';J_R; L_\Lambda,L_{\widehat{\Lambda}})$, the SFT-compactness results of \cite{CPT} implies that $\widetilde{u}_R $ has a subsequence that converges to a pseudoholomorphic building $\widetilde{w}$. Notice that to apply the SFT-compactness theorem in this step we need the assumption that  $\lambda_{\widehat{\Lambda}'}$ is non-degenerate and that all Reeb chords in $\mathcal{T}_{\Lambda}(\lambda_{\widehat{\Lambda}'})$, $\mathcal{T}_{\widehat{\Lambda}}(\lambda_{\widehat{\Lambda}'})$ and $\mathcal{T}_{\Lambda \to \widehat{\Lambda}}(\lambda_{\widehat{\Lambda}'})$ are transverse.
 Because of the stretching the neck process, we have that one of the levels of this building lives in an exact symplectic cobordism from $\lambda_{\widehat{\Lambda}'}$ to $\mathrm{k}_{\lambda_{\widehat{\Lambda}'}}\lambda_0$.

We will see that for topological reasons one of the punctures of this level has to detect a Reeb chord $\widehat{c} \in \mathcal{T}^{\rho}_{\Lambda \to \widehat{\Lambda}}(\lambda_{\widehat{\Lambda}'})$ with action smaller than $C$. Let $\widetilde{w}^k$ for $k\in \{1,...,m\}$ be the levels of the pseudoholomorphic building $\widetilde{w}$. Because the topology of our pseudoholomorphic curve does not change on the breaking we must have the following picture.
\begin{itemize}
\item{The upper level $\widetilde{w}^1$ is composed of one pseudoholomorphic disc, with one positive puncture, which is asymptotic to a Reeb chord $c_0  \in \mathcal{T}^{\rho}_{\Lambda \to \widehat{\Lambda}}((\max(f_\lambda) +\mu) \lambda_0)=\mathcal{T}^{\rho}_{\Lambda \to \widehat{\Lambda}}(\lambda_0)$, and several negative boundary and interior punctures. All of the negative punctures detect contractible Reeb orbits,  Reeb chords from $\Lambda$ to itself that represent the trivial element of  $\pi_1(Y,\Lambda)$, or Reeb chords from $\widehat{\Lambda}$ that represent the trivial element of $\pi_1(Y,\widehat{\Lambda})$, with the exception of one negative boundary puncture that detects a Reeb chord $c_1$ which belongs: either to $\mathcal{T}^{\rho}_{\Lambda \to \widehat{\Lambda}}(\lambda_0)$ in case this level lives in the symplectization of $\lambda_0$, or to $\mathcal{T}^{\rho}_{\Lambda \to \widehat{\Lambda}}(\lambda_{\widehat{\Lambda}'})$ in case this level lives in a cobordism from $(\max(f_\lambda) +\mu) \lambda_0$ to $\lambda_{\widehat{\Lambda}'}$.}
\item{On every other level $\widetilde{w}^k$ there is a special curve which has one positive puncture, which is asymptotic to a Reeb chord $c_{k-1}$ from $\Lambda$ to $\widehat{\Lambda}$ in the homotopy class $\rho$ and possibly several interior and boundary negative punctures. Of the negative boundary punctures there is one that is asymptotic to a Reeb orbit $c_k$ from $\Lambda$ to $\widehat{\Lambda}$ in the homotopy class $\rho$, and all other negative punctures detect contractible Reeb orbits,  Reeb chords from $\Lambda$ to itself that represent the trivial element of  $\pi_1(Y,\Lambda)$, or Reeb chords from $\widehat{\Lambda}$ that represent the trivial element of $\pi_1(Y,\widehat{\Lambda})$.}
\end{itemize}
As a consequence we obtain that the level $\widetilde{w}^k$ living in the exact symplectic cobordism from $\lambda_{\widehat{\Lambda}'}$ to $\mathrm{k}_{\lambda_{\widehat{\Lambda}'}}\lambda_0$ contains a pseudoholomorphic curve with one positive puncture asymptotic to a Reeb chord $\widehat{c} \in \mathcal{T}^{\rho}_{\Lambda \to \widehat{\Lambda}}(\lambda_{\widehat{\Lambda}'})$.

In order to obtain the bound on the action of $\widehat{c}$ we first observe that the Reeb chord $c_0$ has action $\leq C$ for the contact form $(\max(f_\lambda) +\mu)\lambda_0$, since its action as a Reeb chord of $\lambda_0$ is $\frac{C}{\max(f_\lambda) + \mu}$. It follows from Stokes theorem that the sum of the action of the positive punctures of a pseudoholomorphic curve with boundary on a conical exact Lagrangian cobordism is bigger or equal to the sum of the action of the negative punctures. Applying this to the building $\widetilde{w}$ we conclude that all Reeb chords and Reeb orbits detected by components of $\widetilde{w}$ have action $\leq C$.
 It follows that the action of $\widehat{c}$  is $ \leq C$.

Applying this reasoning for all  homotopy classes $\rho \in \Sigma^{\frac{C}{\max(f_\lambda) + \mu}}_{\Lambda \to \widehat{\Lambda}}(\lambda_0)$, we obtain that
\begin{equation} \label{eqimp}
 N_C(\lambda_{\widehat{\Lambda}'}) ,\Lambda,\widehat{\Lambda})\geq \# (\Sigma^{\frac{C}{\max(f_\lambda) + \mu}}_{\Lambda \to \widehat{\Lambda}}(\lambda_0)) > e^{\frac{aC}{\max(f_\lambda) + \mu}+d},
 \end{equation}
 for all $C\geq C_0$. The proposition in this case then follows from combining \eqref{eqimp} and \eqref{eqbijec}.
 
 \textbf{Step 6.} Passing to general $\lambda$.
 
 Let $j$ be a natural number. As $\widehat{\Lambda}'$ is $(\lambda,\Lambda)$ transverse, it is possible to make a $C^\infty$ small perturbation of the  contact form $\lambda$ to a non-degenerate contact form $\lambda(j)$ on $(Y,\xi)$ such that:
\begin{itemize}
\item{$N_C(\lambda,\Lambda,\widehat{\Lambda}')=N_C(\lambda(j),\Lambda,\widehat{\Lambda}')$, for all $C\leq j$,}
\item{ all Reeb chords in $\mathcal{T}_{\Lambda}(\lambda)$, $\mathcal{T}_{\widehat{\Lambda}'}(\lambda)$ and $\mathcal{T}_{\Lambda \to \widehat{\Lambda}'}(\lambda)$ are transverse.}
\end{itemize}

If the perturbation $\lambda(j)$ is close enough to $\lambda$ then for each $\widehat{\Lambda}'$ that is $\delta$-close to $\widehat{\Lambda}$ in the $C^3$-sense we know that there exists a straight exact symplectic cobordism from $(\max(f_\lambda) + \mu) \lambda_0$ to $\Upsilon_{\widehat{\Lambda}'}^*\lambda(j)$.
Applying the steps 3, 4 and 5 to $\Upsilon_{\widehat{\Lambda}'}^*\lambda(j)$ we conclude that :
\begin{equation}
N_C(\lambda(j),\Lambda,\widehat{\Lambda}')=N_C(\Upsilon_{\widehat{\Lambda}'}^*\lambda(j),\Lambda,\widehat{\Lambda}) >e^{\frac{aC}{\max(f_\lambda) + \mu }+d},
\end{equation}
for every $C\geq C_0$. 
Combining this with the fact that $N_C(\lambda,\Lambda,\widehat{\Lambda}')=N_C(\lambda(j),\Lambda,\widehat{\Lambda}')$ for $C\leq j$ and sending $j \to +\infty$ we obtain
\begin{equation}
e^{\frac{aC}{\max(f_\lambda) + \mu}+d}<N_C(\lambda,\Lambda,\widehat{\Lambda}')
\end{equation}
for all $C\geq C_0$.
\qed

\subsection{Positivity of the topological entropy}
Using the results of the previous subsection, we prove Theorem \ref{theorem1}. For this, we first study Reeb chords from a fixed Legendrian knot $\Lambda$ to generic Legendrian fibers $\widehat{\Lambda}^z$ of a tubular neighbourhood $\mathcal{V}_\delta(\widehat{\Lambda})$ of a Legendrian knot $\widehat{\Lambda}$.

By the Weinstein tubular neighbourhood theorem for Legendrian submanifolds (\cite{Ge}), the Legendrian knot $\widehat{\Lambda}$ in $(Y,\xi)$ has a tubular neighbourhood $\mathcal{V}_\delta(\widehat{\Lambda}) \subset Y$ such that $(\mathcal{V}_\delta(\widehat{\Lambda}),\xi|_{\mathcal{V}_\delta(\widehat{\Lambda})})$ admits a contactomorphism to the local model $(S^1 \times \mathbb{D},\ker (\sin(\theta)dx + \cos(\theta) d y))$, where $\theta \in S^1$ and $z=(x,y)$ are coordinates of the disk $\mathbb{D}$, which takes  $\widehat{\Lambda}$ to the core circle $S^1 \times \{0\}$ of $S^1 \times \mathbb{D}$. We denote by $\widehat{\Lambda}^z$ the Legendrian submanifolds obtained by fixing the coordinate $z$ of $\mathbb{D}$. As the $\widehat{\Lambda}^z$ form a fibration of $\mathcal{V}_\delta(\widehat{\Lambda})$ by Legendrian knots we will refer to them as Legendrian fibers of $(\mathcal{V}_\delta(\widehat{\Lambda}),\xi|_{\mathcal{V}_\delta(\widehat{\Lambda})})$. We choose the neighbourhood $\mathcal{V}_\delta(\widehat{\Lambda})$ to be small enough so that all Legendrian fibers  $\widehat{\Lambda}^z$ are $\delta$-close to $\widehat{\Lambda}$ in the $C^3$ sense and $\mathcal{V}_\delta(\widehat{\Lambda})$ is disjoint from $\Lambda$.

In the disc $\mathbb{D}$ we will consider the Lebesgue measure, which is obtained from restricting the Lebesgue on $\mathbb{R}^2$ to $\mathbb{D}$. Therefore, when we say for almost every $z \in \mathbb{D}$ we mean almost every with respect to this measure.
\begin{lemma}\label{generic}
Let $\lambda$ be a contact form on a contact 3-manifold  $(Y,\xi)$, and $\Lambda$ and $\widehat{\Lambda}$ be disjoint Legendrian knots on $(Y,\xi)$. We consider a tubular neighbourhood $\mathcal{V}(\widehat{\Lambda})$ of $\widehat{\Lambda}$ such that $(\mathcal{V}(\widehat{\Lambda}),\xi|_{\mathcal{V}(\widehat{\Lambda})})$ admits a contactomorphism to the local model $(S^1 \times \mathbb{D},\ker (\sin(\theta)dx + \cos(\theta) d y))$ which takes  $\widehat{\Lambda}$ to $S^1 \times \{0\}$.
Then, for almost every $z \in \mathbb{D}$ the Legendrian knot $\widehat{\Lambda}^z$ is $(\lambda,\Lambda)$ transverse.
\end{lemma}
\textit{Proof:} Taking a parametrization $\nu: S^1 \to \Lambda$, we use the Reeb flow of $\lambda$ to define a map $F: S^1 \times \mathbb{R} \to Y$ by:
\begin{center}
$F(s,t)=\phi_{X_\lambda}^t(\nu(s))$.
\end{center}

The set $U = F^{-1}(\mathcal{V}(\widehat{\Lambda}))$ is an open subset  of $S^1 \times \mathbb{R}$.
Using the coordinates $(\theta,z)$ on $\mathcal{V}(\widehat{\Lambda})$ given by the contactmorphism from $(\mathcal{V}(\widehat{\Lambda}),\xi|_{\mathcal{V}(\widehat{\Lambda})})$ to $(S^1 \times \mathbb{D},\ker (dx + \theta d y))$ considered in the statement of the lemma, we define the map $\mathrm{pr}_{\mathbb{D}}: \mathcal{V}(\widehat{\Lambda})  \to \mathbb{D}$ to be the projection to the $z$ coordinate.
Composing $F\mid_{U}$ with $\mathrm{pr}_{\mathbb{D}}$ to obtain a smooth map:
\begin{center}
$\mathrm{pr}_{\mathbb{D}} \circ F\mid_U: U \to \mathbb{D}$.
\end{center}

A point $(s_0,t_0) \in S^1 \times \mathbb{R}$ is a  critical point of $\mathrm{pr}_{\mathbb{D}} \circ F\mid_U$ if, and only if, the closed curve $F_{t_0} : S^1 \to Y$ obtained by restricting $F$ to $S^1 \times \{t_0\}$ is tangent to a Legendrian fiber $\widehat{\Lambda}^z$ at $s_0$.
To see this, note that by the definition of $F$ we have $\partial_t (F\mid_{U}) = X_{\lambda} $ and $\partial_s (F\mid_{U}) \neq 0$. This implies that $\partial_t (\mathrm{pr}_{\mathbb{D}} \circ F\mid_U)$ is always non-zero, because the vector field $X_\alpha$ is never tangent to the Legendrian fibers $\widehat{\Lambda}^z$, and that $\partial_s (\mathrm{pr}_{\mathbb{D}} \circ F\mid_U))=0$ if, and only if, $\partial_s (F\mid_{U})$ is tangent to a Legendrian fiber. It follows that the regular values of $\mathrm{pr}_{\mathbb{D}} \circ F\mid_{U}$ in $\mathbb{D}$ are in bijection with the set of Legendrian fibers $\widehat{\Lambda}^z$ that are $(\lambda,\Lambda)$ transverse.

Applying Sard's theorem to $\mathrm{pr}_{\mathbb{D}} \circ F\mid_{U}$ we conclude that almost every element of $\mathbb{D}$ is a regular value of $\mathrm{pr}_{\mathbb{D}} \circ F\mid_U$, from which the lemma follows.
\qed

We are now ready to prove the main result of this Section.
\begin{thrm} \label{theorem1'}
Let $(Y,\xi )$ be a contact 3-manifold, and $\lambda_0$ be a hypertight contact form on $(Y,\xi)$. Assume that $\lambda_0$ is adapted to the pair of disjoint Legendrian knots $(\Lambda,\widehat{\Lambda})$.  Then, if the strip Legendrian contact homology $LC\mathbb{H}_{st}(\lambda_0,\Lambda \to\widehat{\Lambda})$ has exponential homotopical growth rate with exponential weight $a>0$, it follows that the Reeb flow of any $C^\infty$ contact form $\lambda$ on $(Y,\xi)$ has positive topological entropy.
Moreover, if we denote by $f_\lambda$ the positive function such that $\lambda=f_\lambda \lambda_0$, we have
\begin{equation} \label{eqthrm}
h_{top}(\phi_{X_{\lambda}}) \geq \frac{a}{\max(f_\lambda)}.
\end{equation}
\end{thrm}
\textit{Proof:}
Let $\lambda$ be a contact form on $(Y,\xi)$
and fix $\mu>0$.
We choose $\delta>0$ as in Proposition \ref{mainproposition} and let $\mathcal{V}_\delta(\widehat{\Lambda})$ be a tubular Weinstein tubular neighbourhood of $\widehat{\Lambda}$ whose Legendrian fibres are $\delta$-close to $\widehat{\Lambda}$ in the $C^3$ sense.
It follows from Lemma \ref{generic} that there exists a set $\mathfrak{W} \subset \mathbb{D}$ of full measure in $\mathbb{D}$ such that for every $z \in \mathfrak{W}$ the Legendrian $\widehat{\Lambda}^z$ is $(\lambda,\Lambda)$-transverse. With this choices it follows from Proposition \ref{mainproposition} that the number $N_C(\lambda,\Lambda,\widehat{\Lambda}^z)$ of Reeb chords of $X_{\lambda}$ from $\Lambda$ to the Legendrian fiber $\widehat{\Lambda}^z$ satisfies:
\begin{equation} \label{equationgrowth}
e^{\frac{aC}{\max(f_\lambda) + \mu}}<N_C(\lambda,\Lambda,\widehat{\Lambda}^z).
\end{equation}

Now choose a Riemannian metric on the manifold $Y$ which restricts in $\mathcal{V}_\delta(\widehat{\Lambda})$ to the Euclidean metric $d\theta \otimes d\theta + dx\otimes dx + dy\otimes dy$, for the coordinates $(\theta,x,y)\in S^1 \times \mathbb{D}$ on $\mathcal{V}_\delta(\widehat{\Lambda})$ used in Proposition \ref{mainproposition}.
This Riemannian metric induces a measure of area $Area(\Sigma)$ for all surfaces $\Sigma$ immersed in $Y$.
We want to estimate the area $Area^C(\Lambda)$ of the immersed cylinder $\{\phi_{X_\lambda}^t(\Lambda) \mid t\in [0,C]\}$. This is the image of the map $F_{C,\Lambda}: \Lambda \times [0,C] \to Y$ defined by $F_{C,\Lambda}(p,t)=\phi^t_{X_\lambda}(p)$. Defining $U_{\delta,C}:=F_{C,\Lambda}^{-1}(\mathcal{V}_\delta(\widehat{\Lambda}))$ we have:
\begin{equation}
Area^C(\Lambda) \geq  Area(F_{C,\Lambda}(U)) \geq {Area}(\mathrm{pr}_{\mathbb{D}}(F_{C,\Lambda}(U))),
\end{equation}
where the last area is taken with multiplicities and with respect to the Lebesgue measure in $\mathbb{D}$.
The inequality on the right side is true because the measure of area induced by our Riemannian metric on $\mathcal{V}_\delta(\widehat{\Lambda})$ coincides with the measure of area induced by the Euclidean metric $d\theta \otimes d\theta + dx\otimes dx + dy\otimes dy$.

From the estimate \eqref{equationgrowth} for the counting function $N_C(\lambda,\Lambda,\widehat{\Lambda}^z)$ for almost every $z=(x,y) \in \mathbb{D}$, we get:
\begin{equation}
Area(\mathrm{pr}_{\mathbb{D}}(F_{C,\Lambda}(U))) =  \int_{\mathfrak{W}} N_C(\lambda, \Lambda, \widehat{\Lambda}^z) dx \wedge dy  \geq \int_{\mathfrak{W}}  e^{\frac{aC}{\max(f_\lambda)+ \mu}} dx \wedge dy.
\end{equation}
It follows that
\begin{equation}
Area(\mathrm{pr}_{\mathbb{D}}(F_{C,\Lambda}(U))) \geq  \uppi \ e^{\frac{aC}{\max(f_\lambda)+ \mu}}.
\end{equation}

Combining this inequality with the Fubini type estimate proved in Appendix \ref{appendixA} we have:
\begin{equation}
\limsup_{t\to +\infty} \frac{\log (length(\phi_{X_\lambda}^t (\Lambda)))}{t}  \geq \frac{a}{\max(f_\lambda)+ \mu},
\end{equation}
where $length(\phi_{X_\lambda}^C (\Lambda)))$ is the length of the curve $\phi_{X_\lambda}^C (\Lambda)$ with respect to the Riemannian metric we chose.
We now invoke Yomdin's theorem, that implies that $h_{top}(\phi_{X_\lambda}) \geq \frac{a}{\max(f_\lambda)+ \mu}$. Since in Proposition \ref{mainproposition} the constant $\mu>0$ can be taken arbitrarily small, we obtain inequality \eqref{eqthrm}.
\qed

\section{Proof of Theorem 2} \label{section4}
In this section we prove the following theorem:
\begin{thrm}
Let $M$ be a closed oriented connected 3-manifold which can be cut along a nonempty family of incompressible tori into a family $\{M_i, 0 \leq i \leq q\}$ of irreducible manifolds with boundary such that the component $M_0$ satisfies:
\begin{itemize}
\item{$M_0$ is the mapping torus of a diffeomorphism $h: S \to S$ with pseudo-Anosov monodromy on a surface $S$ with non-empty boundary.}
\end{itemize}
Then $M$ can be given infinitely many different tight contact structures $\xi_k$, such that:
\begin{itemize}
\item There exist disjoint Legendrian knots $\Lambda$, $\widehat{\Lambda}$ on $(M,\xi_k)$, and a contact form $\lambda_k$ on $(M,\xi_k)$ adapted to the pair $\Lambda$ and $\widehat{\Lambda}$, such that $LC\mathbb{H}_{st}(\lambda_k,\Lambda \to \widehat{\Lambda})$ has exponential homotopical growth rate. It follows that every Reeb flow on $(M,\xi_k)$ has positive topological entropy.
\end{itemize}
\end{thrm}
 This result establishes the existence of many  contact 3-manifolds which possess a contact form and a pair of Legendrian knots that satisfy the hypothesis of Theorem \ref{theorem1}. The contact manifolds covered in Theorem \ref{theorem2} are among the ones constructed by Colin and Honda \cite{CH}.

We start with some preliminary notions.
 Let $S$ be a compact surface with non-empty boundary that has negative Euler characteristic. Let $n>0$ be the number of boundary components of $S$, and let $\mathfrak{y}_i$ for $1\leq i \leq n$ be the boundary components of $S$. For each $i\in \{1,...,n\}$ let $V_i$ be a tubular neighbourhood of $\mathfrak{y}_i$, which admits coordinates $(r,\theta) \in [-1,0]\times S^1$ identifying $\mathfrak{y}_i$ with $\{0\}\times S^1$.
 It follows from \cite[Lemma 4.1]{LW} that there exists a symplectic form $\omega$ on $S$ and a Liouville\footnote{Recall that $v$ is a Liouville vector field on $(S,\omega)$ if it satisfies $L_v \omega = \omega$ where $L_v \omega$ is the Lie derivative of $\omega$ by $v$.} vector field $v$ for $\omega$ such that:
 \begin{itemize}
 \item the vector field $v$ is Morse-Smale,
 \item there exist a Riemannian metric $g$ on $S$ and a function $\mathfrak{f}: S \to \mathbb{R}$ satisfying $\mathfrak{f} \leq 1$  and $\mathfrak{f}^{-1}(1) = \partial S$, such that  $v$ is the gradient vector field of  $\mathfrak{f}$ with respect to $g$,
 \item $v$ points outwards along all boundary components of $S$,
 \item the one form $\beta := i_v \omega$ satisfies $\beta= f_i(r)d\theta $ for the coordinates $(r,\theta)$ on the neighbourhood $V_i$  for a positive function $f_i$ that satisfies $f_i'>0$.
 \end{itemize}
Notice that it follows from these properties that all singularities of $v$ have index $0$ or $1$.

In a pseudo-Anosov mapping class choose a diffeomorphism  $h:S\to S$ that preserves $\omega$ and
which coincides with the identity map on all neighbourhoods $V_i$.
We now present a well-known recipe to produce a contact form on the the mapping torus $\Omega(S,h)$ of the map $h$. We first choose a smooth non-decreasing function $ F_0:\mathbb{R} \to [0,1]$   which satisfies $F_0(t)=0$ for $t\in (-\infty,-\frac{1}{100})$ and $F_0(t)=1$ for $t\in (\frac{1}{100},+\infty)$. For $i \in \mathbb{Z}$ we define $F_i(t) = F_0(t-i)$.  Fixing $\epsilon >0$, we define a 1-form $\widetilde{\alpha}$ on $\mathbb{R} \times S$ by
\begin{equation}
\widetilde{\alpha} = dt + \epsilon(1-F_i(t))(h^{i})^*\beta + \epsilon F_i(t)(h^{i+1})^*\beta \ \mbox{ for } \ t \in [i,i+1).
\end{equation}
This defines a smooth 1-form on $\mathbb{R} \times S$, and a simple computation shows that the 1-form $\widetilde{\alpha}$ is a contact form for $\epsilon$ small enough. For $t\in[0,1]$, the Reeb vector field $X_{\widetilde{\alpha}}$ is equal to $\partial_t + v(p,t)$, where $v(p,t)$ is the unique vector tangent to $S$ that satisfies $\omega(v(p,t), \cdot )= F'_0(t)(\beta -  h^* \beta)$.

Consider the diffeomorphism $H:\mathbb{R} \times S \to \mathbb{R} \times S$ defined by 
\begin{equation} \label{eqH}
H(t,p)=(t-1,h(p)). 
\end{equation}
The mapping torus $\Omega(S,h)$ is defined by
\begin{equation}
\Omega(S,h):=(\mathbb{R} \times S )/_{(t,p) \sim H(t,p)},
\end{equation}
and we denote by $\mathfrak{p}:\mathbb{R} \times S \to \Omega(S,h)$ the associated covering map.

Because $H^*\widetilde{\alpha}=\widetilde{\alpha}$, there exists a unique contact form $\alpha$ on $\Omega(S,h)$ such that $\mathfrak{p}^*\alpha =\widetilde{\alpha}$. Notice that in the neighbourhoods $S^1 \times V_i$  with coordinates $(t,r,\theta)$ of the connected components of $\partial \Omega(S,h)$, we have the formula $\alpha= d t + \epsilon f(r)d\theta$ for $\alpha$. This implies that $X_{\alpha}$ is tangent to $\partial \Omega(S,h)$.  The Reeb vector field $X_{\alpha}$ on $\Omega(S,h)$ is transverse to the surfaces $\mathfrak{p}^{-1} (\{t\} \times S)$ for $t\in \mathbb{R}/_{\mathbb{Z}}$. This implies that for any $t\in  \mathbb{R}/ {\mathbb{Z}}$ the surface  $\mathfrak{p}^{-1} (\{t\} \times S)$ is a global surface of section for the Reeb flow of $\alpha$. Moreover, by our expression for $X_{\widetilde{\alpha}}$ the first return map of the Reeb flow of $\alpha $ on any of these sections  is isotopic to $h$.

In the mapping torus $\Omega(S,h)$ consider the surface $S_0 = \mathfrak{p}^{-1} (\{0\} \times S)$. As mentioned above this surface is a copy of $S$ inside $\Omega(S,h)$, and we denote by $\mathfrak{i}_0 : S \to \Omega(S,h)$ the embedding associated to $S_0$.

\subsection{A special Legendrian knot in $\Omega(S,h)$} \label{sectionLegendrian}
We begin with the following lemma.
\begin{lemma}\label{lemmacurve}
There exists an embedded  curve $\eta: S^1 \to S$ contained in the interior of $S$ such that $\int_\eta \beta = 0$, and which is not homotopic to $\partial S$. Moreover, given $\delta >0$, we can choose $\eta$ such that $|\int_0^s \eta^* \beta | < \delta$ for\footnote{Here we are identifying $S^1$ with $\mathbb{R}/ \mathbb{Z}$.} all $s\in [0,1]$.
\end{lemma}
\textit{Proof:} As $v$ is a Liouville vector field on $(S,\omega)$ we know that $-v$ contracts the area form $\omega$ on $S$. This implies that $-v$ has no singularities of source type. Notice that at every point of $S$ the vector $-v$ belongs to the kernel of $\beta$.

Since $-v$ is a gradient Morse-Smale vector field, $-v$ has no closed orbits. We know from Morse theory that the flow of $-v$ can be used to produce a deformation retraction from $S$ to the set $\Delta$ defined as the union of the singularities of the flow of $-v$ and the unstable manifolds of its saddle singularities. Because $\Delta$ has the same homotopy type as $S$ it must contain a simple closed piecewise smooth curve $\gamma$ which is contained in the interior of $S$ and that cannot be homotoped to $\partial S$. Because $\gamma$ consists of trajectories of $-v$, we know that vectors tangent to $\gamma$ belong to the kernel of $\beta$. We orient $\gamma$ using a piecewise smooth parametrisation which we denote also by $\gamma:S^1 \to S$. Then $\int_0^s \gamma^*\beta= 0$ for all $s\in [0,1]$ where again $S^1 = \mathbb{R}/ \mathbb{Z}$ .

 Fix a point $p_0 \in \gamma$ which is not a vertex point. We can then smoothen $\gamma$ in small closed neighbourhoods of the vertices of $\gamma$ which are all disjoint from $p_0$ and produce a smooth embedded curve $\gamma'$ on $S'$. As $\gamma'$ coincides with $\gamma$ in a neighbourhood of $p_0$, we pick the orientation in $\gamma'$ which coincides with the orientation chosen for $\gamma$ in the region where the two curves coincide. We choose a parametrisation $\gamma': S^1 \to S$ of $\gamma'$ which induces the same orientation as the one we just defined. Given $\delta> 0$, we can ensure, by doing the smoothing in a sufficiently small neighbourhood of the vertices, that $q := \int_{\gamma'} \beta  \in (-\frac{\delta}{4}, \frac{\delta}{4})$ and that for all $s\in [0,1]$ we have $|\int_0^s \gamma'^*\beta | < \frac{\delta}{2}$.

Using that near the point $p_0$ the curve $\gamma'$ is tangent to $\ker\beta$, we find, for sufficiently small $\delta>0$, a $C^\infty$ small perturbation $\widetilde{\gamma}$ of the parametrised curve $\gamma'$ supported in a neighbourhood $V_{p_0}$ of $p_0$ such that, for the small region $U \subset V_{p_0}$ bounded by $\widetilde{\gamma}$ and $\gamma'$, and oriented such that $\partial U = \gamma' - \widetilde{\gamma}$, we have
\begin{equation}
\int_U \omega = q.
\end{equation}
Now Stokes' theorem and the fact that $q = \int_{\gamma'} \beta$ imply that $\int_{\widetilde{\gamma}} \beta=0$. This perturbation $\widetilde{\gamma}$, which aims at correcting the change in the integral after the smoothing of the corners, can be made explicitly if one uses Darboux coordinates in a neighbourhood of $p_0$, and can be taken so small enough that for all $s\in [0,1]$ we have $|\int_0^s \widetilde{\gamma}^*\beta | < \delta$.

Taking $\eta: S^1 \to S$ to  be the parametrised curve $\widetilde{\gamma}$  we have shown that $\eta$ has all the properties stated in the lemma. \qed

\

Let $\delta < \frac{1}{100}$ and take the parametrised curve $\eta$ obtained in Lemma \ref{lemmacurve}.
Because $\int_{\eta} \beta = 0$ we have that the curve $\eta$ is the Lagrangian projection of a Legendrian curve $\Lambda_0$ in $(\Omega(S,h),\alpha)$. To make this precise we first define the function $t(s):= -\int_0^s \beta(\eta'(s)) $. We then define the curve $\Lambda_0: S^1 \to \mathbb{R} \times S $ by $ \Lambda_0(s):=(t(s),\eta(s))$. As  $\delta < \frac{1}{100}$, and $|t(s)| <\delta$ we conclude that image of $\Lambda_0$ is contained in the set $[-\delta,\delta] \times S $ in which the contact form $\alpha$ is given by $dt + \epsilon \beta$. Using this identity we compute that $\alpha(\Lambda_0'(s))=0$ at every $s\in [0,1]$ which implies that $\Lambda_0$ is a Legendrian curve of $\alpha$. From our construction it is clear that $\Lambda_0$ is a graph over the curve $\eta$ which is embedded in $\{0\}\times S$. We will denote by $\Lambda$ the Legendrian knot $\mathfrak{p}(\Lambda_0)$ in $(\Omega(S,h),\alpha)$.

\

We are now ready to prove the following lemma.
\begin{lemma} \label{lemmachords}
If $\check{\Lambda}$ is a Legendrian curve in $(\Omega(S,h),\ker \alpha)$ which is sufficiently close to $\Lambda$ in the $C^\infty $ topology, then there exists no Reeb chord $c$ of $\alpha$ from $\check{\Lambda}$ to itself such that $[c]$ is the trivial element in $\pi_1(\Omega(S,h),\check{\Lambda})$.
\end{lemma}
\textit{Proof:} If a Legendrian $\check{\Lambda}$ is sufficiently close to $\Lambda$ in the $C^\infty $ topology then there exists a knot $\check{\Lambda}_0$ in $\mathbb{R}\times S$ which is a connected component of $\mathfrak{\pi}^{-1}(\check{\Lambda})$, which is also contained in $[-\frac{1}{100}, \frac{1}{100}] \times S$, and which is also a graph over an embedded curve $\check{\eta}$ in $S$.

Assume there exists  a Reeb chord $c$  of $\alpha$ from $\check{\Lambda}$ to itself such that $[c]$ is the trivial element in $\pi_1(\Omega(S,h),\check{\Lambda})$. Then, the fact the $[c]$ is the trivial element in $\pi_1(\Omega(S,h),\check{\Lambda})$ would imply that $c$ lifts to a Reeb chord of $\widetilde{\alpha}$ from $\check{\Lambda}_0$ to itself in $\mathbb{R}\times S$. Since the Reeb vector field of $\widetilde{\alpha}$ is $\partial_t + v(p,t)$ and  $\check{\Lambda}_0$ is a graph over an embedded curve $\check{\eta}$ in $S$,  it is clear that there is no Reeb chord of $ \widetilde{\alpha}$ from $\check{\Lambda}_0$ to itself. This contradiction proves the lemma.
\qed

\

We let $\widetilde{\Lambda}_0$ be a Legendrian knot in $(\mathbb{R}\times S,\ker \widetilde{\alpha})$ obtained as a $C^\infty$ small perturbation of $\Lambda_0$ which is disjoint from $\Lambda_0$ and is contained in $[-\frac{1}{100}, \frac{1}{100}] \times S$, and denote by $\widehat{\Lambda}$ the Legendrian curve $\mathfrak{p}(\widetilde{\Lambda}_0)$ in $(\Omega(S,h), \ker \alpha)$. Choosing the perturbation generically, we can choose  $\widetilde{\Lambda}_0$ so that:
\begin{itemize}
\item $\widetilde{\Lambda}_0$ is also the graph over an embedded curve $\widehat{\eta}$ in $S$,
\item all the Reeb chords of $\alpha$ from $\Lambda$ to $\widehat{\Lambda}$ are transverse.
\end{itemize}

There is another genericity condition  that we want to impose on the Reeb chords from $\Lambda$ to $\widehat{\Lambda}$. We want to guarantee that no Reeb orbits of our Reeb flow intersect Reeb chords from  $\Lambda$ to $\widehat{\Lambda}$. To guarantee this we have to perturb the contact form $\alpha$ in the interior of $\Omega(S,h)$.

\begin{lemma} \label{adapted}
There exists a contact form $\widehat{\alpha}$ on $(\Omega(S,h),\ker \alpha)$ which is obtained by perturbing $\alpha$ in the interior of $\Omega(S,h)$ and that satisfies:
\begin{itemize}
\item $\widehat{\alpha}$ has no contractible Reeb orbits,
\item all Reeb chords of $\widehat{\alpha}$ from $\Lambda$ to itself represent non-trivial elements in $\pi_1(\Omega(S,h),{\Lambda})$,
\item all Reeb chords of $\widehat{\alpha}$ from $\widehat{\Lambda}$ to itself represent non-trivial elements in $\pi_1(\Omega(S,h),{\widehat{\Lambda}})$,
\item the Reeb orbits of $\widehat{\alpha}$ do not intersect the Reeb chords of $\widehat{\alpha}$ from ${\Lambda}$ to $\widehat{\Lambda}$,
\item there exist no Reeb chords of $\mathfrak{p}^*\widehat{\alpha}$ from $\Lambda_0$ to itself, or from $\widetilde{\Lambda}_0$ to itself.
\end{itemize}
\end{lemma}
\textit{Proof:} If $\widehat{\alpha}$ is a sufficiently small perturbation of ${\alpha}$ in the interior of $\Omega(S,h)$ then we can guarantee that the $\partial_t $ component of the vector field $\pi^*X_{\widehat{\alpha}}$ is bigger than $\frac{1}{2}$.

In order to prove the second and third statements we first pick $N>0$ such that the curves ${\Lambda}_0$ and $\check{\Lambda}_0$ are contained in  $[-N,N] \times S$. It follows that any Reeb chord of  $\pi^*{\widehat{\alpha}}$ from ${\Lambda}_0$ to itself has action $\leq 4N$. Indeed, for any point $p \in  {\Lambda}_0$ the trajectory $\phi^t_{X_{\pi^*{\widehat{\alpha}}}}$ stays outside $[-N,N] \times S$ for $t \geq 4N$. The same is true for $\check{\Lambda}_0$. This implies that any Reeb chord of $\widehat{\alpha}$ from $\Lambda$ to itself that represents the trivial element in $\pi_1(\Omega(S,h),{\Lambda})$ must have action $\leq 4N$, the same being true for $\widehat{\Lambda}$.

We claim that if $\widehat{\alpha}$ is sufficiently close to $\alpha$, then there exists no Reeb chord of $\widehat{\alpha}$ from $\Lambda$ to itself that represent the trivial element in $\pi_1(\Omega(S,h),{\Lambda})$, the same being true for $\widehat{\Lambda}$. Indeed, otherwise we could take a sequence of contact forms $ \widehat{\alpha}_n$ converging to $ \widehat{\alpha}$ all of them having Reeb chords $c_n$ from $\Lambda$ to itself that represent the trivial element in $\pi_1(\Omega(S,h),{\Lambda})$ with action $\leq 4N$. Since the action of the Reeb chords $c_n$ is bounded by $4N$, the Arzel\'a-Ascoli theorem implies that $c_n$ has a convergent subsequence that converges to a Reeb chord $c$ of $\alpha$ from $\Lambda$ to itself that must be trivial in  $\pi_1(\Omega(S,h),{\Lambda})$. As this contradicts Lemma \ref{lemmachords}, we have proved our claim for $\Lambda$.  The claim for $\widehat{\Lambda}$ is proved in the same way.

The last claim of the lemma is proved by an approximation argument identical to the one used to prove the second and third claims.

We pick $\widehat{\alpha}$ to be sufficiently close to $\alpha$, so that there exists no Reeb chord of $\widehat{\alpha}$ from $\Lambda$ to itself that represent the trivial element in $\pi_1(\Omega(S,h),{\Lambda})$, and no Reeb chord of $\widehat{\alpha}$ from $\widehat{\Lambda}$ to itself that represent the trivial element in $\pi_1(\Omega(S,h),\widehat{\Lambda})$. Moreover by choosing our perturbation generically we can choose $\widehat{\alpha}$ such that:
\begin{itemize}
\item the Reeb orbits of $\widehat{\alpha}$ do not intersect the Reeb chords of $\widehat{\alpha}$ from ${\Lambda}$ to $\widehat{\Lambda}$,
\item all the Reeb chords of $\widehat{\alpha}$ from $\Lambda$ to $\widehat{\Lambda}$ are transverse.
\end{itemize}

To complete the proof of the lemma we must now show that $\widehat{\alpha}$ has no contractible Reeb orbits. But since the $\partial_t$ coordinate of the vector field $\pi^*X_{\widehat{\alpha}}$ is always positive, we know that the flow of $\pi^*X_{\widehat{\alpha}}$ has no periodic orbits. It follows that $\widehat{\alpha}$ has no contractible Reeb orbits.
\qed

\subsection{Another covering of $\Omega(S,h)$}

We now construct another covering $\widetilde{\mathfrak{p}}: \mathbb{R} \times S \to \Omega(S,h)$ of the mapping torus $\Omega(S,h)$. We first define
\begin{equation}
\widetilde{\mathfrak{p}}(t,p):= \phi^t_{X_{\widehat{\alpha}}}(p).
\end{equation}
Because $S_0$ is a global surface of section for the Reeb flow of ${\widehat{\alpha}}$ we conclude that $\widetilde{\mathfrak{p}}$ is indeed a covering map.

Using the homotopy lifting theorem we know that there exists a unique diffeomorphism $\mathbf{\Psi}: \mathbb{R} \times S \to \mathbb{R} \times S$ that is the identity when restricted to $\{0\} \times S$ and fits in the commutative diagram
\begin{center}
\begin{tikzpicture}
  \matrix (m) [matrix of math nodes,row sep=4em,column sep=3em,minimum width=3em]
  {
     \mathbb{R} \times S & \mathbb{R} \times S \\
     \Omega(S,h)\\};
  \path[-stealth]
    (m-1-1) edge node [left] {$\mathfrak{p}$} (m-2-1)
            edge node [above] {$\mathbf{\Psi}$} (m-1-2)
         (m-1-2) edge node [below right] {$\widetilde{\mathfrak{p}}$} (m-2-1);
\end{tikzpicture}
\end{center}
It is clear that the pullback of $X_{{\widehat{\alpha}}}$ by $\widetilde{\mathfrak{p}}$ is the vector field $\partial_t$, and as observed before this pullback is of the form $\mathfrak{b}(p,t)\partial_t + \widehat{v}(p,t)$ where $\widehat{v}(p,t)$ is tangent to $S$ and $\mathfrak{b}(p,t)$ is a positive function. Because $\mathfrak{b}(p,t)$ is invariant by the map $H$ we know that there exists a constant $\mathbf{k} > 1$ such that $\frac{1}{\mathbf{k} }<\mathfrak{b}(p,t)< \mathbf{k} $.

From the fact   $\mathbf{\Psi}: \mathbb{R} \times S \to \mathbb{R} \times S$ is the identity when restricted to $\{0\} \times S$ we know that $\mathbf{\Psi}$ is isotopic to the identity. From the expression defining $\widetilde{\mathfrak{p}}$  and the fact that  $\mathbf{\Psi}$ fits in the above diagram  we know that $\mathbf{\Psi}^*\partial_t = \mathfrak{b}(p,t)\partial_t + \widehat{v}(p,t)$, and therefore that  $\mathbf{\Psi}$ conjugates the flows of $\partial_t$ and $\mathfrak{b}(p,t)\partial_t + \widehat{v}(p,t)$. This property of $\mathbf{\Psi}$ and the fact $\frac{1}{\mathbf{k} }<\mathfrak{b}(p,t)< \mathbf{k} $ imply that there exist constants $\mathfrak{A} > 1$ and $\mathfrak{B}>0$ such that
\begin{equation} \label{estimatepsi}
\frac{1}{\mathfrak{A}}t + \mathfrak{B}  < \mathrm{pr}_1(\mathbf{\Psi}(t,p)) < \mathfrak{A}t + \mathfrak{B},
\end{equation}
where $\mathrm{pr}_1: \mathbb{R} \times S \to \mathbb{R}$ is the projection onto the first coordinate.

\subsection{Contact 3-manifolds containing $(\Omega(S,h),\widehat{\alpha})$ as a component}
We now proceed to construct contact manifolds containing $(\Omega(S,h),\widehat{\alpha})$.
Let $M$ be a closed connected oriented  3-manifold which can be cut along a non-empty family of incompressible tori into a family $\{M_i, 0 \leq i \leq q\}$ of irreducible manifolds with boundary, such that the component $M_0$ is diffeomorphic to $\Omega(S,h)$.
Then one can construct hypertight contact forms on $M$ which match with $\alpha$ in the component $M_0$. More precisely, we have the following result due to Colin and Honda, and Vaugon:
\begin{prop} \label{colinhonda} (\cite{CH,Vaugon})
Let $M$ be a closed connected oriented  3-manifold which can be cut along a non-empty family of incompressible tori into a family $\{M_i, 0 \leq i \leq q\}$ of irreducible manifolds with boundary, such that the component $M_0$ is diffeomorphic to $\Omega(S,h)$. Then, there exists an infinite family $\{\xi_k, k \in \mathbb{Z}\}$ of non-diffeomorphic contact structures on $M$  such that
\begin{itemize}
\item for each $k \in \mathbb{Z}$ there exists a hypertight contact form $\lambda_k$ on $(M,\xi_k)$ which coincides with $\widehat{\alpha}$ on the component $M_0$.
\end{itemize}
\end{prop}

We briefly recall the construction of the contact forms $\lambda_k$, and refer the reader to \cite{CH,Vaugon} for the details. For $i \geq 1$, we apply \cite[Theorem 1.3]{CH} to obtain a hypertight contact form $\alpha_i$ on $M_i$ which is compatible with the orientation of $M_i$, and whose Reeb vector field $X_{\alpha_i}$ is tangent to the boundary of $M_i$. On the special piece $M_0$ we consider the contact form $\alpha_0$ equal to $\widehat{\alpha}$ constructed above.

Let $\{\mathfrak{T}_j \mid 1 \leq \ j \leq m \}$ be the family of incompressible tori along which we cut $M$ to obtain the pieces $M_i$. Then the contact forms $\alpha_i$ give a hypertight contact form on each component of $M \setminus \bigcup_{j \geq1}^m \mathbb{V}(\mathfrak{T}_j) $, where $\mathbb{V}(\mathfrak{T}_j)$ is a small open neighborhood of $\mathfrak{T}_j$. This gives a contact form $\widehat{\lambda}$ on $M \setminus \bigcup_{j \geq1}^m \mathbb{V}(\mathfrak{T}_j) $.
 Using an interpolation process (see \cite[Section 7]{Vaugon}), one can construct contact forms on the neighborhoods $\overline{\mathbb{V}(\mathfrak{T}_j)}$ which coincide with $\widehat{\lambda}$ on $\partial \overline{\mathbb{V}(\mathfrak{T}_j)}$. The interpolation process is not unique and can be done in many different ways as to produce  an infinite family of distinct contact forms $\{ \lambda_k \ | \ k\in \mathbb{Z}\}$ on $M$ that extend $\widehat{\lambda}$,  and which are associated to contact structures $\xi_k:= \ker \lambda_k$ that are all non-diffeomorphic. The contact topological invariant used to distinguish the contact structures $\xi_k$ is the \textit{Giroux torsion} (see \cite[Section 7]{Vaugon}).

\

 The Legendrian curves $\Lambda$ and $\widehat{\Lambda}$ in $(\Omega(S,h),\ker \alpha)$  are obviously also Legendrian curves of $(M,\xi_k)$. Because the Reeb flow of $\lambda_k$ leaves the region $\Omega(S,h)$ invariant, we know that all Reeb chords of $\lambda_k$ from $\Lambda$ to $\widehat{\Lambda}$ are  contained in $\Omega(S,h)$ and are therefore Reeb chords of $\alpha$. We conclude that $\widehat{\Lambda}$ is $(\lambda_k, \Lambda)$-transverse.

\begin{prop} \label{propadapted}
 All Reeb chords of $\lambda_k$ from $\Lambda$ to itself represent non-trivial elements in $\pi_1(M,\Lambda)$, and likewise all Reeb chords of $\lambda_k$ from $\widehat{\Lambda}$ to itself represent non-trivial elements in $\pi_1(M,\widehat{\Lambda})$.
\end{prop}

\textit{Proof:} We will show that Lemma \ref{adapted} implies the proposition. We present the proof of the statement for the Legendrian curve $\Lambda$. The proof for $\widehat{\Lambda}$ is identical.

By contradiction suppose there is a smooth map $\mathfrak{g}: \overline{\mathbb{D}^2} \to M$ such that $\mathfrak{g}|_{S^1} $ is the concatenation of a Reeb chord   $c$ of $\lambda_k$ from $\Lambda$ to itself with a path $\gamma \subset \Lambda$. Since $\mathfrak{g}|_{S^1} $ is contractible in $M$ and is contained in one component of the JSJ-decomposition of $M$, we conclude that $\mathfrak{g}|_{S^1} $ must be contractible in $M_0$; for this fact see for example the proof of \cite[Corollaire 1.6]{CH}. However, this contradicts Lemma \ref{adapted}.
 \qed

 \

 We have the following important corollary of Proposition \ref{propadapted}.
\begin{coro}
The contact form $\lambda_k$ on $(M,\xi_k)$ is adapted to the disjoint Legendrian knots ${\Lambda}$ and $\widehat{\Lambda}$.
\end{coro}
\textit{Proof:} We have to show that the triple $(\lambda_k,\Lambda,\widehat{\Lambda})$ satisfies conditions (a), (b), (c) and (d) from Section~\ref{sectionstrip}. Condition (a) follows from Proposition~\ref{colinhonda} which says that the contact form $\lambda_k$ is hypertight. Conditions (b), (c) and (d) follow from Proposition~\ref{propadapted}. \qed

\subsection{Exponential homotopical growth of $LC\mathbb{H}_{st}(\tau,\Lambda \to \widehat{\Lambda})$}
To study the growth rate of $LC\mathbb{H}_{st}(\tau,\Lambda \to \widehat{\Lambda})$ we will consider some special relative homotopy classes of paths from $\Lambda$ to $\widehat{\Lambda}$.

\begin{defi} \label{defrelnielsen}
Let $c_1$ and $c_2$ be Reeb chords from $\Lambda$ to $\widehat{\Lambda}$. We say that $c_1$ and $c_2$ are in the same \textbf{relative Nielsen class} if, and only if, there exists a smooth strip $\mathfrak{e}:[0,1]\times[0,1] \to \Omega(S,h) $ such that:
\begin{itemize}
\item $\mathfrak{e}({0}\times [0,1])$ is a path in $\Lambda$ and $\mathfrak{e}({1}\times [0,1])$ is a path in $\widehat{\Lambda}$,
\item  $\mathfrak{e}( [0,1] \times {0})= c_1$ and $\mathfrak{e}( [0,1] \times {1})= c_2$.
\end{itemize}
\end{defi}

Relative Nielsen classes are just homotopy classes of paths from $\Lambda$ to $\widehat{\Lambda}$ in the mapping torus $\Omega(S,h) $. Our first step is to prove that the relative Nielsen classes generate a partition of $LCH_{st}(\tau,\Lambda \to \widehat{\Lambda})$ in subcomplexes because they can be regarded as elements in the set $\Sigma_{\Lambda \to \widehat{\Lambda}}$ of homotopy classes of paths from $\Lambda$ to $\widehat{\Lambda}$  in $M$.

\begin{lemma} \label{relnielsen}
Let $c_1$ and $c_2$ be Reeb chords from $\Lambda$ to $\widehat{\Lambda}$, and $\mathfrak{e}:[0,1]\times[0,1] \to M$ such that:
 \begin{itemize}
\item{ $\mathfrak{e}(\{0\}\times [0,1])$ is a path in $\Lambda$ and $\mathfrak{e}(\{1\}\times [0,1])$ is a path in $\widehat{\Lambda}$,}
\item{ $\mathfrak{e}( [0,1] \times \{0\})=$ is a parametrisation of $c_1$ and $\mathfrak{e}( [0,1] \times \{1\})$ is a parametrisation of $ c_2$.}
\end{itemize}
Then, there exists a strip $\mathfrak{e}':[0,1]\times[0,1] \to \Omega(S,h)$ such that $\mathfrak{e}'(\partial([0,1]\times[0,1])) = \mathfrak{e}(\partial([0,1]\times[0,1]))$.
\end{lemma}
\textit{Proof:} The proof is very similar to the one of Proposition \ref{propadapted} above, so we will only give an outline.

By genericity we can assume that the map of $\mathfrak{e}$ is transverse to $\partial(\Omega(S,h))$. The set $\mathfrak{e}^{-1}(\partial(\Omega(S,h)))$ consists of a finite collection of circles $w_1$,...,$w_k$. For each $j \in \{1,...,k\}$ the circle $\mathfrak{e}(w_j)$ is contractible.

The assumption that $\partial(\Omega(S,h))$ is incompressible implies that the circles $\mathfrak{e}(w_1)$,...,$\mathfrak{e}(w_k)$ are also contractible in $\partial(\Omega(S,h))$.  Each circle $w_j$ bounds a disk $\mathfrak{d}_j$ in $\mathbb{D}$. We can choose a subset $\mathfrak{E} \subset \{1,...,k\}$ such that:
\begin{itemize}
\item for every $i \in   \{1,...,k\}$ there exists $j \in \mathfrak{E}$ such that$\mathfrak{d}_i \subset \mathfrak{d}_j$,
\item if $j \neq i $ are both elements in $\mathfrak{E}$ then $\mathfrak{d}_j$ and $\mathfrak{d}_i$ do not intersect.
\end{itemize}

By a cut and paste procedure identical to the one used in Proposition~\ref{propadapted} we can cut off the discs $\mathfrak{e}(\mathfrak{d}_j)$ for each $j \in \mathfrak{E}$, and replace them by other discs to obtain the map $\mathfrak{e}' :[0,1]\times[0,1] \to \Omega(S,h)$ which coincides with $\mathfrak{e}$ outside the disks $\mathfrak{d}_j$. The map $\mathfrak{e}'$ is the desired one. \qed

Recall from Section \ref{sectionstrip} that the differential $\partial_{st}$ of the strip Legendrian contact homology $LC\mathbb{H}_{st}(\tau,\Lambda \to \widehat{\Lambda})$ counts pseudoholomorphic strips $\widetilde{u}: \mathbb{R} \times [0,1] \to \mathbb{R} \times M$ in the symplectization of $\lambda_k$ with Fredholm index $1$ satisfying the boundary conditions:
\begin{itemize}
\item $\widetilde{u}(\mathbb{R} \times \{0\}) \subset \mathbb{R} \times \Lambda$,
\item $\widetilde{u}(\mathbb{R} \times \{1\}) \subset \mathbb{R} \times \widehat{\Lambda}$.
\end{itemize}

It follows from Lemma \ref{relnielsen} that the relative Nielsen classes can be seen as elements in $\Sigma_{\Lambda \to \widehat{\Lambda}}$. More precisely, denoting by $\mathfrak{R}$ the set of relative Nielsen classes,  we have a map 
\begin{equation} \label{eqI}
\mathcal{I}:\mathfrak{R} \to  \Sigma_{\Lambda \to \widehat{\Lambda}},
\end{equation}
 defined as follows: given a relative Nielsen class $\rho$, we pick a Reeb chord $c\in \rho$ and define $\mathcal{I}(\rho)$ to be the class of $[c] \in \Sigma_{\Lambda \to \widehat{\Lambda}}$. It is easy to see that $\mathcal{I}$ is well defined and the Lemma~\ref{relnielsen} implies that $\mathcal{I}$ is injective.

 Since $\lambda_k$ coincides with $\widehat{\alpha}$ on $M_0$ we know that $M_0$ is invariant by the Reeb flow of $\lambda_k$. It follows that all the Reeb chords in $\mathcal{T}_{\Lambda \to \widehat{\Lambda}}(\lambda_k)$ are contained in the component $M_0=\Omega(S,h)$, and therefore belong to elements in $\Sigma_{\Lambda \to \widehat{\Lambda}}$ which are in the image of our map $\mathcal{I}$.
We can therefore write $LC\mathbb{H}_{st}(\tau,\Lambda \to \widehat{\Lambda})$ as a direct sum:

\begin{equation}
LC\mathbb{H}_{st}(\tau,\Lambda \to \widehat{\Lambda}) = \bigoplus_{\varrho \in \mathfrak{R}} LC\mathbb{H}^{\mathcal{I}(\varrho)}_{st}(\tau,\Lambda \to \widehat{\Lambda})
\end{equation}

\subsubsection{The relative Nielsen classes}

We will use the coverings $\widetilde{\mathfrak{p}}: \mathbb{R} \times S \to \Omega(S,h)$ and ${\mathfrak{p}}: \mathbb{R} \times S \to \Omega(S,h)$ to obtain information about the relative Nielsen classes. The curve $\Lambda_0$ constructed in Section \ref{sectionLegendrian} is a lift of $\Lambda$ for the covering ${\mathfrak{p}}: \mathbb{R} \times S \to \Omega(S,h)$. We then fix in $\mathbb{R} \times S$ the lift $\mathbf{\Psi}(\Lambda_0) \subset \widetilde{\mathfrak{p}}^{-1}(\widehat{\Lambda})$ of $\Lambda$ for the covering $\widetilde{\mathfrak{p}}: \mathbb{R} \times S \to \Omega(S,h)$.

We want to describe the set $\widetilde{\mathfrak{p}}^{-1}(\widehat{\Lambda})$ of lifts of $\widehat{\Lambda}$ to $\mathbb{R} \times S$ via $\widetilde{\mathfrak{p}}$. We begin by fixing the lift $\widetilde{\Lambda}_0$ of $\widehat{\Lambda}$ for $\mathfrak{p}$ constructed in Section \ref{sectionLegendrian}. We know that
\begin{equation}
 {\mathfrak{p}}^{-1}(\widehat{\Lambda}) = \bigcup_{n \in \mathbb{Z}} H^n(\widetilde{\Lambda}_0),
\end{equation}
where $H$ is the map defined in \eqref{eqH}.
It follows that
\begin{equation}
 \widetilde{\mathfrak{p}}^{-1}(\widehat{\Lambda}) = \bigcup_{n\in \mathbb{Z}} \mathbf{\Psi}(H^n(\widetilde{\Lambda}_0)).
\end{equation}
We define $\widehat{\Lambda}_n:= \mathbf{\Psi}(H^{-n}(\widetilde{\Lambda}_0))$.

Given a Reeb chord $c$ from $\Lambda$ to $\widehat{\Lambda}$ we define by $\widetilde{c}$ the unique lift of $c$ via the covering map $\widetilde{\mathfrak{p}}$ which has its starting point in $\mathbf{\Psi}(\Lambda_0)$.
It is not difficult to see that if $c_1$ and $c_2$ are Reeb chords from $\Lambda$ to $\widehat{\Lambda}$ that are in the same relative Nielsen class, then $\widetilde{c_1}$ and $\widetilde{c_2}$ have to have endpoints in the same lift $\widehat{\Lambda}_n$ of $\widehat{\Lambda}$. We will see, however, that this condition is far from sufficient to guarantee that $c_1$ and $c_2$ are in the same relative Nielsen class.

We claim that every $\widehat{\Lambda}_n$ is a graph over an embedded curve $\widehat{\eta}_n$ of $S$. If this was not the case, it would follow that there exists a trajectory of the vector field $\partial_t$ which starts and ends in $\widehat{\Lambda}_n$. Applying the map $H^{-n} \circ \Psi^{-1}$ to this trajectory we would obtain a trajectory of the Reeb flow of $\widetilde{\alpha}$ that starts and ends in $\Lambda_0$, something that Lemma \ref{adapted} tells us does not exist. From Lemma \ref{adapted} we also know that $\mathbf{\Psi}(\Lambda_0)$ is a graph over an embedded curve $\eta_0$ in $S$.

\begin{remark} \label{above}
We claim that there exists $n_0$ such that if $n\geq n_0$ then the curve $\widehat{\Lambda}_n$  is above $\mathbf{\Psi}({\Lambda}_0)$ in the sense that for all $(t_0,p_0) \in  \mathbf{\Psi}({\Lambda}_0)$ and $(t_n,p_n) \in \widehat{\Lambda}_n$ we have
\begin{equation}
t_n > t_0.
\end{equation}
To see this let $N>0$ be such that  $\Lambda_0 \subset [-N,N] \times S$ and $\mathbf{\Psi}(\Lambda_0) \subset [-N,N] \times S$. We then know that $H^{-n}(\Lambda_0) \subset [-N+n,N+n] \times S$ . Combining this with \eqref{estimatepsi} we see that there exists $n_0$ such that if $n\geq n_0$ then $\widehat{\Lambda}_n = \mathbf{\Psi}( H^{-n}(\widetilde{\Lambda}_0)) \subset [2N,+\infty) \times S$.
\end{remark}

Let $\mathrm{pr}_2: \mathbb{R} \times S \to S$ be the projection on the second coordinate. The curves $\eta_0$ and $\widehat{\eta}_n= \mathrm{pr}_2(\widehat{\Lambda}_n)$ are embedded curves in $S$. Because $\partial_t= \mathfrak{p}^*X_{\widehat{\alpha}}$ we know that if $n \geq n_0$ the Reeb chords from $\Lambda_0$ to $\widehat{\Lambda}_n$ are in one-to-one correspondence with the intersection points of $\eta_0$ and $\widehat{\eta}_n$. By the same reason  we know that  the transversality of all the Reeb chords from $\Lambda$ to $\widehat{\Lambda}$ implies the transversality of $\eta_0$ and $\widehat{\eta}_n$ for every $n \geq n_0$. We now proceed with the following characterization of the relative Nielsen classes.
\begin{prop}\label{nielsen}
Let $c_1$ and $c_2$ be Reeb chords in $\mathcal{T}_{\Lambda \to \widehat{\Lambda}}(\widehat{\alpha})$ with $p_1:=\mathrm{pr}_2 (\widetilde{c}_1)$ and $p_2:=\mathrm{pr}_2 (\widetilde{c}_2)$. Then $c_1$ and $c_2$ are in the same relative Nielsen class if, and only if, $\widetilde{c}_1$ and $\widetilde{c}_2$ have end points in the same $\widehat{\Lambda}_n$, and there exists a map $v:[0,1] \times [0,1] \to S$, such that:
\begin{eqnarray}
v([0,1] \times \{0\})= p_1, \\
v([0,1] \times \{1\})= p_2, \\
v(\{0\} \times [0,1]) \subset \eta_0, \\
v(\{1\} \times [0,1]) \subset \widehat{\eta}_n.
\end{eqnarray}
\end{prop}

\textit{Proof:} Suppose $c_1$ and $c_2$ are in the same relative Nielsen class. We take the map $\mathfrak{e} :[0,1] \times [0,1] \to \Omega(S,h)$ given in Definition \ref{defrelnielsen}, and consider its lift $\widehat{\mathfrak{e}}:[0,1] \times [0,1] \to \mathbb{R} \times S$ to $\mathbb{R}\times S$ via $\widetilde{\mathfrak{p}}$, that satisfies $\widehat{\mathfrak{e}}([0,1] \times \{0\})= \widetilde{c}_1$ and $\widehat{\mathfrak{e}}([0,1] \times \{1\})= \widetilde{c}_2$. It is easy to see that taking $v=\mathrm{pr}_2 \circ \widehat{\mathfrak{e}}$ gives a strip in $S$ satisfying the conditions in the statement of the proposition. This proves one implication.

To prove the reverse implication take a map $v[0,1] \times [0,1] \to S$ satisfying the conditions in the statement of the proposition. By taking the path $v(\{0\} \times [0,1])\subset \eta_0$ there exists a unique function $g_0 :[0,1] \to \mathbb{R}$ such that the path $\gamma_0 (s)$ defined by $ \gamma_0 (s) = (v(0,s),g_0(s))$ is a path in $\mathbf{\Psi}(\Lambda_0)$. Analogously there exists a function $g_1:[0,1] \to \mathbb{R}$ such that $ \gamma_1 (s) = (v(1,s),g_1(s))$ is a path in $\widehat{\Lambda}_n$. Take $f:[0,1] \times [0,1] \to \mathbb{R}$ to be an homotopy between $g_0$ and $g_1$ that is, $f(0,s) = g_0(s)$ and $f(1,s)= g_1(s)$. Now we can define the strip $\widehat{\mathfrak{e}}(r,s):=(v(r,s),f(r,s))$ in $\mathbb{R} \times S$, and considering $\mathfrak{e}:=\widetilde{\mathfrak{p}} \circ \widehat{\mathfrak{e}}$ we obtain a strip in $\Omega(S,h)$ which satisfies the conditions of the definition of relative Nielsen classes for $c_1$ and $c_2$. This finishes the second implication and the proof of the proposition.
\qed

Proposition \ref{nielsen} gives a complete description of the relative Nielsen classes. It also shows how to identify different relative Nielsen classes of Reeb chords by looking at properties of intersection points of the curves $\widehat{\eta}_n$ and $\eta_0$. This is the crucial link that will allow us to use the pseudo-Anosov monodromy of $h^*$ to estimate the growth of the number of relative Nielsen classes.
Among the relative Nielsen classes, the subset of relative Nielsen classes with an odd number of chords will be of special importance to us. We will call them fundamental relative Nielsen classes and denote their set by $\mathfrak{R}^f$.

By the discussion above we can partition the set $\mathfrak{R}^f$ in subsets $\mathfrak{R}^f_n$ defined by the following rule: $\varrho \in \mathfrak{R}^f_n$ if, and only if, for every Reeb chord $c \in \varrho$, the lift $\widetilde{c}$ has its endpoint in $\widehat{\Lambda_n}$.
Our next step will be to estimate the cardinality of $\mathfrak{R}^f_n$.

We start by noticing that $\mathrm{pr}_2 (H^{-n}(\widetilde{\Lambda}_0)) = h^{-n}(\widehat{\eta})$. Since $\mathbf{\Psi}$ is isotopic to the identity we conclude that the curve $\widehat{\eta}_n = \mathrm{pr}_2(\mathbf{\Psi}(H^{-n}(\widetilde{\Lambda}_0)))$ is isotopic to $h^{-n}(\widehat{\eta})$ in the space of embedded curves in $S$.

\begin{remark} \label{remarkgrowth}
Recall from \cite[Expos\'e 1]{FLP} that for isotopy classes $\mathfrak{o}$ and $\mathfrak{q}$ of embedded closed curves in $S$ we define
\begin{equation}
\mathfrak{Int}(\mathfrak{o},\mathfrak{q})
\end{equation}
to be the minimal number of intersection points of curves $\mathfrak{o}$ and $\mathfrak{q}$.
Because $h$ has pseudo-Anosov monodromy and the curves  $\eta$ and $\widehat{\eta}$ are not homotopic to boundary
components of $S$, it follows from \cite[Expos\'e 1 and Expos\'e 11]{FLP} that there exists real numbers $\overline{a}>0$, $\overline{b}$ and $n_1$ such that
\begin{equation}
\mathfrak{Int}([\eta],[h^{-n}(\widehat{\eta})]) \geq e^{\overline{a}n+\overline{b}},
\end{equation}
for every $n\geq n_1$, where $[\eta]$ denotes the isotopy class of curves in $S$ isotopic to $\eta$, and $[h^{-n}(\widehat{\eta})]$ denotes the isotopy class of curves in $S$ isotopic to $[h^{-n}(\widehat{\eta})]$. Since $\eta_0$ is isotopic to $\eta$, and $\widehat{\eta}_n$ is isotopic to $h^{-n}(\widehat{\eta})$,
\begin{equation}
\mathfrak{Int}([\eta_0],[\widehat{\eta}_n]) = \mathfrak{Int}([\eta],[h^{-n}(\widehat{\eta})]).
\end{equation}
\end{remark}

We are now ready to prove the main result of this section.
\begin{thrm}
There exist numbers $\overline{a}>0$, $\overline{b}$ and $n_2$ such that
$\sharp(\mathfrak{R}^f_n) \geq  e^{\overline{a}n+\overline{b}}$ for every $n\geq n_2$.
\end{thrm}

\textit{Proof:} We endow $S$ with a hyperbolic metric $g$ having $\partial S$ as a geodesic boundary. Since $\eta_0$ is a  simple closed curve, \cite[ Lemma 2.6]{BC} implies that there is a homeomorphism $\psi: S \to S$ homotopic to the identity  for which $\psi(\eta_0)$ is a geodesic of the hyperbolic metric $g$.

As $\psi(\widehat{\eta}_n)$ is an embedded closed curve in $S$ it admits an isotopy to an embedded hyperbolic geodesic $\gamma$. The number of intersection points of $\gamma$ and $\psi(\eta_0)$ is equal to the number $\mathfrak{Int}([\psi(\eta_0)],[\psi(\widehat{\eta}_n)])=\mathfrak{Int}([\eta_0],[\widehat{\eta}_n])=\mathfrak{Int}([\eta],[h^{-n}(\widehat{\eta})])$, because pairs of hyperbolic geodesics realise the minimal number of intersection points in their isotopy classes.
  We let $\mathfrak{m}_n$ be the number of intersection points of $\gamma$ and $\psi(\eta_0)$, and  denote by $\{p^n_1,....,p^n_{\mathfrak{m}_n}\}$ the set of the intersection points of $\gamma$ and $\psi(\eta_0)$.

There exists a subset $\widetilde{S}$ of the Poincar\'e disc $(\mathbb{D},g_{-1})$ which is the universal cover of $(S,g)$. Let $\pi: \widetilde{S} \to S$ the covering map. Given an embedded closed curve $\mathfrak{q}$ in $S$ which is not homologous to $\psi(\eta_0)$, let $\overline{\mathfrak{q}}$ be a lift of $\mathfrak{q}$ in $\widetilde{S}$, and take a closed subinterval $I_\mathfrak{q}$ of $\overline{\mathfrak{q}}$ such that $\pi(\partial I)=p_0 \notin \psi(\eta_0)$ and that covers every point $x\neq p_0$ of $\gamma$ exactly once, in the sense that the intersection of $\pi^{-1}(x)$ and $I_\mathfrak{q}$ has one element. We call $I_\mathfrak{q}$ a fundamental interval of $\mathfrak{q}$.

From now on we suppose $n\geq n_1$ so that $\gamma$ and $\psi(\widehat{\eta}_n)$ are not homologous to $\psi(\eta_0)$. Consider a lift $\overline{\gamma}$ in $\mathbb{D}$, and take a fundamental interval $I_\gamma$ of $\overline{\gamma}$. It is well-known that two distinct geodesics of the hyperbolic plane can intersect in at most one point. It follows that $\overline{\gamma}$ cannot intersect a lift of $\psi(\eta_0)$ more than once. Therefore, $I_\gamma$ intersects exactly $\mathfrak{m}_n$ different lifts $\{\kappa_1,...,\kappa_{\mathfrak{m}_n}\}$ of $\psi(\eta_0)$.

Denote by $\gamma_t$ an isotopy between $\gamma$ and $\psi(\widehat{\eta}_n)$ where $t \in[0,1]$, and $\gamma_0 = \gamma$ and $\gamma_1 =\psi(\widehat{\eta}_n)$. Using the isotopy $\gamma_t$, we can construct a path $I_t$ of fundamental intervals of the curves $\gamma_t$.  This generates an isotopy of $I_\gamma$ to a fundamental interval $I_{\psi(\widehat{\eta}_n)}$ of $\psi(\widehat{\eta}_n)$ through fundamental intervals of $\gamma_t$. It follows from the properties of fundamental intervals that $\pi(\partial(I_t))$ is disjoint from $\psi(\eta_0)$ for all $t \in [0,1]$. We then conclude that $I_{\psi(\widehat{\eta}_n)}$ must also intersect the same $\mathfrak{m}_n$ different lifts $\{\kappa_1,...,\kappa_{\mathfrak{m}_n}\}$ of $\psi(\eta_0)$ that are intersected by $I_\gamma$, though it may intersect also others lifts of $\psi(\widehat{\eta}_n)$.

The set $\mathcal{A}$ of intersection points of $\eta_0$ and $\widehat{\eta}_n$ is in bijective correspondence with the set $\mathcal{L}$ of intersection points of $\psi(\eta_0)$ and $\psi(\widehat{\eta}_n)$. Because of the properties of fundamental intervals, there also exists a bijection between the set $\mathcal{L}$ and the set $\mathcal{B}$ of intersection points of $I_{\psi(\widehat{\eta}_n)}$ with the geodesics $\{\kappa_1,...,\kappa_{\mathfrak{m}_n}\}$.
There exists then a bijection $\varphi: \mathcal{A} \to \mathcal{B}$.
As observed  above, the set $\mathcal{A}$ is in bijective correspondence with the set of Reeb chords of $\widetilde{\mathfrak{p}}^*(X_{\widehat{\alpha}})$ from $\Lambda_0$ to $\widehat{\Lambda}_n$.

Taking now $p_1, p_2 \in \mathcal{A}$, we claim that there is a strip $v$ satisfying the four conditions of Proposition \ref{nielsen} if, and only if, $\varphi(p_1)$ and $\varphi(p_2)$ lie in the same $\kappa_j$.
To prove one direction of the claim notice that if there exists such a strip $v$ then we can take a lift $\overline{v}$ of $\psi\circ v$ in the universal cover $\widetilde{S}$. By the boundary conditions that are satisfied by $\overline{v}$ combined with the fact that $\psi(\eta_0)$ and $\psi(\widehat{\eta}_n)$ are embedded curves in $S$, it follows that $\varphi(p_1)$ and $\varphi(p_2)$ have to lie in the same $\kappa_j$.
For the reverse direction of the claim we remark that if $\varphi(p_1)$ and $\varphi(p_2)$ lie in the same $\kappa_j$ we can construct a strip $\overline{v}$ satisfying $\overline{v}([0,1] \times \{0\})= \varphi(p_1)$, $\overline{v}([0,1] \times \{1\})= \varphi(p_2)$, $\overline{v}(\{0\} \times [0,1]) \subset \kappa_j$ and $\overline{v}(\{1\} \times [0,1]) \subset \mathfrak{i}_n$, where $\mathfrak{i}_n$ denotes the lift of $\psi(\widehat{\eta}_n)$ that contains $I_{\psi(\widehat{\eta}_n)}$. Taking $v=\psi^{-1}\circ \pi \circ \overline{v}$ we obtain the desired strip satisfying the conditions of Proposition \ref{nielsen}.

Combining the last paragraph and Proposition \ref{nielsen} we conclude that to each $\kappa_j$ is associated a different relative Nielsen class $ \varrho_j$ in $\mathfrak{R}_n$, and that the intersection points between $I_{\psi(\widehat{\eta}_n)}$ and $\kappa_j$ are in bijective correspondence with the Reeb chords in $\varrho_j$. It follows that there are at least $\mathfrak{m}_n$ different relative Nielsen classes in $\mathfrak{R}_n$.

We must now prove that $ \varrho_j$ is a fundamental relative Nielsen class. We know that $I_{\gamma}$ intersects $\kappa_j$ exactly once, and we have an isotopy $I_t$ between $I_\gamma$ and $I_{\psi(\widehat{\eta}_n)}$ such that $\partial(I_t)$ never intersects $\kappa_j$. Since $I_{\gamma}$ and $I_{\psi(\widehat{\eta}_n)}$ are both transversal to $\kappa_j$, we conclude that $I_{\psi(\widehat{\eta}_n)}$ must intersect $\kappa_j$ an odd number of times. Combining this with the observation in the last paragraph we conclude that $ \varrho_j$ is a fundamental relative Nielsen class. Thus, there are at least $\mathfrak{m}_n$ different fundamental relative Nielsen classes in $\mathfrak{R}_n$.

Taking the numbers $\overline{a}>0$, $\overline{b}$ in Remark \ref{remarkgrowth}, and defining $n_2=\max\{n_0,n_1\}$ for the numbers $n_0$ in Remark \ref{above} and $n_1$ in Remark \ref{remarkgrowth}, it follows that if $n\geq n_2$
\begin{equation}
\mathfrak{m}_n \geq  e^{\overline{a}n+\overline{b}}.
\end{equation}
Since $\# \mathfrak{R}_n \geq \mathfrak{m}_n$, we have proved the theorem.
\qed

\

We have now all the ingredients needed to establish the exponential homotopical growth rate of $LC\mathbb{H}_{st}(\tau,\Lambda \to \widehat{\Lambda})$.
Let $N>0$ be a number such that $\mathbf{\Psi}(\Lambda_0)  \subset [-N,N] \times S$ and $\widetilde{\Lambda}_0 \subset [-N,N] \times S$. From the definition of the map $H$ it follows that for every integer $n$
\begin{equation} \label{incrusion}
H^{-n}(\widetilde{\Lambda}_0) \subset [-N + n,N+n] \times S.
\end{equation}
 Using formula \eqref{estimatepsi} we conclude that
 \begin{equation} \label{inclusion}
 \widehat{\Lambda}_n = \mathbf{\Psi}(H^{-n}(\widetilde{\Lambda}_0)) \subset (-\infty,\mathfrak{A}(N+n) + \mathfrak{B}] \times S,
 \end{equation}
where $\mathfrak{A}>1$ and $\mathfrak{B}>0$.

For $\varrho \in \mathcal{I}(\mathfrak{R}^f_n)$, we know that every Reeb chord in $\varrho$ lifts to a trajectory of the vector field $\partial_t$ in $\mathbb{R}\times S$ that starts in  $\mathbf{\Psi}(\Lambda_0)$ and ends in
$ \widehat{\Lambda}_n$. From \eqref{incrusion} and \eqref{inclusion} it follows that all trajectories of $\partial_t$ starting in $\mathbf{\Psi}(\Lambda_0)$ and ending in
$ \widehat{\Lambda}_n$ have length $\leq \mathfrak{A}(N+n) + \mathfrak{B} +N$. We thus conclude that all Reeb chords in $\varrho$ have action $ \leq \mathfrak{A}(N+n) + \mathfrak{B} + N$.
The homotopy class $\varrho$ contains an odd number of Reeb chords. The computation
\begin{eqnarray}
\dim LC\mathbb{H}^{\varrho}_{st}(\tau,\Lambda \to \widehat{\Lambda}) = \dim\ker(\partial_{st})- \dim Im(\partial_{st})= \\
= \dim LCH^{\varrho}_{st}(\tau,\Lambda \to \widehat{\Lambda}) - 2\dim Im(\partial_{st}))
\end{eqnarray}
implies that the numbers $\dim(LC\mathbb{H}^{\varrho}_{st}(\tau,\Lambda \to \widehat{\Lambda}))$  and
$\dim(LCH^{\varrho}_{st}(\tau,\Lambda \to \widehat{\Lambda}))$ have the same parity. Therefore $LC\mathbb{H}^{\varrho}_{st}(\tau,\Lambda \to \widehat{\Lambda}) \neq 0$. Our discussion in this paragraph shows that for the map $\mathcal{I}$ defined in \eqref{eqI} we have
\begin{equation}
\mathcal{I}(\mathfrak{R}^f_n) \subset \Sigma^{\mathfrak{A}(N+n) + \mathfrak{B} + N}_{\Lambda \to \widehat{\Lambda}}(\lambda_k),
\end{equation}
for $ \Sigma^{\mathfrak{A}(N+n) + \mathfrak{B}+N}_{\Lambda \to \widehat{\Lambda}}(\lambda_k)$ as defined in Section \ref{sectiongrowth}.

Combining the last two paragraphs we obtain:
\begin{equation}
\#(\Sigma^{C}_{\Lambda \to \widehat{\Lambda}}(\lambda_0)) \geq \#\mathcal{I}(\mathfrak{R}^f_{\lfloor \frac{C-\mathfrak{B} -N }{\mathfrak{A}}-N\rfloor}) \geq e^{\overline{a}\lfloor\frac{C-\mathfrak{B}-N}{\mathfrak{A}}-N\rfloor + \overline{b}},
\end{equation}
if $C$ is sufficiently large; here we are denoting by $\lfloor \mathrm{r} \rfloor$ the greatest integer $\leq \mathrm{r}$.
It is an elementary exercise to see that taking $a:=\frac{\overline{a}}{2\mathfrak{A}}$, there exist numbers $b$, and $C_0$ such that
\begin{equation}
e^{\overline{a}\lfloor\frac{C-\mathfrak{B}-N}{\mathfrak{A}}-N\rfloor + \overline{b}} \geq e^{aC+ b}
\end{equation}
for all $C\geq C_0$.
We have thus established the following theorem.
\begin{thrm} \label{theoremgrowth}
The strip Legendrian contact homology $LC\mathbb{H}_{st}( \lambda_k,\Lambda \to \widehat{\Lambda})$ has exponential homotopical growth rate with exponential weight $a= \frac{\overline{a}}{2\mathfrak{A}}$.
\end{thrm}

We now proceed to deduce Theorem \ref{theorem2} from Theorem \ref{theoremgrowth}. \\
\textit{Proof of Theorem 2:}
For each $k\in \mathbb{Z}$ we have proved that  $LC\mathbb{H}_{st}(\lambda_k,\Lambda \to \widehat{\Lambda})$ has exponential homotopical growth rate with exponential weight $a= \frac{\overline{a}}{2\mathfrak{A}}$. Using that the contact manifolds $(M,\xi_k)$ are all non-diffeomorphic by Proposition \ref{colinhonda}, we establish the theorem.
\qed

Combining Theorem \ref{theorem1} and Theorem \ref{theorem2}, we have that for every contact form $\lambda$ associated to $(M,\xi_k)$, the Reeb flow of $\lambda$ has positive topological entropy.

\section{Concluding remarks} \label{section5}

We mention that the methods developed in the present paper can be used to obtain estimates for the topological entropy in other families of contact manifolds. One other class of examples is obtained by applying the Foulon-Hasselblatt integral surgery (introduced in \cite{FH}) on the Legendrian lift of a separating geodesic of a hyperbolic surface; we refer the reader to \cite{A,FH} for the precise definition of the Foulon-Hasselblatt surgery. By the analysis in \cite{FH,HT}, we know that for most integral surgeries the resulting contact 3-manifold $(M',\xi')$ is not a Seifert fibre space, but an ``exotic'' graph manifold. In the author's Ph.D. thesis \cite{A1} the following theorem is proved:

\begin{thrm}
Let $(M, \ker \lambda_F)$ be the contact manifold associated with the hypertight contact form $\lambda_F$ obtained via the integral Foulon-Hasselblat surgery on the Legendrian $L_\sigma \subset T_1 S$, where $\sigma$ is a separating geodesic in the closed hyperbolic surface $S$ and $L_\sigma$ is its Legendrian lift. Then there exist disjoint Legendrian curves $\Lambda$ and $\Lambda'$ on $(M,\ker(\lambda_F))$ such that $\lambda_F$ is adapted to the pair $(\Lambda,\Lambda')$ and $LC\mathbb{H}_{st}(\lambda_F,\Lambda \to \Lambda')$ has exponential homotopical growth rate.
\end{thrm}

In a recent work \cite{Dahinden} Dahinden has obtained an extension of the results of \cite{MS} to positive contactomorphisms. He showed that if the homology of the based loop space of a manifold $Q$ is rich then every positive contactomorphism in $(T_1 Q, \xi_{geo})$ has positive topological entropy. We believe that a similar result should also hold for the contact manifolds covered by Theorem \ref{theorem2}.

In \cite{ReebAnosov} the author showed that if a contact 3-manifold $(M,\xi)$ admits a contact form $\lambda_0$ whose Reeb flow is Anosov, then every Reeb flow on $(M,\xi)$ has positive topological entropy. A question we believe is interesting is: does there exist a pair of Legendrian knots $\Lambda$ and $\widehat{\Lambda}$ in  $(M,\xi)$ such that, $\lambda_0$ is adapted to $\Lambda$ and $\widehat{\Lambda}$ and $LC\mathbb{H}_{st}(\lambda_0,\Lambda \to \widehat{\Lambda})$ has exponential growth?

It would be extremely interesting to know if there exist overtwisted contact 3-manifolds on which every Reeb flow has positive topological entropy. Unfortunately, there seems to be no available technology that could help to answer this question.

\begin{appendices}

\section{A Fubini type estimate } \label{appendixA}

\setcounter{lemma}{0}

In this appendix we prove a Fubini type estimate which is necessary in the proof of Theorem \ref{theorem1}. We begin by fixing some notation. Let $M$ be a closed oriented 3-manifold, $X$ be a smooth vector field without singularities on $M$, $\mathbb{L}$ be an embedded knot in $M$, and $g$ be a Riemannian metric on $M$.
We denote by $\phi^t_X$ the flow of the vector field $X$. For a vector $v$ tangent to $M$ we let $|v|_g$ be the norm of $v$ with respect to $g$.

The Riemannian metric $g$ provides us with a way to measure area for any surface immersed on $M$.
We fix a parametrisation\footnote{For this parametrisation we are using the identification $S^1 = \mathbb{R} /_{2\pi \mathbb{Z}}$.} $\mathbb{L}:S^1 \to M$ of the knot $\mathbb{L}$.
Define a map $F_{C,\mathbb{L}}: [0,C] \times S^1 \to M$ by the formula
\begin{equation}
F_{C,\mathbb{L}}(t,\theta) = \phi^t_X(\mathbb{L}(\theta)).
\end{equation}
We denote by $Area^C(\mathbb{L})$ the area with respect to $g$ of the immersed surface $F_{C,\mathbb{L}}([0,C] \times S^1)$.

We will denote by $\widehat{g}$ the pullback metric $(F_{C,\mathbb{L}})^*g$ of $g$ to $[0,C] \times S^1$ via the map $F_{C,\mathbb{L}}$.
 Associated to the Riemannian metric $\widehat{g}$ to $[0,C] \times S^1$ we have an area form $\mu_{\widehat{g}}$ on $[0,C] \times S^1$. Because $\widehat{g}$ is the pullback metric $(F_{C,\mathbb{L}})^*g$ we have the following equality:
\begin{equation} \label{eq74}
Area^C(\mathbb{L}) = \int_{[0,C] \times S^1} \mu_{\widehat{g}} = \int_{[0,C] \times S^1} |\mu_{\widehat{g}}(\partial_t, \partial_\theta)| dt d\theta.
\end{equation}

For a vector $v$ tangent to $[0,C] \times S^1$ we denote by $|v|_{\widehat{g}}$ its norm with respect to the metric $\widehat{g}$.
It is a classical fact that $|\mu_{\widehat{g}}(\partial_t, \partial_\theta)| \leq  |\partial_t |_{\widehat{g}}
 |\partial_\theta |_{\widehat{g}} $.
 It follows from the definitions of $F_{C,\mathbb{L}}$ and $\widehat{g}$ that
 \begin{eqnarray}
 |\partial_t |_{\widehat{g}} = |D(F_{C,\mathbb{L}})\partial_t|_g = | X|_g,    \label{eq75}           \\
  |\partial_\theta |_{\widehat{g}} = |D(F_{C,\mathbb{L}})\partial_\theta|_g, \label{eq76}
 \end{eqnarray}
 where $D(F_{C,\mathbb{L}})$ denotes the differential of $F_{C,\mathbb{L}}$.
Since $M$ is compact there exists a constant $\mathrm{K}>0$ such that $|X|_g \leq \mathrm{K}$ on all points of $M$. Combining this with \eqref{eq75} and \eqref{eq76} we obtain:
\begin{equation} \label{eq77}
\int_{[0,C] \times S^1} |\mu_{\widehat{g}}(\partial_t, \partial_\theta)| dt d\theta \leq  \int_{[0,C] \times S^1}  |\partial_t |_{\widehat{g}}  |\partial_\theta |_{\widehat{g}}  d\theta dt \leq \mathrm{K} \int_0^C \Big( \int_{S^1}  |\partial_\theta |_{\widehat{g}} d\theta \Big) dt.
\end{equation}

It is immediate to see that for each $t\in [0,C]$, the integral $\int_{S^1}  |\partial_\theta |_{\widehat{g}} d\theta$ equals the length $length(\phi^t_X(\mathbb{L}))$ of the curve $\phi^t_X(\mathbb{L})$ with respect to $g$. Combining this with  \eqref{eq74} and \eqref{eq77} we get:
\begin{equation} \label{eq78}
Area^C(\mathbb{L}) \leq \mathrm{K} \int_0^C length(\phi^t_X(\mathbb{L})) dt.
\end{equation}
For the statement of the next lemma we keep the notation introduced so far in the appendix.
\begin{lemma}
Suppose that there exist numbers $a>0$, $b$ and $C_0$ such that:
\begin{equation}
Area^C(\mathbb{L}) \geq e^{aC+b},
\end{equation}
for all $C\geq C_0$. Then
\begin{equation}
\limsup_{t \to +\infty} \frac{\log length(\phi^t_X(\mathbb{L}))}{t} \geq a.
\end{equation}
\end{lemma}
\textit{Proof:}
We argue by contradiction, assuming that  $\limsup_{t \to +\infty} \frac{length(\phi^t_X(\mathbb{L}))}{t} <a$.
Then there would be $t_0>0$ and $\epsilon<0$ such that:
\begin{equation}\label{eq81}
length(\phi^t_X(\mathbb{L})) < e^{(a-\epsilon)t}
\end{equation}
for all $t\geq t_0$. Integrating both sides of the \eqref{eq81} for $t$ between $0$ and $C$ and invoking \eqref{eq78} , we would conclude that:
\begin{equation} \label{eq82}
Area^C(\mathbb{L}) < \frac{e^{C(a-\epsilon)}-e^{t_0(a-\epsilon)}}{a-\epsilon} + \mathrm{K}\int_{0}^{t_0}  length(\phi^t_X(\mathbb{L}))dt
\end{equation}
for all $C\geq t_0$. It is immediate that the right side of \eqref{eq82} becomes a lot smaller than $e^{C(a-\frac{\epsilon}{2})}$ for sufficiently large $C$. This would imply that  $Area^C(\mathbb{L}) \leq e^{C(a-\frac{\epsilon}{2})}$ for large $C$ and would lead to a  contradiction, since we assumed $Area^C(\mathbb{L}) \geq e^{aC+b}$ for $C$ sufficiently large. \qed

\section{ Proof of Lemma \ref{lemma:changeofposition} } \label{AppendixB}

In this appendix we present, for the convenience of the reader, a proof of the following lemma.
\setcounter{lemma}{2}
\begin{lemma} \label{lemma:changeofposition'}
Let $\lambda$ be a contact form on $(Y,\xi)$ and $\Lambda$ and $\widehat{\Lambda}$ be disjoint connected Legendrian submanifolds in $(Y,\xi)$. Let $\mathcal{V}(\widehat{\Lambda})$ be a tubular neighbourhood of $\widehat{\Lambda}$ that does not intersect $\Lambda$. Then, given $\epsilon>0$ there exists $\delta>0$ such that for every $\widehat{\Lambda}'$ that is $\delta$-close to $\widehat{\Lambda}$ in the $C^3$-sense there exists a contactomorphism $\Upsilon_{\widehat{\Lambda}'}: (Y,\xi) \to (Y,\xi)$ which satisfies
\begin{itemize}
\item [(1)] $\Upsilon_{\widehat{\Lambda}'}(\widehat{\Lambda}) = \widehat{\Lambda}'$,
\item [(2)] $\Upsilon_{\widehat{\Lambda}'}$ is $\epsilon$-close to the identity in the $C^1$-sense,
\item [(3)] the support of is contained in $\mathcal{V}(\widehat{\Lambda})$.
\end{itemize}
\end{lemma}

We start by recalling some necessary notions. We follow the conventions of \cite[Section 2.3]{Ge}. Given a contact manifold $(Y,\xi)$ and a contact form $\lambda$ on $(Y,\xi)$, one associates to each function $H: [0,1]\times Y \to  \mathbb{R}$ a \textit{contact Hamiltonian vector field} $X_{\lambda,H}$ characterised by
\begin{eqnarray}
\lambda(X_{\lambda,H}) = H, \\
d\lambda(X_{\lambda,H}, \cdot) = dH(X_{\lambda})\lambda - dH.
\end{eqnarray}
 It follows directly from this definition that the $C^1$-norm of $X_{\lambda,H}$ depends continuously on the $C^2$-norm of $H$. Therefore the $C^1$-distance of the time one map $\Upsilon_{H,\lambda}$ of the flow of $X_{\lambda,H}$ to the identity map of $Y$ also depends continuously on the $C^2$-norm of $H$.
We are now ready to prove the lemma.
 
 \

\textit{Proof of Lemma \ref{lemma:changeofposition}:}

\textbf{Step 1:} We first take $\mathcal{U}(\widehat{\Lambda})$ to be a tubular neighbourhood of $\widehat{\Lambda}$ whose closure is contained in  $\mathcal{V}(\widehat{\Lambda})$. It is clear that there exists $\delta_0>0$ such that:
\begin{itemize}
\item [(A)] every $\widehat{\Lambda}'$ that is $\delta_0$-close to $\widehat{\Lambda}$ in the $C^3$ sense is contained in $\mathcal{U}(\widehat{\Lambda})$. 
\end{itemize}

\textbf{Step 2:}
Let $\mathcal{L}_{\widehat{\Lambda}'}: [0,1] \times \widehat{\Lambda} \to Y$ be the Legendrian isotopy between $\widehat{\Lambda}$ and $\widehat{\Lambda}'$ contained in $\mathcal{U}(\widehat{\Lambda})$ . 
We will associate to  $\mathcal{L}_{\widehat{\Lambda}'}$ a function $H_{\mathcal{L}_{\widehat{\Lambda}'}}:[0,1]\times Y \to \mathbb{R}$. First we define it over the trace of $\mathcal{L}_{\widehat{\Lambda}'}$ via the formula 
\begin{equation} \label{eqham}
H_{\mathcal{L}_{\widehat{\Lambda}'}}(t,\mathcal{L}(t,p))= \lambda(\partial_t \mathcal{L}_{\widehat{\Lambda}'}(t,p)).
\end{equation}
 For each $\mu>0$, it is elementary to see that there exists $0<\delta_1 \leq \delta_0$ such that  if $\mathcal{L}_{\widehat{\Lambda}'}$ is $\delta_1$-small in the $C^3$-sense then (A) holds and the extension of $H_{\mathcal{L}_{\widehat{\Lambda}'}}$ can be chosen to satisfy:
\begin{itemize}
\item [(B)] $H_{\mathcal{L}_{\widehat{\Lambda}'}}$ is $\mu$-small in the $C^2$ norm,
\item [(C)] the support of $H_{\mathcal{L}_{\widehat{\Lambda}'}}$ is contained in $\mathcal{V}(\widehat{\Lambda})$.
\end{itemize}

\textbf{Step 3:}Associated to $\lambda$ and the function $H_{\mathcal{L}_{\widehat{\Lambda}'}}$ we have the contact Hamiltonian vector field $X(\lambda,H_{\mathcal{L}_{\widehat{\Lambda}'}})$ as defined above.  We denote by $\Upsilon_{\widehat{\Lambda}'}$ to be the time-one map of $X(\lambda,H_{\mathcal{L}_{\widehat{\Lambda}'}})$.

As remarked above the $C^1$ distance of the time-one map of the flow of $X(\lambda,H_{\mathcal{L}_{\widehat{\Lambda}'}})$ to the identity depends continuously on the $C^2$ norm of  $H_{\mathcal{L}_{\widehat{\Lambda}'}}$.
Therefore, given $\epsilon>0$ we fix $\mu>0$ such that if $H_{\mathcal{L}_{\widehat{\Lambda}'}}$ is $\mu$-small in the $C^2$ norm then $\Upsilon_{\widehat{\Lambda}'}$ is $\epsilon$-close to the identity in the $C^1$-sense. 
We now choose $\delta>0$ so that (A), (B) and (C) hold for this choice of $\mu$. It then follows that $\Upsilon_{\widehat{\Lambda}'}$ is $\epsilon$-close to the identity in the $C^1$-sense which implies (2).

It is shown in proof of \cite[Theorem 2.6.2]{Ge} that \eqref{eqham} implies that   $\Upsilon_{\widehat{\Lambda}'}$ satisfies (1). That $\Upsilon_{\widehat{\Lambda}'}$ satisfies (3) follows from the fact that $H_{\mathcal{L}}$ satisfies~(B).
\qed

\end{appendices}

\end{document}